\documentclass[a4paper,12pt]{amsart}
\usepackage{amsmath,amsthm,amssymb, mathdots}
\usepackage{amstext, amscd, setspace,mathtools, enumerate, graphics, latexsym, siunitx}
\usepackage{xcolor}
\usepackage[T1]{fontenc}
\usepackage[all]{xy}
\usepackage{mathrsfs}

\usepackage{bm}
\usepackage[hidelinks]{hyperref}
\usepackage{tikz}

\usetikzlibrary{positioning}
\usetikzlibrary{arrows.meta}

\usepackage[margin=1in]{geometry}

\makeatletter
\renewenvironment{proof}[1][{\bf \proofname}]{\par
  \pushQED{\qed}%
  \normalfont \topsep6\p@\@plus6\p@\relax
  \trivlist
  \itemindent\z@ % original has \normalparindent
  \item[\hskip\labelsep
    #1\@addpunct{.}]\ignorespaces
}{%
  \popQED\endtrivlist\@endpefalse
}
\makeatother

\newtheorem{theorem}{Theorem}[section]
\newtheorem{corollary}{Corollary}[theorem]
\newtheorem{proposition}[theorem]{Proposition}
\newtheorem{lemma}[theorem]{Lemma}
\newtheorem{remark}[theorem]{Remark}

\theoremstyle{definition}

\theoremstyle{definition}
\newtheorem{example}{Example}[section]

\title{On the Linear Algebraic Monoids Associated to Congruence of Matrices}
\author{Himadri Mukherjee and Gunja Sachdeva}
\address{BITS Pilani, K.K. BIRLA GOA Campus, NH 17B, Bypass Road, Zuarinagar, Sancoale, Goa 403726, India.}
\email{himadrim@goa.bits-pilani.ac.in}
\address{BITS Pilani, K.K. BIRLA GOA Campus, NH 17B, Bypass Road, Zuarinagar, Sancoale, Goa 403726, India.}
\email{gunjas@goa.bits-pilani.ac.in}

\begin{document}

\maketitle

\begin{abstract}

This paper discusses the generalized congruence equation $X^tAX=B$, for $X \in M_n(k)$ over any field $k$, through the action of monoid $Sol_A \times Sol_B :=  \{X \ | \ X^tAX = A\} \times  \{X \ | \ X^tBX = B\}$. We have completely characterized for what matrices $A$, the monoid $Sol_A$ is a Lie group. We have given the structure of the Lie group $Sol_A$ and $Sol_{A^2}$, and their Lie algebras when $A$ is $n \times n$ nilpotent matrix of nilpotency $n$. In this case, we have also proved that the invariants of $Sol_A$ for any $n$, and $Sol_{A^2}$ for $n$ even, are finitely generated.
\end{abstract}

\section{Introduction}
For a field $k$, and matrices $A$, $B$ in $M_n(k)$, we are interested in finding $Sol(A, B) = \{X \in M_n(k) \ | \ X^t A X = B \}$ such that $Sol(A,B) \cap GL_n(k)$ is non-empty. This problem is called the congruence problem of matrices, which is a part of a more general class of problems called the simultaneous equivalence or similarity of matrices \cite{Simultaneous}. This is an important problem in linear algebra and representation theory \cite{wild}, stated as follows:\\

   {\it  For a collection $ \mathcal{X} = \{(A_1, \ldots, A_d) \ | \ A_i \in M_n(k)\}$, consider the action of $GL_n(k)$ (or $GL_n(k) \times GL_n(k))$ on $\mathcal{X}$ by $g \ast (A_1, \ldots, A_d) = (gA_1 g^{-1}, \ldots, g A_d g^{-1})$ (or by $(g, h) \ast (A_1, \ldots, A_d) = (gA_1h, \ldots, g A_dh)$). The problem is to classify the orbits of this action.}\\
   
For an algebraically closed field $k$, and $d=2$, the simultaneous equivalence problem is called the classification of the matrix pencils $P(A,B)=\lambda A + \mu B$, up to similarity. The matrix pencil classification problem has been fully solved by the cumulative efforts of several authors spanning the last century. Beginning with, the classification of singular matrix pencils by Kronecker \cite{Aitken}. The classification of a special class of matrix pencils, known as the ``palindromic pencils" $\lambda A + \mu A^t$, answers the matrix congruence problem through the canonical forms \cite{Teran_history}. These canonical forms for the palindromic pencils are not easy to work with as these are given as a direct sum of six types of matrices (see \cite{Aitken}). This resulted in several articles with further simplification approaches to this problem. Currently, there are at least four distinct types of canonical forms of equivalent pencils; each of these attempts at either reducing the number of canonical forms or simplifying the structure of the canonical forms. All of these are of a similar essence in the sense that they are given in terms of direct sums of a set of matrices.  \\

Recently, F. D. Teran has systematically studied the linear equation $X^tA + AX = 0$, related to the congruence and the conjugate transpose problem in \cite{Teran_Lie}, and the general equation related to the star conjugate in \cite{Teran_Lie_star}, and \cite{Teran_Lie_Star_two}. He has used solutions of these linear equations to understand various invariants of the orbits, such as the dimension, etc (see \cite{Teran_orbits}). One can then refer an excellent paper by Teran \cite{Teran_history} for the historical perspective on this topic. The area of matrix pencils also has practical applications in PDE and perturbations theory. R. Bapat et. al. in \cite{bapat_pencil} and in \cite{bapat_pencil_2} have studied perturbations of matrix pencils and theories regarding matrix polynomials in general. R. Bhatia in \cite{bhatia_pencil} has studied perturbations of matrix pencils with real spectrum. The field of perturbations of matrix pencils with constraints has given rise to many recent works for their utility in various areas of applied mathematics (see \cite{Safique}).\\

 Towards the congruence matrix equations, in a series of papers \cite{W1}, \cite{W2}, Williams has studied the congruence problem for finite fields, $\mathbb{R}$ and $\mathbb{C}$. Moreover, Horn, Sergeichuk, and others have results on canonical forms for congruence and $^*$congruence matrix problems \cite{W4}, or the congruence of two-dimensional subspaces \cite{W3}. Special cases of the congruence question, especially the question of *-congruence \cite{Horn1} or unitary congruence \cite{Horn2} have also attracted a lot of attention from the research community. The projective congruence question ($X^tAX=\lambda B$ for any non-zero $\lambda \in k$) has also garnered interest (see \cite{Projective}), with recent work by Edelman et al. in \cite{Complex_case} focusing on the geometric structure of solutions to the matrix equation $X^tAX=A$ and $X^*AX=A$ for complex matrices. 
 The theory of linear algebraic monoids is an interesting and emerging field of research \cite{Putcha}, \cite{Renner}. Such monoids arising out of matrix equations are also studied in other articles like \cite{Complex_case}. In the paper \cite{integer}, the authors calculated the number of solutions of a general congruence problem, and also calculated the number of solutions coming from a finite ring.\\
 
 If there exists an inner automorphism of $GL_n(k)$ that takes $Sol_A$ to $Sol_B$, then one can observe that $A$, $B$ are congruent to each other. Hence the multiplicative structure of these monoids are important to be understood. We classify the matrices $A$ such that $Sol_A \cap GL_n = Sol_A$, which allows us to look at the action of $Sol(A) \times Sol(B)$ on $M_n(k)$ by $(g,h)\cdot X = gXh$. Let the orbit decomposition be, $Sol(A,B) = \displaystyle\bigcup_{X \in Sol(A, B)} \mathcal{O}_X$, where $\mathcal{O}_X = \{gXh \ | \ (g, h) \in Sol(A) \times Sol(B)\}.$ Using Orbit-Stablizer theorem, $\mathcal{O}_X = Sol(A) \times Sol(B)/ G_X.$ Now, either one finds the stabilizer $G_X$ directly through some clever linear algebra or employs classical invariant theory to find the orbits. Once the question of finite generation of the invariants of $Sol_A$ is understood, the quotient space leads to another form of geometric invariants that solves the matrix congruence problem.  That is, we will fully understand the structure of $Sol(A, B)$ both algebraically and geometrically.\\
 
Unfortunately, these groups, in general, may not be reductive. So, the question of the invariants poses a difficult problem of finding quotients by non-reductive groups. In this regard, there are a number of recent developments dealing with the question of invariant theory of non-reductive groups, see \cite{observable}, \cite{Kirwan}, and \cite{non-reductive}. The idea is to reduce the problem of the invariants of a non-reductive group to that of an unipotent group. \\

For a matrix $A$ such that $Sol_A$ is a group, we look at the Lie algebra $sol_A = \{X \ | \ A^t X + XA = 0\}$. This gives us an exact sequence $0 \rightarrow \mathfrak{a} \rightarrow sol_A \rightarrow sol_A / \mathfrak{a} \rightarrow 0$, where $\mathfrak{a}$ is the unipotent radical. Let $N$ be the normal subgroup with Lie algebra isomorphic to $\mathfrak{a}$. Once this sequence is ``exponentiated", it leads to an exact sequence of the corresponding groups. This is the key step to identifying the equations of the group. Indeed, it gives us $Sol_A \simeq N \rtimes G/N$ where $N$ is such that $Lie(N) = \mathfrak{a}$. As a result, we have simplified the question of understanding the equations of $Sol_{A}$ to that of understanding the equations of the unipotent radical and that of the reductive quotient. \\
\\
\section{Notations and terminologies}
\begin{tabular}{p{10em} p{50em}} 
$k$ & algebraically closed field\\
$A, B$ & any $n \times n$ matrices with entries in $k$\\
$O_n$ & $n \times n$ zero matrix\\
$A \sim B$ &  there exist an invertible matrix $P$ such that $P^{-1}AP = B$\\
$A \simeq B$ & there exists an invertible matrix $P$ such that $P^tAP=B$\\
$\text{Sol}(A, B)$ & $\{X \in M_n(k) \ | \ X^tAX = B\}$\\
$\text{Sol}_A$ & $\{X \in M_n(k) \ | \ X^tAX = A\}$\\
$\text{Sol}_A^\ast$ & $\{X \in M_n(k) \ | \ X^\ast AX=A\}$\\ 
$\text{Sol}^\ast(A,B)$ & $\{X \in M_n(k) \ | \ X^\ast AX=B\}$\\
$G_A$ & $GL_n(k) \cap \text{Sol}_A$\\
$A_n$ & $n \times n$ upper triangular nilpotent matrix with nilpotency $n$\\
$A_d$ & $n \times n$ diagonal matrix\\
$A^2_n$ & $n \times n$ upper triangular nilpotent matrix with nilpotency $n-1$\\
$Ker(A)$ & $\{v \in k^n \ | \ Av = 0\}$\\
$\text{sol}_A$ & Lie algebra of $\text{Sol}_A = \{X \in M_n(k) \ | \ X^tA +AX = O_n\}$\\
$[X,Y]$ & $XY -YX$ for $X, Y \in \text{sol}_A$\\
$G^\circ$ & connected component of the group at identity
\end{tabular}
\vskip1 cm
This manuscript is organized in the following way: in $\S$ \ref{Sec3}, we classify all matrices $A$ such that the monoid $Sol_{A}$ coincides with its group of invertible elements. In $\S$ \ref{Sec4}, we look at the types $\mathcal{A}_{2k+1}$, $\mathcal{B}_{2k+1}$ (i.e. for $A_n$), and work towards finding the structure of these groups. The $\S$ \ref{Sec5} contains the structure of $Sol_{A_n^2}$; this is done by first understanding the Lie algebra in the sub-sections; we then use it to find the structure of the groups for $A$, where $A$ is congruent to a direct sum of two (possibly same) from the classes $\mathcal{A}_{2k+1}, \mathcal{B}_{2k}$. In $\S$ \ref{Sec6}, we discuss the finite generation of the invariants and examples. At the end, in $\S$ \ref{Sec7}, we have an appendix where we prove an algebraic identity in a non-commutative involutive polynomial algebra.

\section{The Congruence Monoid}\label{Sec3}

 In the relation $X^tAX=B$, we call $A$, to be the coefficient matrix. For any matrix $A$, we know that $Sol_{A}$ is a linear algebraic monoid, which acts on the set $Sol(A, B)$, whenever this set is non-empty. We will exploit this monoid structure and the related action on the set $Sol(A, B)$. Furthermore, define the group of the monoid $Sol_{A}$ as $G_A:= GL_n \cap Sol_A$, which we will henceforth call the \textit{congruence group} of the matrix $A$. The main goal of this section is to obtain a complete characterization of the coefficient matrices $A$ for which $G_A = Sol_A$. Recall an action of the monoid $Sol_A \times Sol_B$ on the set $Sol(A, B)$ by $(g,h)\cdot X = gXh$, where the right-hand side product is the product of matrices. 
  This action can be thought of as a restriction of the action of $Sol_A \times Sol_B$ on the matrices $M_n(k)$, and further, the action also restricts to an action by the group $G_A \times G_B$ on $M_n(k)$. It is important to understand the structure of $Sol_A$ because of the following observation: if $A=P^tBP$ with $P$ invertible and $X \in Sol_A$, then $PXP^{-1} \in Sol_B$, which clearly defines an isomorphism of monoids. In other words, $A\simeq B$ if and only if $Sol_{A}$ is isomorphic to $Sol_B$ as monoids by an inner-automorphism of $GL_n$. By Skolem-Noether theorem \cite{Skolem_Noether}, one can improve the statement by saying that $A$, $B$ are congruent if and only if $Sol_A$ and $Sol_B$ are isomorphic by an algebra automorphism of $M_n(k)$. Also note that, the monoids $Sol_{A}$ and $Sol_{A^t}$ are naturally isomorphic (equal as sets), and so are the monoids $Sol_A^*$ and $Sol_{A^*}^*$.\\

To obtain a necessary and sufficient condition on \textit{coefficient matrix} $A$, for the monoid $Sol_{A}$ to be a group, we first prove the following lemma:

\begin{lemma}\label{primary}
 Let $V$ be a finite dimensional vector space over $\mathbb{C}$, and $N$ be the matrix ${\footnotesize\begin{pmatrix}
        1 & i \\
        i & -1
    \end{pmatrix}}.$ The following are true:
   \begin{enumerate} \item \label{group_one} If $v \in V$ such that $v^tv=0$ and $B_v=vv^t \neq 0$, then $B_v \simeq {\footnotesize\begin{pmatrix}
        N & 0\\
        0 & 0
    \end{pmatrix}}$, where 0 represents blocks of zero with appropriate size.
    \item \label{group_two} Moreover, if $v \in V$  such that $vv^t$ is nilpotent, then either $v =0$ or $vv^t \simeq {\footnotesize\begin{pmatrix}
        N & 0\\
        0 & 0
    \end{pmatrix}}$.
    \item \label{group_three}
If $A$ and $B$ are matrices of dimension $n\times 2$ and $2 \times m$ respectively, and let $X$ be the matrix ${\footnotesize\begin{pmatrix}
    0 & 0\\
    -i & 1
\end{pmatrix}}$. Then the following are true:
\begin{itemize}
    \item If $AN=0$, then $AX=A.$
    \item If $NB=0$, then $X^tB=B$. 
    \item If $A= {\footnotesize\begin{pmatrix}
        N & A_2\\
        A_3 & A_4
    \end{pmatrix}}$ with $NA_2=0$ and $A_3N=0$, then the matrix ${\footnotesize\begin{pmatrix}
        X & 0\\
        0 & I
    \end{pmatrix}}$ is an element of $Sol_{A}$.
\end{itemize}
\end{enumerate}
   
\end{lemma}
\begin{proof}

\begin{itemize}
    \item To prove (\ref{group_one}), first, we observe that it is enough to consider $V$ to be two-dimensional. Let us assume that the vector $v = xe_1+ye_2$, where $e_1,e_2$ are standard basis vectors. Since $v^tv=0$, we have $x^2+y^2=0$, and the condition $vv^t \neq 0$ implies at least one of $x,y \neq 0$. But we also have $x^2+y^2=0$ so we conclude both $x,y \neq 0$. Now let us take $y=ix$ (the proof is the same for $y=-ix$), then we can write the matrix $B_v$ in this basis as, ${\footnotesize\begin{pmatrix}
        x^2 & ix^2 \\
        ix^2 & -x^2
    \end{pmatrix}}.$ This matrix is clearly congruent to ${\footnotesize\begin{pmatrix}
        1 & i \\
        i & -1
    \end{pmatrix}}.$ 
    \item For (\ref{group_two}), if $v \neq 0$ then observe that $B_v^2=0$ since rank of $B_v$ is one and it is nilpotent. Now $B_v=0$ clearly gives $v=0$. So we must assume $B_v \neq 0$, then $v^tv=0$ and thus by the Lemma (\ref{group_one}) we have the result.
    \item For the first statement note that $AN=0$ implies that $C_1=iC_2$, where $C_i$ are columns of $A$. Then $AX= (i C_2 C_2)=(C_1 C_2)=A$. For the second statement, note that $NB=0$ implies $R_1+iR_2=0$, where $R_i$ are rows of $B$. Then $X^tB= {\footnotesize\begin{pmatrix}-iR_2\\ R_2 \end{pmatrix}} = {\footnotesize\begin{pmatrix} R_1\\ R_2 \end{pmatrix}}=B$. For the third statement, consider $X^tAX= {\footnotesize\begin{pmatrix}
        X^t & 0\\
        0 & I \end{pmatrix}} {\footnotesize\begin{pmatrix}
        N & A_2\\
        A_3 & A_4
    \end{pmatrix}} {\footnotesize\begin{pmatrix}
        X & 0\\
        0 & I \end{pmatrix}} = {\footnotesize\begin{pmatrix}
        X^t & 0\\
        0 & I \end{pmatrix}} {\footnotesize\begin{pmatrix}
        NX & A_2\\
        A_3X & A_4
    \end{pmatrix}}$. But since $A_3N=0$, we have $A_3X=A_3$ and $NX=N$. We can rewrite $X^tAX = {\footnotesize\begin{pmatrix}
        X^t & 0\\
        0 & I \end{pmatrix}} {\footnotesize\begin{pmatrix}
        N & A_2\\
        A_3 & A_4
    \end{pmatrix}} = {\footnotesize\begin{pmatrix}
        X^tN & X^tA_2\\
        A_3 & A_4
    \end{pmatrix}}$. Now $X^tN=N$ and since $NA_2=0$, we have $X^tA_2=A_2$. This gives $X^tAX= {\footnotesize\begin{pmatrix}
        X^tN & X^tA_2\\
        A_3 & A_4
    \end{pmatrix}} = {\footnotesize\begin{pmatrix}
        N & A_2\\
        A_3 & A_4
    \end{pmatrix}} =A$, showing that $X$ is a solution. This completes the proof of (\ref{group_three}).
\end{itemize}

\end{proof}
\begin{theorem}\label{group_main}
    If $A$ is a matrix over any field $k$ such that $Sol_{A}$ is a group, then $vv^t$ is nilpotent for all vectors $v \in ker(A) \cap ker(A^t)$.
\end{theorem}
\begin{proof}
  Let us assume that $Sol_{A}$ is a group, and there exists a vector $v$ such that $Av=0$, $v^tA=0$, and $vv^t$ is not nilpotent. From these conditions we will construct a solution $X$ which is not invertible, thus arriving at a contradiction. Let us choose a parameter $c \in k$ and observe the matrix $X=I- c vv^t$, as we have assumed $vv^t \neq 0$, this matrix is a solution as $(I-cvv^t)A(I-cvv^t)= (A-0)(I-cvv^t)=A$. Now we can choose a $c$ such that this matrix $I-cvv^t$ is not invertible, namely $c$ as the reciprocal of a non-zero eigenvalue of $vv^t$. This completes the proof.
\end{proof}

In the following example, we show that the condition in the above theorem is not sufficient for $Sol_{A}$ to be a group. 
\begin{example}
    If $A$ is the complex matrix 
${\footnotesize\begin{pmatrix}
    1 & i \\
    i & -1
\end{pmatrix}}$, then the only vector (up to a scalar multiple)  in the set $ker(A) \cap ker(A^t)$ is $v={\footnotesize\begin{pmatrix}
    1 \\
    i
\end{pmatrix}}$, and one can check that $vv^t$ is nilpotent, but $Sol_{A}$ contains the matrix ${\footnotesize\begin{pmatrix}
    0 & 0 \\
    -i & 1
\end{pmatrix}}$.
\end{example}

\begin{theorem}\label{group}
    For a complex matrix $A$, $Sol_{A}$ is a group if and only if the intersection $ker(A)\cap ker(A^t)$ is trivial.
\end{theorem}
\begin{proof}
 
    If the intersection contains only the zero vector, let us assume that $Sol_{A}$ is not a group. Then, there exists $X \in Sol_A$, which is noninvertible. That is, there is a non-zero vector $v$ such that $Xv=0$, which implies $X^tAXv=Av=0$. Similarly, one sees that $v^tA=0$ together, we obtained $v \in ker(A) \cap ker(A^t)$, which is a contradiction. Conversely, assume that $Sol_{A}$ is a group, then from the Theorem \ref{group_main}, $ker(A) \cap ker(A^t)$ contains vectors $v$ such that $B_v=vv^t$ is nilpotent. Now, by the (\ref{group_two}) of Lemma \ref{primary}, we know that either $v=0$ or $B_v \simeq {\footnotesize\begin{pmatrix}
        N & 0\\
        0 & 0
    \end{pmatrix}}$. If $v=0$, then we are done. If not, then there exist invertible matrix $P$ such that $B_v =P {\footnotesize\begin{pmatrix}
        N & 0\\
        0 & 0
    \end{pmatrix}}P^t$. Now put $A' = P^tAP$ and $v' = P^{-1}v$, then $v'v'^t=B_{v'} = {\footnotesize\begin{pmatrix}
        N & 0\\
        0 & 0
    \end{pmatrix}}$. Furthermore, observe that $v' \in \ker(A') \cap ker( A'^t)$ and $Sol_{A'}$ isomorphic to $ Sol_A$. Therefore, by abuse of notation, we will call $A'$ as $A$ and $B_{v'}$ as $B_v$. So, if we prove that $Sol_{A}$ is not a group, then we are done. Note that $AB_v=B_vA=0$, writing $A$ in a block form as below \[{\footnotesize \begin{pmatrix}
        A_1 & A_2\\
        A_3 & A_4
    \end{pmatrix}},\]
    we obtain $A_1N=NA_1=0$, $NA_2=0$, and $A_3N=0$. Then $A_1N=NA_1=0$ implies that $A_1$ is a scalar multiple of $N$; we may assume that it is a nonzero scalar multiple (we deal with the $A_1$=0 case later). So once we multiply $A$ with an appropriate scalar, we can assume that $A_1=N$. Thus we have $A$ in the same form as in the (\ref{group_three}) of the Lemma \ref{primary}, so we have an element in $Sol_{A}$, \[Y={\footnotesize\begin{pmatrix}
        X & 0 \\
        0 & I
    \end{pmatrix}}\] where $X$ is ${\footnotesize\begin{pmatrix}
        0 & 0 \\
        -i & 1
    \end{pmatrix}}$. Clearly, $Y$ is not invertible. Thus, $Sol_{A}$ is not a group in this case. If $A_1=0$, then we could have taken \[Y={\footnotesize\begin{pmatrix}
        0 & 0 \\
        0 & I
    \end{pmatrix}}\] and see that it is a solution.  This completes the proof of the theorem.
\end{proof}
Let us call a matrix pencil $\lambda A + \mu B$ to be completely singular if there is a nonzero vector $v$ (independent of $\lambda, \mu$) such that $(\lambda A + \mu B )v =0$ for all $\lambda$, and $\mu$. Note that it is a stronger condition than the matrix pencil $\lambda A= \mu B$ just being singular. With this, we can write an equivalent statement for the above theorem. Since this statement is about the matrix pencils, which are notably related to the question of congruence \cite{Aitken}, it is an interesting result.

\begin{theorem}\label{pencil_version}
Let $A$ be a matrix in $M_n(\mathbf{C})$. Then $Sol_{A}$ is a group if and only if the palindromic pencil $\lambda A + \mu A^t$ is not completely singular.
\end{theorem}
\begin{proof}
The proof is immediate, as $ker(A) \cap ker(A^t)$ is trivial if and only if the palindromic pencil $\lambda A + \mu A^t$ is not completely singular.
\end{proof}

Similarly, one can prove results for the monoids arising out of the equation $X^*AX=A$, for a complex coefficient matrix $A$. We will characterize the case when $Sol_A^*$ is a group. Furthermore, we will deal with the structure of these groups in special cases as a part of another paper\cite{Next}. 
\begin{theorem}\label{complex}
    $Sol_A^*$ is a group if and only if $ker(A) \cap ker(A^*)$ is trivial.
\end{theorem}
\begin{proof}
    First, let us prove that if $ker(A)\cap ker(A^*)$ is trivial, then $Sol_A^*$ is a group. If not, then there exists $X \in Sol_A^*$, which is not invertible. Let $v$ be a non-zero vector such that $Xv=0$. As a result, we see that $X^*AXv=Av=0$ and $v^*XAX=v^*A=0$, or $v \in ker(A) \cap ker(A^*)$, a contradiction. Now for the other direction, let us assume that $Sol_A^*$ is a group and $v \in ker(A) \cap ker(A^*)$ be a non-trivial vector. Then $B_v= vv^*$ is a non-nilpotent matrix because $B_v^2 \neq 0$. Let $\lambda$ be a non-zero eigenvalue of $B_v$,  and let $t$ be a real number. Consider the matrix $X_t= tB_v - I_n$, and observe that it belongs to $Sol_A^*$ for any $t$. Let $f(t)$ be the determinant of the matrix $X_t$; it is clear that $f(t) =t^n \phi(1/t)$ where $\phi$ is the characteristic polynomial of $B_v$. Now we can take $t=1/\lambda$ and observe that $X_t$ is not invertible for this particular choice of $t$. This contradicts the assumption that $Sol_A^*$ is a group.
 \end{proof}
 Similar to the Theorem \ref{pencil_version}, one can write a pencil version of the above theorem.
 \begin{theorem}\label{pencil_version_two}
For $A \in M_n(\mathbb{C})$, $Sol^*_A$ is a group if and only if the palindromic pencil $\lambda A + \mu A^*$ is not completely singular.
\end{theorem}
 
 We can immediately observe from the above theorems that for a symmetric/skew-symmetric matrix $A$, $Sol_{A}$ is a group if and only if it is invertible, and for a Harmitian/skew Hermitian matrix $A$, $Sol_A^*$ is a group if and only if $A$ is invertible. One also observes that for any complex number $z$ of absolute value one and for any $X \in Sol_A^*$, we have $zX \in Sol_A^*$. As a result, we have a group action of the circle group $S^1$ on the monoid $Sol_A^*$. In the real case, this action degenerates to the action of $\mathbb{Z}_2$ on $Sol_{A}$. In the example below we show that the scenarios in $Sol_{A}$ and $Sol_A^*$ are vastly different. 
 \begin{example}
     Let $A$ be the matrix ${\footnotesize\begin{pmatrix}
    1 & i \\
    i & -1
\end{pmatrix}}$, we know that $Sol_{A}$ is not a group. Indeed, it is a monoid with a set of matrices of the type  ${\footnotesize\begin{pmatrix}
    1-iz & i-iw \\
    z & w
\end{pmatrix}}$ and ${\footnotesize\begin{pmatrix}
    -1-iz & -i-iw \\
    z & w
\end{pmatrix}}$, but $Sol_A^*$ is a group isomorphic to $\Big\{X = {\footnotesize\begin{pmatrix}
    x & y \\
    -y & x
\end{pmatrix}} \, \mid x,y \, \in \mathbb{C} \text{ and } \overline{X}X=1\Big\}. $
 \end{example}

\medskip

%********************************************************************************************************
For $A \in M_n(k)$ where $k$ is an algebraically closed field. Let us write the matrix $A$ as a sum of its diagonal part and nilpotent part as $A_d+ A_{\text{nil}}$. Then we prove a necessary and sufficient condition for $Sol_A$ to be a group.
\begin{lemma}\label{jordan_chevalley_one}
    Let $N$ be a nilpotent triangular matrix (either upper or lower). Then for any vector $v$, $Nv=N^tv$ if and only if $v \in ker(N) \cap ker(N^t)$.
\end{lemma}
\begin{proof}
    We choose a basis for the vector space so that the nilpotent matrix $N$ is in its canonical form. That is, for the chosen basis, the matrix is of the form, 
    ${\footnotesize\begin{pmatrix}
        M_{n_1} & 0 & 0 \cdots & 0\\
        0 & M_{n_2} & 0 \cdots & 0\\
        \vdots & \vdots & \cdots &\vdots \\
        0 & 0 & \cdots & M_{n_k}
    \end{pmatrix}}$ where $M_{n_i}$ is an $n_i \times n_i$ matrix which is either zero or of the type,
    ${\footnotesize\begin{pmatrix}
        0 & 1 & 0 \cdots & 0 &0\\
        0 & 0 & 1 \cdots & 0 & 0\\
        \vdots & \vdots & \cdots & 1 & 0\\
        0 & 0 & \cdots & 0 & 1\\
        0 & 0 & \cdots & 0 &0 
    \end{pmatrix}}$. Also note that for the same basis, the matrix $N^t$ has the form, ${\footnotesize\begin{pmatrix}
        M_{n_1}^t & 0 & 0 \cdots & 0\\
        0 & M_{n_2}^t & 0 \cdots & 0\\
        \vdots & \vdots & \cdots &\vdots \\
        0 & 0 & \cdots & M_{n_k}^t
    \end{pmatrix}}$. So we just have to prove the result for each block $M_{n_i}$. Now for such a block, it is clear that $M_{n_i}v = M_{n_i}^tv$ if and only if $ker(M_{n_i}) \cap ker (M_{n_i}^t)$ is trivial. This completes the proof of lemma.
\end{proof}

\begin{theorem}
Let $A=A_d+A_{\text{nil}}$ be a matrix. The monoid $Sol_{A}$ is a group if and only if $ker(A_d)\cap ker(A_{\text{nil}}) \cap ker(A_{\text{nil}}^t)$ is trivial. In particular, if $Sol_{A_{\text{nil}}}$ is a group, then $Sol_{A}$ is a group for any diagonal matrix $A_d$.
\end{theorem}
\begin{proof}
We know that $Sol_{A}$ is a group if and only if $ker(A) \cap ker(A^t)$ contains only the trivial vector. So it's enough to prove equality between the sets $ker(A) \cap ker(A^t)$ and $ker(A_d)\cap ker(A_{\text{nil}}) \cap ker(A_{\text{nil}}^t)$. Let $v$ be a vector such that $v \in ker(A) \cap ker(A^t)$, then $A_d v+ A_{\text{nil}}v =0 = A_dv+ A_{\text{nil}}^tv$. As a result, we have $A_{\text{nil}}v=A_{\text{nil}}^tv$, which implies $v \in ker(A_{\text{nil}}) \cap ker(A_{\text{nil}}^t)$ by the Lemma \ref{jordan_chevalley_one}. Then $A_dv =0$, as a result we have $ker(A)\cap ker(A^t) = ker(A_{\text{nil}}) \cap ker(A_{\text{nil}}^t) \cap ker(A_d)$. Now, the second assertion follows from this.
\end{proof}

\begin{example}
    For the matrix, $A={\footnotesize\begin{pmatrix}
        0 & 0 & 1\\
        0 & 1 & 0\\
        0 & 0 & 0
    \end{pmatrix}}$, $Sol_{A}$ is a group even though for the associated nilpotent $N= {\footnotesize\begin{pmatrix}
        0 & 0 & 1\\
        0 & 0 & 0\\
        0 & 0 & 0
    \end{pmatrix}}$, $Sol_N$ is not a group as the diagonal matrix is non-singular when restricted to the $ker(N) \cap ker(N^t)=Ke_2$. 
\end{example}
\textbf{Remark.} The above theorem has obvious restrictions in terms of applicability, as the nature of $Sol_{A}$ being a group or not is preserved under similarity relation. It would be interesting to find such a simple result for a general case without the assumption that the matrix is in its Jordan normal form.
%***********************************************************************************************************

\subsection{Congruence matrices} As a consequence of our characterization, one can now classify all the complex matrices $A$ for which $Sol_{A}$ is a group by using the Aitken-Turnbull canonical forms for congruence. Let us first recall the canonical forms of congruence matrices from \cite{Aitken}, (the reader might find this article \cite{Teran_history} useful to recall the canonical forms). The canonical forms of the congruence class of matrices are related to the equivalence of palindromic matrix pencils, $\mu A + \lambda A^t$. 

\begin{theorem}(\cite[Theorem 18]{Aitken})\label{Aitken}  For $c \neq \pm 1$, define $J(c) =  {\footnotesize\begin{bmatrix}
    0&&&1\\
    &&1&c\\
    &\iddots &\iddots&\\
    1&c&&0
\end{bmatrix}}_{k \times k},$ and $I(c) = {\footnotesize\begin{bmatrix}
    0&&&c\\
    &&c&1\\
    &\iddots&\iddots&\\
    c&1&&0
\end{bmatrix}}_{k \times k}$. Then, each matrix $A \in \mathbb{C}^{n \times n}$ is congruent to a direct sum of blocks of the following types:
    \begin{enumerate}
    
        \item[(a)] $\mathcal{A}_1 = 0$; and $\mathcal{A}_{2k+1} = {\tiny \begin{bmatrix}
            
                O_{k+1}
            & \begin{bmatrix}
                1 & & \\
                0&\ddots&\\
                &\ddots&1\\
                &&0
           \end{bmatrix}_{(k+1)\times k}\\
            \begin{bmatrix}
               0&1&&\\
               &\ddots&\ddots&\\
               &&0&1
            \end{bmatrix}_{k \times (k+1)} &
                 O_k
            
        \end{bmatrix}}_{(2k+1)\times (2k+1)}$
        \smallskip
        
        \item[(b)] $\mathcal{B}_{2k}(c)={\footnotesize \begin{bmatrix}
            O_k & J(c)\\
            I(c) & O_k
        \end{bmatrix}_{(2k)\times (2k)}}$
\smallskip

        \item[(c)] $\mathcal{C}_{2k+1} = {\tiny \begin{bmatrix}
            0&&&&&&1\\
            &&&&&1&1\\
            &&&&\iddots&\iddots&\\
            &&&1&1&&\\
            &&1&-1&&&\\
            &\iddots&\iddots&&&&\\
            1&-1&&&&&0
        \end{bmatrix}_{(2k+1)\times (2k+1)}}$  
        \smallskip
        
        \item[(d)] $\mathcal{D}_{2k}= {\footnotesize\begin{bmatrix}
            O_k & J(1)\\
            J(-1) & O_k
        \end{bmatrix}}_{(2k)\times (2k)}$\  ($k$ even)
        \smallskip
        
        \item[(e)] $\mathcal{E}_{2k} = {\tiny \begin{bmatrix}
            0&&&&&&1\\
            &&&&&1&1\\
            &&&&\iddots&\iddots&\\
            &&&1&1&&\\
            &&-1&1&&&\\
            &\iddots&\iddots&&&&\\
            -1&1&&&&&0
        \end{bmatrix}}_{(2k) \times (2k)}$
        \smallskip
        
        \item[(f)] $\mathcal{F}_{2k}= {\footnotesize\begin{bmatrix}
            O_k & I(1)\\
            I(-1) & O_k
        \end{bmatrix}}_{(2k)\times (2k)}$ \  ($k$ odd).
    \end{enumerate}
   
\end{theorem}

\begin{proposition} \label{path}
    Let $A_n$ be the matrix ${\footnotesize\begin{pmatrix}
0 &1  &  & &0\\ 
 & 0 & 1 &  &\\ 
 &  &  \ddots&\ddots &\\ 
 & &  &\ddots &1\\ 
 0&  &  & &0
\end{pmatrix}}$, and recall the matrices, $\mathcal{A}_n$, and $\mathcal{B}_n$ as defined above. Then, we have
    \begin{enumerate}
        \item $A_n \simeq \begin{cases}
            \mathcal{A}_n; & \text{ if } n \text{ is odd},\\
            \mathcal{B}_n; & \text{ if } n \text{ is even}.
        \end{cases}$
        \item $A_n^k \simeq \begin{cases}
            \bigoplus_{\alpha \text{ times}} \mathcal{A}_{n_1} \oplus \bigoplus_{\beta \text{ times}} \mathcal{B}_{n_1 -1}; & \text{ if } n_1 \text{ is odd},\\
            \bigoplus_{\alpha \text{ times}} \mathcal{B}_{n_1} \oplus \bigoplus_{\beta \text{ times}} \mathcal{A}_{n_1 -1}; & \text{ if } n_1  \text{ is even},
        \end{cases}$
    \end{enumerate}
    where $n_1 = \lfloor\tfrac{n-1}{k}\rfloor, \ \alpha = (n-kn_1), \ \beta = \tfrac{n-\alpha(n_1 +1)}{n_1}.$
\end{proposition}
\begin{proof}
To prove the first assertion, let us assume that $n=2k$. Now consider the directed path as in the following figure. 
\begin{center}

    \begin{tikzpicture}[roundnode/.style={circle, fill, inner sep=1pt, minimum size=1mm}]
\node[circle,fill,inner sep=1pt,label={north:1}]      (1) {};
\node[roundnode,label={south:2}]      (2)       [right=of 1] {};
\node[roundnode,label={north:3}]      (3)       [right=of 2] {};
\node[roundnode,label={south:4}]      (4)       [right=of 3] {};
\node[roundnode,label={north:}]      (5)       [right=of 4] {};
\node[roundnode,label={south: $n-3$}]      (g1)       [right=of 5] {};
\node[roundnode,label={north: $n-2$}]      (g2)       [right=of g1] {};
\node[roundnode,label={south: $n-1$}]      (g3)       [right=of g2] {};
\node[roundnode,label={north: $n$}]      (n)       [right=of g3] {};
\draw[->] (1) -- (2);
\draw[->] (2) -- (3);
\draw[->] (3) -- (4);
\draw[dashed] (4) -- (5);
\draw[dashed] (5) -- (g1);
\draw[->] (g1) -- (g2);
\draw[->]  (g2) -- (g3);
\draw[->] (g3) -- (n);
\end{tikzpicture}
\end{center}
Note that the adjacency matrix of this path is precisely the matrix $A_n$. Similarly, let us consider the directed path of the matrix $\mathcal{B}_{2k}$, and we obtain the following.

\begin{center}
        \begin{tikzpicture}[roundnode/.style={circle, fill, inner sep=1pt, minimum size=1mm}]
\node[circle,fill,inner sep=1pt,label={north:1}]      (1) {};
\node[roundnode,label={south:$2k$}]      (2)       [right=of 1] {};
\node[roundnode,label={north:2}]      (3)       [right=of 2] {};
\node[roundnode,label={south:$2k-1$}]      (4)       [right=of 3] {};
\node[roundnode,label={north:}]      (5)       [right=of 4] {};
\node[roundnode,label={south: $k-1$}]      (g1)       [right=of 5] {};
\node[roundnode,label={north: $k+2$}]      (g2)       [right=of g1] {};
\node[roundnode,label={south: $k$}]      (g3)       [right=of g2] {};
\node[roundnode,label={north: $k+1$}]      (n)       [right=of g3] {};
\draw[->] (1) -- (2);
\draw[->] (2) -- (3);
\draw[->] (3) -- (4);
\draw[dashed] (4) -- (5);
\draw[dashed] (5) -- (g1);
\draw[->] (g1) -- (g2);
\draw[->]  (g2) -- (g3);
\draw[->] (g3) -- (n);
\end{tikzpicture}
\end{center}
Clearly, this is isomorphic to the path of $A_{2k}$, so that must be a relabeling that gives this isomorphism. The isomorphism is given by the permutation matrix $\sigma_{n}$ whose $ij^\text{th}$ entry is $$(\sigma_{n})_{ij} = \begin{cases}
       1, & \text{ if } i = 2a+1, \ j = a+1;\\
       1, & \text{ if } i = 2a, \ j = n-a+1;\\
       0 & \text{ otherwise.}
   \end{cases}$$
   As a result, we obtain \[\sigma_{n}^{-1}A_{n} \sigma_{n}=\sigma_{n} ^t A_{n} \sigma_n = \mathcal{B}_{n}. \]\\
   Now for $n=2k+1$-odd, consider the path with labels as follows. 
   
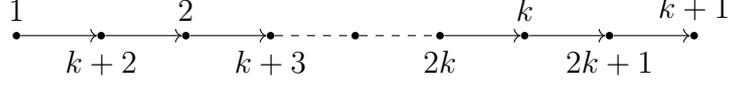
\begin{center}
        \begin{tikzpicture}[roundnode/.style={circle, fill, inner sep=1pt, minimum size=1mm}]
\node[circle,fill,inner sep=1pt,label={north:1}]      (1) {};
\node[roundnode,label={south:$k+2$}]      (2)       [right=of 1] {};
\node[roundnode,label={north:2}]      (3)       [right=of 2] {};
\node[roundnode,label={south:$k+3$}]      (4)       [right=of 3] {};
\node[roundnode,label={north:}]      (5)       [right=of 4] {};
\node[roundnode,label={south: $2k$}]      (g1)       [right=of 5] {};
\node[roundnode,label={north: $k$}]      (g2)       [right=of g1] {};
\node[roundnode,label={south: $2k+1$}]      (g3)       [right=of g2] {};
\node[roundnode,label={north: $k+1$}]      (n)       [right=of g3] {};
\draw[->] (1) -- (2);
\draw[->] (2) -- (3);
\draw[->] (3) -- (4);
\draw[dashed] (4) -- (5);
\draw[dashed] (5) -- (g1);
\draw[->] (g1) -- (g2);
\draw[->]  (g2) -- (g3);
\draw[->] (g3) -- (n);
\end{tikzpicture}
   \end{center}
   Observe that this directed path is isomorphic to the path $P_n$ so let us take the relabeling matrix $\sigma_n$. It is an easy check to see that the matrix $\sigma_n$ is given by its entries as follows, $ij^\text{th}$ entry is $$(\sigma_n)_{ij} = \begin{cases}
       1, & \text{ if } i = 2a, \ j = a+\lfloor\tfrac{n}{2}\rfloor +1;\\
       1, & \text{ if } i = 2a+1, \ j = a+1;\\
       0 & \text{ otherwise,}
   \end{cases}$$
   and one checks immediately that for this matrix, we obtain 
   \[\sigma_n ^t A_n \sigma_n = \mathcal{A}_{n}.\]
   This completes the proof of (1). For the second assertion, consider a directed path and use the observation that the powers of the matrix $A_n$, give the adjacency matrix of the same power of the directed path. That is, $A_n^k$ is the adjacency matrix of $P_n^k$ where $P_n^k$ is the $k$-th power of the path graph $P_n$. One can verify that $P_n^k = \alpha P_{n_1} \cup \beta P_{n_1 -1}$, where $\alpha = \left(n-k\lfloor\tfrac{n-1}{k}\rfloor\right), \ \beta = \tfrac{n-\alpha(n_1 +1)}{n_1}$, and $aP_n$ is the union of $\alpha$ copies of $P_n$. Using (1), we can acquire permutation matrices $\sigma_{n_1}$ and $\sigma_{n_1-1}$ that give the labeled isomorphism of the paths in terms of a direct sum to the paths below.
   
   
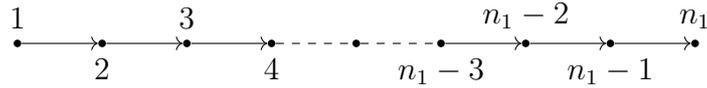
\begin{center}

    \begin{tikzpicture}[roundnode/.style={circle, fill, inner sep=1pt, minimum size=1mm}]
\node[circle,fill,inner sep=1pt,label={north:1}]      (1) {};
\node[roundnode,label={south:2}]      (2)       [right=of 1] {};
\node[roundnode,label={north:3}]      (3)       [right=of 2] {};
\node[roundnode,label={south:4}]      (4)       [right=of 3] {};
\node[roundnode,label={north:}]      (5)       [right=of 4] {};
\node[roundnode,label={south: $n_1-3$}]      (g1)       [right=of 5] {};
\node[roundnode,label={north: $n_1-2$}]      (g2)       [right=of g1] {};
\node[roundnode,label={south: $n_1-1$}]      (g3)       [right=of g2] {};
\node[roundnode,label={north: $n_1$}]      (n)       [right=of g3] {};
\draw[->] (1) -- (2);
\draw[->] (2) -- (3);
\draw[->] (3) -- (4);
\draw[dashed] (4) -- (5);
\draw[dashed] (5) -- (g1);
\draw[->] (g1) -- (g2);
\draw[->]  (g2) -- (g3);
\draw[->] (g3) -- (n);
\end{tikzpicture}
\end{center}
and 
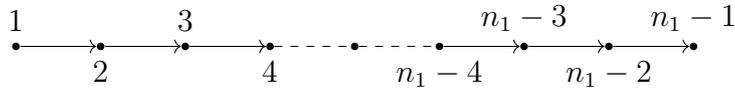
\begin{center}

    \begin{tikzpicture}[roundnode/.style={circle, fill, inner sep=1pt, minimum size=1mm}]
\node[circle,fill,inner sep=1pt,label={north:1}]      (1) {};
\node[roundnode,label={south:2}]      (2)       [right=of 1] {};
\node[roundnode,label={north:3}]      (3)       [right=of 2] {};
\node[roundnode,label={south:4}]      (4)       [right=of 3] {};
\node[roundnode,label={north:}]      (5)       [right=of 4] {};
\node[roundnode,label={south: $n_1-4$}]      (g1)       [right=of 5] {};
\node[roundnode,label={north: $n_1-3$}]      (g2)       [right=of g1] {};
\node[roundnode,label={south: $n_1-2$}]      (g3)       [right=of g2] {};
\node[roundnode,label={north: $n_1-1$}]      (n)       [right=of g3] {};
\draw[->] (1) -- (2);
\draw[->] (2) -- (3);
\draw[->] (3) -- (4);
\draw[dashed] (4) -- (5);
\draw[dashed] (5) -- (g1);
\draw[->] (g1) -- (g2);
\draw[->]  (g2) -- (g3);
\draw[->] (g3) -- (n);
\end{tikzpicture}
\end{center}

Now let us consider the permutation matrix $\sigma$ which is a block diagonal matrix with the diagonals, $\sigma_{n_1}$ for $\alpha$ many times, and $\sigma_{n_1-1}$ for $\beta$ many times. That is $$\sigma = diag(\underbrace{\sigma_{n_1}, \ldots, \sigma_{n_1}}_\text{$\alpha$ times}, \underbrace{\sigma_{n_1-1}, \ldots, \sigma_{n_1-1}}_\text{$\beta$ times}).$$ 
\end{proof}
\begin{theorem}
    Let $A$ be a matrix that is congruent to a direct sum $\displaystyle\bigoplus_{\alpha \in \Gamma} M_{\alpha}$, where $\Gamma$ is a finite set and $M_{\alpha}$ is one of the six classes of matrices in the Aitken-Turnbull classification \ref{Aitken}. Then $Sol_A$ is a group if and only if $\mathcal{A}_1$ is not a summand.
\end{theorem}
\begin{proof}
    We will apply our Theorem \ref{group} to prove the result. Note that the condition $ker(A) \cap ker(A^t)$ obeys the direct sum. Also, observe that the summands which are invertible do not at all contribute to this condition. Except for the classes $\mathcal{A}_{2k+1}$ and $\mathcal{B}_{2k}$, rest of the matrices are invertible. Hence we can ignore the classes $\mathcal{C}_{2k+1}, \mathcal{D}_{2k}, \mathcal{E}_{2k}, \mathcal{F}_{2k}$, as in the Theorem \ref{Aitken}. Now, using the above Proposition \ref{path}, we know that $\mathcal{A}_n$ and $\mathcal{B}_n$ are congruent to $A_n$ for $n>1$. Furthermore, from Theorem \ref{group}, we can check that $A_n$ is a group for $n>1$. Clearly, if $\mathcal{A}_1$ is a summand, then $Sol_A$ is not a group. This proves the Theorem.
\end{proof}
\begin{corollary} \label{group_nilpotent}
    $Sol_{A_n^k}$ is a group whenever $\lfloor \frac{n-1}{k} \rfloor $ is strictly greater than one.
\end{corollary}

\section{The structure of $Sol_{A_n}$}\label{Sec4}
Using corollary \ref{group_nilpotent} for $k =1$, we know that for the matrix $A_n$, $Sol_{A_n}$ is a group for all $n >2$. Also one can check $Sol_{A_2}$ is also a group. In this section, we will find the structure of $Sol_{A_n} = \{g \in M_n(\mathbf{C}) \ | \ g^tA_ng = A_n\}$ by utilizing the structure of its Lie algebra. Before dealing with general $n$, we will discuss one prototypical example.\\

\subsection{Example.} We now provide a full running example demonstrating the core idea to write down the equations of the groups. This example can help the reader understand this section's main theme. Let us choose $n$ to be equal to 8. First we observe that the Lie algebra $sol_{A_8} = \{X | X^tA_8+A_8X = 0\}$ of $\text{Sol}_{A_8}$ consists of matrices as below where $x_0, x_1,x_2,x_3 $ are free parameters $${X = \tiny\begin{pmatrix}
    x_0 & 0& x_1&0&x_2&0&x_3&0\\
    0&-x_0&0&0&0&0&0&0\\
    0&0&x_0&0&x_1&0&x_2&0\\
    0&-x_1&0&-x_0&0&0&0&0\\
    0&0&0&0&x_0&0&x_1&0\\
    0&-x_2&0&-x_1&0&-x_0&0&0\\
    0&0&0&0&0&0&x_0&0\\
    0&-x_3&0&-x_2&0&-x_1&0&-x_0  
\end{pmatrix}}$$
\smallskip

As a consequence of the above, we can take the generators of the Lie algebra $sol_{A_8}$ as $$h=E_{11}-E_{22}+E_{33}-E_{44}+E_{55}-E_{66}+E_{77}-E_{88}$$ $$e_1= E_{13}+E_{35}+E_{57}-E_{42}-E_{64}-E_{86}$$ $$e_2= E_{15}+E_{37}-E_{62}-E_{84}$$ $$e_3=E_{17}-E_{82}.$$ Now observe that $e_1,e_2,e_3$ generate a Lie ideal, which is also the unipotent radical. Then, let us take the normal subgroup $N$ corresponding to this ideal, which is obtained by exponentiating this ideal. Or equivalently throwing in the non-linear parts of the generating equations. Let us take the diagonal subgroup as $D$ which is the set of matrices $diag(t,t^{-1},t, t^{-1},t , t^{-1}, t, t^{-1})$ for $t \in \mathbb{C}$. As a result, we have the group $Sol_{A_8}$ as a semi-direct product of $D$ and $N$, and putting these together, we have the structure of a general element of the group, $\text{Sol}_{A_8}$, as below,
 $${g = \tiny \begin{pmatrix}
    x_0 & 0& x_0x_1&0&x_0x_2&0&x_0x_3&0\\
    0&x_0^{-1}&0&0&0&0&0&0\\
    0&0&x_0&0&x_0x_1&0&x_0x_2&0\\
    0&-x^{-1}_0x_1&0&x_0^{-1}&0&0&0&0\\
    0&0&0&0&x_0&0&x_0x_1&0\\
    0&-x^{-1}_0(x_2 -x_1^2)&0&-x^{-1}_0x_1&0&x_0^{-1}&0&0\\
    0&0&0&0&0&0&x_0&0\\
    0&-x^{-1}_0(x_3+2x_1x_2-x_1^3)&0&-x^{-1}_0(x_2 -x_1^2)&0&-x^{-1}_0x_1&0&x_0^{-1}  
\end{pmatrix}}$$
\smallskip

\subsection{Lie algebra of $\text{Sol}_{A_n}$}
The Lie algebra of $Sol_{A_n}$ is related to the works of F. D. Teran et. al. \cite[Lemma 4]{Teran_Lie}, where the authors solve the tangent space equation $XA_n+A_nX^t$ to find the dimension of the congruence orbits. In this section, we will find the structure of the Lie algebra of $Sol_{A_n}$, denoted as $sol_{A_n}$. Recall $sol_{A_n}= \{X \ | \ X^tA_n + A_nX=0\} $.
\begin{lemma}\label{odd_case_A}
When $n \equiv 1(\text{mod } 2).$ Then the following types of matrices with $x_i \in \mathbf{C}$ constitute the solutions of the equation, $ X^tA_n + A_nX=0$.

{\tiny$$X = \begin{pmatrix}
x_0 &x_1  & 0 & x_2&0&\cdots&x_{\tfrac{n-1}{2}}& 0\\ 
 0& -x_0 &0& 0 & 0& 0&\ddots &0\\ 
0 & -x_1 &  x_0&x_1 & 0&x_2&\ddots&\vdots\\ 
 0& 0 &  0&\ddots &\ddots&\ddots&\ddots&\vdots\\ 
0&-x_2&\ddots&\ddots&\ddots&\ddots&0&0\\
 \vdots& \ddots & \ddots & \ddots&\ddots&x_0&x_1&0\\
0&\ddots&0&0&0&0&-x_0&0\\
0 &-x_{\tfrac{n-1}{2}}&\ldots&0&-x_2&0&-x_1&x_0
\end{pmatrix}$$}

\end{lemma}
\begin{proof}
These are indeed valid solutions for the equation $X^tA_n + A_nX =0$ can be shown with some simple calculation. Therefore, we will demonstrate the generation process below.

Suppose $X$ is a solution to the equation $X^tA_n+A_nX=0$. Let $C_i(X)$ represent the $i$th column of $X$, and $R_i(X)$ represent the $i$th row. Observe that for $i > 1$, we have $C_i(X^tA_n) = C_{i-1}(X^t) = R_{i-1}(X)$ and $C_1(X^tA_n) = 0$. Similarly, $R_i(XA) = R_{i+1}(X)$ for $i \leq n-1$, and $R_n(X^tA_n) = 0$. Consequently, we have $(X^tA_n+XA_n)_{i,j} = (X^tA_n)_{i,j}+(A_nX)_{i,j} = X_{j-1, i} + X_{i+1, j}$ for $1<j\leq n$ and $1 \leq i < n$. This leads to the following equations:
\[ X_{j-1,i}+X_{i+1,j}=0 \mbox{ for } 1< j \leq n  , 1 \leq i < n\]
\[X_{i,1}=0 \mbox{ for any } i \]
\[X_{i,n}=0 \mbox{ for any } i.\]
Then, we have \begin{equation*} X_{i,j}=- X_{j-1,i-1}=X_{i-2,j-2}.\end{equation*}
Due to this observation, if $i$ is even and $j$ is odd (assuming $i<j$), we find that $X_{i,j} = X_{i+2,j+2} = X_{i+2k, j+2k}$ for any $k$ within the range. Consequently, for a suitable $k$, we will have $j+2k = n$ (since $n$ is odd) and therefore $X_{i+2k, j+2k} = 0$. If $i$ is even and $j$ is odd and $i > j$, then we can consider $X_{i,j} = X_{i-2,j-2} = X_{i-2k, j-2k}$, with a similar reasoning as before. There will be some $k$ such that $j-2k = 1$, which leads to $X_{i,j} = X_{i-2k,1} = 0$.

When $i$ is odd and $j$ is even, first assume $i < j$. We then have $X_{i,j} = X_{i-2k,j-2k}$ for any $k$ in the range. For a certain $k$, we will have $i-2k = 1$, resulting in $X_{i,j} = X_{1,j-2k} = x_{j-2k/2}$. If $i > j$, we find that $X_{i,j} = -X_{j-1,i-1}$, which corresponds to the case where the first index is odd, the second even, and the first one is smaller. This situation has already been addressed in the first part of the proof. This concludes the proof of generation.
\end{proof}
Clearly, each $X \in \text{sol}_{A_n}$ can be further decomposed as $\mathfrak{d} + \mathfrak{n}$, where $\mathfrak{d}$ is the diagonal matrix $\text{diag}(x_0, -x_0, \ldots, x_0, -x_0, x_0)$ and $\mathfrak{n}$ is the matrix $X$ with $X_{ii} = 0$.
Let us define the Lie algebra elements $$h:=\sum_{i=1}^n (-1)^{i+1}E_{ii};$$ and $$e_k: =\sum_{l=1}^{(n-2k+1)/2}( E_{2l-1,2k+2l-2} - E_{2k+2l-1, 2l}), \ \ \ \text{ for } k = 1, \ldots, \tfrac{n-1}{2}.$$ 
\begin{lemma}\label{A_brackets_odd}
    The Lie brackets are as follows:
    \begin{center}
        \begin{tabular}{c|c c}
             & $h$&$e_k$ \\
             \hline
            $h$ & $0$&$2e_k$\\
            $e_{k'}$ &$-2e_{k'}$&$0$\\
        \end{tabular}
    \end{center}
    This makes the algebra $sol_{A_n}$ nilpotent.
\end{lemma}
\begin{proof}
$[h, e_k] = \sum_{i} (-1)^{i+1} \sum_{l} ( [E_{ii}, E_{2l-1, 2k+2l-2}] - [E_{ii}, E_{2k+2l-1, 2l}]) $. Once we simplify the first bracket expressions we obtain, $ \sum_{l} (E_{2l-1, 2k+2l-2} + E_{2l-1,2k+2l-2}) = 2 \sum_{l} E_{2l-1, 2k+2l-2}$. Similarly, the second bracket simplifies to $2 \sum_{l} E_{2k+2l-1, 2l} $  Thus $[h,e_k]=2e_k$.

For the Lie bracket $[e_k, e_{k'}]$, observe that the indices in $e_k=\sum_{l=1}^{(n-2k+1)/2}( E_{2l-1,2k+2l-2} - E_{2k+2l-1, 2l})$ are (odd, even) so when we take bracket of two of these types, it results in zero. 
\end{proof}

\begin{lemma}\label{even_case_A}
When $n \equiv 0 (\text{ mod } 2).$ Then the following matrices constitute the space of solutions of the equation $ X^tA_n + A_nX=0$.
$$ X = {\tiny\begin{pmatrix}
x_0 &0&x_1  & 0 & x_2&0&\cdots&x_{\tfrac{n-2}{2}}& 0\\ 
 0& -x_0 &0& 0 & 0& 0&0&\ddots &0\\ 
0 & 0 &  x_0&0&x_1 & 0&x_2&\ddots&0\\ 
 0& -x_1 &  0&-x_0&\ddots &\ddots&\ddots&\ddots&\vdots\\ 
 0&0&\ddots&\ddots&\ddots&\ddots&\ddots&\ddots&\vdots\\
\vdots&-x_2&\ddots&\ddots&\ddots&\ddots&\ddots&\ddots&\vdots\\
\vdots&\ddots&\ddots&\ddots&\ddots&\ddots&\ddots&\ddots&\vdots\\
0&\ddots&\ddots&0&0&0&0&x_0&0\\
0 &-x_{\tfrac{n-2}{2}}&\ldots&0&-x_2&0&-x_1&0&-x_0
\end{pmatrix}}$$
\end{lemma}
\begin{proof}
The proof of this case is similar to the previous Lemma \ref{odd_case_A}, so it is omitted.
\end{proof}
Let us define the Lie algebra elements $$h:=\sum_{i=1}^n (-1)^{i+1}E_{ii};$$
and
$$e_k: =\sum_{l=1}^{(n-2k-1)/2}( E_{2l-1,2k+2l-1} - E_{2k+2l, 2l}), \ \ \text{ for } k = 1, \ldots, \tfrac{n-2}{2}.$$
\begin{lemma}
    The Lie brackets are as follows:
    \begin{center}
        \begin{tabular}{c|c c}
             & $h$&$e_k$ \\
             \hline
            $h$ & $0$&$0$\\
            $e_{k'}$ &$0$&$0$\\
        \end{tabular}
    \end{center}
    This makes the algebra $sol_{A_n}$ abelian, and hence nilpotent.
\end{lemma}
\begin{proof}
    This proof is exactly similar to that of the Lemma \ref{A_brackets_odd}, so it is omitted.
\end{proof}
With the help of above lemmas, we can conclude the following:
\begin{proposition}\label{radical_A}
    For any $n$, the algebra $sol_{A_n}$ is solvable; hence, the radical (the maximal solvable ideal) coincides with the algebra.
\end{proposition}
\subsection{Lie group of $\text{Sol}_{A_n}$}
In this section, we will show that the group $Sol_{A_n}$ splits as a semidirect product of two subgroups $D, N$. 
   \subsubsection{ The case of $n \equiv 1(\text{mod}\ 2)$}  Define $D = \{\text{diag}(x_0, x_0^{-1}, x_0, \cdots, x_0^{-1}, x_0) \ | \ x_0 \in \mathbb{C}^* \},$ and $ N = n \times n$ matrices whose columns $C_i$, with value 1 at $i^{th}$ position, are defined as 
    $$\text{For $i$ even}, \ \ C_i = {\footnotesize\begin{bmatrix}
     x_{\tfrac{i}{2}}  \\
     \vdots\\
     0\\
     x_1\\
     1\\
     -x_1\\
     0\\
     \vdots\\
     -x_{\tfrac{n-i+1}{2}}
\end{bmatrix}}; \ \ \text{For $i$ odd}, \ \ C_i = {\footnotesize\begin{bmatrix}
     0\\
     \vdots\\
     0\\
     1\\
     0\\
     \vdots\\
     0
\end{bmatrix}},$$
where $x_i \in \mathbb{C}$. Then, for any $i, \widetilde{C}_i$ represents the truncated column or $A_nC_i$, which is obtained from column $C_i$ by shifting each value at $k^{th}$ index to $(k-1)^{th}$ index, and append zero at the last entry.
\begin{lemma} \label{diagonal}
    The above-defined $D$ and $N$ are the subgroups of $\text{Sol}_{A_n}.$
\end{lemma}
\begin{proof} \begin{enumerate}
    \item \textbf{Claim.} $D$ is a subgroup of $Sol_{A_n}$.\\
    
    We need to prove for any $d \in D,$ we have $d^tA_nd = A_n.$ It is same as proving $dA_n = A_nd^{-1}$, where $d^{-1} = \text{diag}(x_0^{-1}, x_0, x_0^{-1} \cdots, x_0, x_0^{-1})$. Clearly $(dA_n)_{i,j} = d_{i, j-1},$ for $j \geq 2$ and $(A_nd^{-1})_{i,j} = d^{-1}_{i+1, j}$ for $i\leq n-1.$ Therefore $$(dA_n)_{i,j} = \begin{cases}
              0 & \text{if } i \neq j-1\\
              x_0 & \text{if } i \text{ is odd}\\
              x_0^{-1} & \text{if } i \text{ is even}
\end{cases}$$
and $$(A_nd^{-1})_{i,j} = \begin{cases}
             0 & \text{if } i+1 \neq j\\
              x_0^{-1} & \text{if } i \text{ is even}\\
              x_0 & \text{if } i \text{ is odd}
\end{cases}$$
which are equal. Furthermore, since $D$ consists of the diagonal matrices, for any $d_1, d_2 \in D$, we have $d_1d_2^{-1} \in D.$\\
\item \textbf{Claim.} $N$ is a subgroup of $Sol_{A_n}$.\\
We need to prove for any $X \in N$, we have $X^tA_nX = A_n.$ It's enough to prove $C_i^t \widetilde{C}_j = 1$ for $j = i+1$, and otherwise it is zero. For $j = i+1$ we get
$$C_i^t \widetilde{C}_{i+1} = {\footnotesize \begin{bmatrix}
     x_{\tfrac{i}{2}} &
     \cdots&
     0&
     x_1&
     1&
     -x_1&
     0&
     \cdots&
     -x_{\tfrac{n-i+1}{2}}
\end{bmatrix} \begin{bmatrix}
     0\\
     \vdots\\
     0\\
     1\\
     0\\
     \vdots\\
     0
\end{bmatrix} }= 1.$$
Because for both $C_i^t$ and $\widetilde{C}_{i+1}$, the value 1 is at the $i^\text{th}$ place. When $i$ and $j$ are even such that $i > j$, then $C_i^t \widetilde{C}_j = 0$ because
\[ {\tiny\begin{matrix}
\\ \\ (j-1)^{\text{th}} \text{ place}  \\ \\ \\ \\ i^{\text{th}} \text{ place}  \\  \\ \\ \\
\end{matrix}\begin{matrix}
\\ \\ \\ \\ \longrightarrow \\  \\ \\ \\  \longrightarrow \\ \\ \\ \\ \\ \\ 
\end{matrix}\begin{bmatrix}
     x_{\tfrac{i}{2}}  \\ 0\\ \vdots \\x_{\tfrac{i-j+2}{2}}\\\vdots\\ 0\\x_1 \\1 \\  -x_1\\0 \\ \vdots\\ 0 \\ -x_{\tfrac{n-i+1}{2}}
\end{bmatrix}^t\begin{bmatrix}
     0\\x_{\tfrac{j-2}{2}}  \\\vdots \\ 0\\x_1 \\1 \\  -x_1\\0 \\ \vdots\\ -x_{\tfrac{i-j+2}{2}} \\ \vdots \\ -x_{\tfrac{n-j+1}{2}} \\0
\end{bmatrix} \begin{matrix}
\\ \\ \\ \\ \\ \\ \longleftarrow  \\ \\ \\ \\ \\ \longleftarrow \\ \\ \\ \\ \\
\end{matrix}\begin{matrix}
\\ \\ \\(j-1)^{\text{th}} \text{ place}  \\ \\ \\ \\ \\  i^{\text{th}} \text{ place}  \\  \\ \\
\end{matrix}}\]
$$= x_{\tfrac{i-j+2}{2}} + x_{\tfrac{i-j-1}{2}}(-x_1) + \cdots + x_1(-x_{\tfrac{i-j-1}{2}}) - x_{\tfrac{i-j+2}{2}},$$
which is clearly zero. For the other remaining case, we can similarly see the product is zero.\\

For any $X, Y \in N$, we need to show $XY^{-1} \in N.$ The rows of any element $X$ in $N$ are as follows:
$$R_i = \begin{cases}
    [0, -x_{\tfrac{i-1}{2}}, \ldots, 0, -x_1, 1, x_1, 0, \ldots, x_{\tfrac{n-i}{2}}, 0] & \text{ if } i \text{ is odd}\\
    [0, \ldots, 0, 1, 0, \ldots, 0] & \text{ if } i \text{ is even}
\end{cases}$$
Then the product $$(XY)_{ij} = \begin{cases} x_{\tfrac{j-i+1}{2}} +y_{\tfrac{j-i+1}{2}} & \text{ if $i$-odd, $j$-even and $i<j$} \\
-x_{\tfrac{i-j+1}{2}} -y_{\tfrac{i-j+1}{2}} & \text{ if $i$-odd, $j$-even and $i>j$}\\
\delta_{ij} & \text{ if } i \equiv j (\text{mod} 2)\\
0 &  \text{ if $i$-even, $j$-odd},
\end{cases}$$
which clearly shows $XY \in N.$ Now for any $X \in N$, the \textbf{Inverse of $X$} is denoted by $X^{-1}$ and is defined as $$(X^{-1})_{ij} = \begin{cases}
    -X_{ij} & \text{ if } i \neq j;\\
    1 & \text{ if } i=j.
\end{cases}$$
Clearly $X^{-1} \in N$. This completes the proof of the lemma.
\end{enumerate}
\end{proof}

\subsubsection{The case of $n \equiv 0(\text{mod } 2)$}  Define $D = \{\text{diag}(x_0, x_0^{-1}, \cdots, x_0, x_0^{-1}) \ | \ x_0 \in \mathbb{C}^*\},$ and $ N = n \times n$ matrices whose columns $C_i$, with value 1 at $i^{th}$ position, are 
    $$\text{For $i$ even}, \ \ C_i = {\footnotesize\begin{bmatrix}
     0\\
     \vdots\\
     0\\
     1\\
     0\\
     -x_1 + x_1^\ast\\
     0\\
     \vdots\\
     -x_{\tfrac{n-i}{2}} + x_{\tfrac{n-i}{2}}^\ast
\end{bmatrix}}; \ \ \text{For $i$ odd}, \ \ C_i = {\footnotesize\begin{bmatrix}
     x_{\tfrac{i-1}{2}}\\
     \vdots\\
     0\\
     x_1\\
     0\\
     1\\
     0\\
     \vdots\\
     0
\end{bmatrix}},$$
where $x_i \in \mathbb{C}$. Put $x_N^\ast = \sum_{a+b= N}x_a(x_b - x_b^\ast)$ such that $x_1^\ast = 0.$
\begin{lemma}
    The above-defined $D$ and $N$ are the subgroups of $\text{Sol}_{A_n}.$
\end{lemma}
\begin{proof}
By following the same proof as in Lemma \ref{diagonal}, we can show $D$ is a subgroup of $\text{Sol}_{A_n}.$ Furthermore, we prove $N$ is a subgroup of $Sol_{A_n}$.\\

Recall $\tilde{C}_i = A_nC_i$. We show for every matrix $X \in N$ it satisfies $X^tA_nX =  A_n,$ which implies $N \subseteq \text{Sol}_{A_n}.$ For $j = i+1$, we have
$$C_i^t \widetilde{C}_{i+1} = 
    {\tiny \begin{bmatrix}
     x_{\tfrac{i-1}{2}} &
     \cdots&
     0&
     x_1&
     0&
     1&
     0&
     \cdots&
     0
\end{bmatrix} \begin{bmatrix}
     0\\
     \vdots\\
     0\\
     1\\
     -x_1 + x_1^\ast\\
     0\\
     \vdots\\
     -x_{\tfrac{n-i}{2}} + x_{\tfrac{n-i}{2}}^\ast
\end{bmatrix}} = 1, \ \ \text{ if } i \text{ is odd};$$
$$= \ {\tiny \begin{bmatrix}
     0&
     \cdots&
     0&
     1&
     0&
     -x_1+x_1^\ast&
     0&
     \cdots&
     -x_{\tfrac{n-i}{2}} + x_{\tfrac{n-i}{2}}^\ast 
\end{bmatrix} \begin{bmatrix}
     0\\
     x_{\tfrac{i-2}{2}}\\
     \vdots\\
     x_1\\
     0\\
     1\\
     0\\
     \vdots\\
     0
\end{bmatrix} }=1 \ \ \text{ if } i \text{ is even}.$$
Furthermore, for $j \neq i+1$, due to the positions of zeros, we directly get
$C_i^t \widetilde{C}_j = 0$ except the cases $i$ is odd, $j$ is even such that $i>j$; and $i$ is even, $j$ is odd such that $i<j$. For the case $i$-odd, $j$-even, $i>j$, we get
$$C_i^t\widetilde{C}_j ={\tiny  \begin{matrix}
\\ \\ \\ \\ (j-1)^{\text{th}} \text{ place} \\ \\ \\ \\ \\ i^{\text{th}} \text{ place}  \\  \\ \\ \\
\end{matrix}\begin{matrix}
\\  \\ \\ \\\longrightarrow \\  \\ \\ \\ \\  \longrightarrow \\ \\ \\ \\ 
\end{matrix}\begin{bmatrix}
     x_{\tfrac{i-1}{2}}  \\ 0\\ \vdots \\x_{\tfrac{i-j+1}{2}}\\\vdots\\ 0\\x_1 \\0\\1 \\  0\\ \vdots \\ 0
\end{bmatrix}^t\begin{bmatrix}
     0\\ \vdots \\ 0 \\1 \\ 0\\ -x_1+x_1^\ast\\0 \\ \vdots\\ -x_{\tfrac{i-j+1}{2}}+ x_{\tfrac{i-j+1}{2}}^\ast \\ \vdots \\ -x_{\tfrac{n-j+1}{2}}+ x_{\tfrac{n-j+1}{2}}^\ast \\0
\end{bmatrix} \begin{matrix}
\\ \\ \longleftarrow \\ \\ \\ \\ \\ \\ \longleftarrow \\ \\ \\ 
\end{matrix}\begin{matrix}
\\ \\ \\(j-1)^{\text{th}} \text{ place} \\ \\ \\ \\ \\ \\  i^{\text{th}} \text{ place}  \\  \\ \\
\end{matrix}} =0.$$
Similarly, one can check for other cases where the product is zero.
For matrices $X$ and $Y$ in $N$, we need to show that $XY \in N$. The rows of any matrix $X \in N$ denotes as $R_i$ where 1 is at $i^{th}$ position and are given as follows:
    $$R_i = \begin{cases} [0, \cdots,0, 1, 0, x_1, 0, x_2, \cdots, x_{\tfrac{n-i-1}{2}}, 0]; & i \text{ is odd}\\
    [0, -x_{\tfrac{i}{2}-1} + x^\ast_{\tfrac{i}{2}-1},\cdots, -x_1+x_1^\ast,0, 1, 0, \cdots, 0]; & i \text{ is even}.
    \end{cases}$$
    Then $ij^{th}$ of the product matrix $XY$ is 
    $${\footnotesize(XY)_{ij} = \begin{cases}
        \sum_{a+b = \tfrac{j-i}{2}} x_a y_b + x_{\tfrac{j-i}{2}} + y_{\tfrac{j-i}{2}}; & i<j, \ \ i,j \text{ are odd}\\
        \sum_{a+b = \tfrac{i-j}{2}} (-x_a+x_a^\ast)(- y_b+ y_b^\ast) -x_{\tfrac{i-j}{2}}+x^\ast_{\tfrac{i-j}{2}} - y_{\tfrac{i-j}{2}} +y^\ast_{\tfrac{i-j}{2}} & i>j, \ \ i,j \text{ are even}\\
        0; & \text{ else}
    \end{cases}}$$
    Put $z_n = \sum_{a+b = n} x_a y_b + x_n+ y_n$. To prove $XY \in N$, we need to prove 
    $$-z_n +z^\ast_n = \sum_{a+b = n} (-x_a+x_a^\ast)(- y_b+ y_b^\ast) -x_n + x_n^\ast - y_n+ y_n^\ast.$$
    That is,
    $$z^\ast_n = \sum_{a+b = n} (-x_a+x_a^\ast)(- y_b+ y_b^\ast) +\sum_{a+b = n} x_a y_b + x_n^\ast + y_n^\ast$$
    or $$(-1)\sum_{a+b= n}z_a(-z_b + z_b^\ast) = \sum_{a+b = n} (-x_a+x_a^\ast)(- y_b+ y_b^\ast) +\sum_{a+b = n} x_a y_b + x_n^\ast + y_n^\ast,$$
    which holds because of Theorem \ref{star_one}.
 \end{proof}
For $X \in N$, we denote the inverse as matrix $X^{-1}$ whose columns are as follows:
$$\text{For $i$ even}, \ \ C_i = {\footnotesize\begin{bmatrix}
     0\\
     \vdots\\
     0\\
     1\\
     0\\
     x_1\\
     0\\
     \vdots\\
     x_{\tfrac{n-i}{2}}
\end{bmatrix}}; \ \ \text{For $i$ odd}, \ \ C_i = {\footnotesize\begin{bmatrix}
     -x_{\tfrac{i-1}{2}}+ x_{\tfrac{i-1}{2}}^\ast\\
     \vdots\\
     0\\
     -x_1+x_1^*\\
     0\\
     1\\
     0\\
     \vdots\\
     0
\end{bmatrix}}.$$
Clearly $ X^{-1} \in N$. Furthermore, 
$$(XX^{-1})_{ij} = \begin{cases}
1; & i=j\\
        -x_{\tfrac{j-i}{2}}+x^*_{\tfrac{j-i}{2}}+\sum_{a+b = \tfrac{j-i}{2}} x_a (-x_b+x^*_b) + x_{\tfrac{j-i}{2}}; & i<j, \ \ i,j \text{ are odd}\\
       -x_{\tfrac{i-j}{2}}+x^\ast_{\tfrac{i-j}{2}}+ \sum_{a+b = \tfrac{i-j}{2}} (-x_a+x_a^\ast)x_b +x_{\tfrac{i-j}{2}} & i>j, \ \ i,j \text{ are even}\\
        0; & \text{ else}
    \end{cases}$$ and thus by using the definition of $x_i^*$ it is indeed the inverse. Hence $N$ is a subgroup.
\subsubsection{$\text{Sol}_{A_n}$ as a semi-direct product}
\begin{lemma}\label{connected_lemma}
    The constructed subgroup $N$ is contractible to a point, in particular, it is path-connected. 
\end{lemma}
\begin{proof}
    We will show that there is a continuous path from any element $X \in N$ to the identity element. Let us take an arbitrary element $X \in N$ and write it as $I + Y$, separating the diagonal elements and the rest. Observe that the columns of $Y$ are given as in the Lemma \ref{diagonal}, except that zeros have replaced the diagonal 1s. We define the matrix $Y_t$ by replacing the variables $x_i$ with $tx_i$ in $Y$. It is easy to check that $X_t= I+Y_t$ is an element of $Sol_{A_n}$. This gives a straight-line homotopy that shows that $N$ is contractible to $I$.
\end{proof}
\begin{theorem}
    $N$ is a normal subgroup of $\text{Sol}_{A_n}$. Furthermore, $\text{Sol}_{A_n} \cong N \rtimes D.$
\end{theorem}
\begin{proof}
    Consider the map $Sol_{A_n} \to D$, which evaluates all the non-diagonal entries to zero; by definition, the fiber at the identity is the subgroup $N$. Now note that $N$ has the Lie algebra as $Lie(N) = \mathfrak{a}$ an ideal in the Lie algebra $sol_{A_n}$; this proves that $N$ is a normal subgroup. Also, note that $Lie (Sol_{A_n}/N) \simeq Lie (D)$ and $N$ is connected from Lemma \ref{connected_lemma}, we conclude that $Sol_{A_n}/N \simeq D$.

\end{proof}
\begin{remark}
    One can see that the groups $Sol_{A_n}$ are abelian when $n$ is even and are non-abelian when $n$ is odd. This is because in the exact sequence as below, \[ 1 \to \mathbb{C}^{\lceil\frac{n-2}{2}\rceil} \to Sol_{A_n} \to \mathbb{C}^ *\to 1\] when $n$ is even the exact sequence splits due to the trivial action of the $\mathbb{C}^*$, but it does not split when $n$ is odd.
\end{remark}
\begin{corollary}
    The fundamental group of $Sol_{A_n}$ is $\mathbb{Z}$, and $\text{dim}(Sol_{A_n}) = \lceil\frac{n-2}{2}\rceil + 1.$
\end{corollary}
\begin{proof}
    Note that $N$ has been proved to be contractible, and $Sol_{A}$ is a semi-direct product of $\mathbb{C}^*$ and $N$, which gives us the result.
\end{proof}
\section{The structure of $Sol_{A^2_n}$ with $n \geq 4$}\label{Sec5}
Recall, $A_n^2 = {\tiny \begin{pmatrix}
0 &0  & 1 & &&0\\ 
 & 0 & 0 &1&&\\ 
 &  &  \ddots&\ddots &\ddots&\\ 
& &&\ddots &\ddots&1\\ 
 &  &  & &\ddots&0\\
 0& &&& &0
\end{pmatrix}.}$ Using corollary \ref{group_nilpotent} for $k =2$, we know that $Sol_{A^2_n} = \{X \in M_n(\mathbb{C}) \ | \ X^tA^2_nX = A^2_n\}$ is a group for all $n \geq 4$. In this section, we will find the structure of $Sol_{A^2_n}$ again by utilizing the structure of its Lie algebra.\\

For this case, we would like to break down the matrices into two by two block matrix forms (after attaching a column and a row of zeros for odd $n$). This way, the problem will reduced to the case of $A_n$ albeit with non-commuting variables. This approach vastly simplifies the computations and also makes the results more compact. Therefore, we first describe the process of reducing to the block matrix form as follows.
\begin{enumerate}
    \item \textbf{When $n$ is even.} 
    {\tiny$$\left(\begin{array}{cccc}
   X^t_{1,1}  & X^t_{3,1} &\cdots & X^t_{n-1, 1} \\
    X^t_{1,3} & X^t_{3,3} & \cdots&X^t_{n-1, 3}\\
    \vdots &\vdots &\vdots &\vdots\\
    X^t_{1, n-1} & X^t_{3, n-1} &\cdots &X^t_{n-1,n-1}
\end{array}\right)\left(\begin{array}{ccccc}
   O_2  & I_2 &O_2 &\cdots&O_2 \\
    \vdots&\ddots&\ddots &\ddots&\vdots \\
    \vdots&\ddots&\ddots &\ddots&\vdots \\
    O_2&\ddots&\ddots&O_2&I_2\\
    O_2&\cdots &\cdots& O_2&O_2
\end{array}\right)\left(\begin{array}{cccc}
   X_{1,1}  & X_{1,3} &\cdots & X_{1, n-1} \\
    X_{3,1} & X_{3,3} & \cdots&X_{3, n-1}\\
    \vdots &\vdots &\vdots &\vdots\\
    X_{n-1, 1} & X_{n-1, 3} &\cdots &X_{n-1,n-1}
\end{array}\right) $$}
{\tiny $$ = \begin{pmatrix}O_2  & I_2 &O_2 &\cdots&O_2 \\
    \vdots&\ddots&\ddots &\ddots&\vdots \\
    \vdots&\ddots&\ddots &\ddots&\vdots \\
    O_2&\ddots&\ddots&O_2&I_2\\
    O_2&\cdots &\cdots& O_2&O_2
    \end{pmatrix},$$}
where $X_{i, j} = {\footnotesize \left(\begin{array}{cc}
    x_{i,j} & x_{i+1,j} \\
     x_{i, j+1}& x_{i+1, j+1}
\end{array}\right)}, \ O_2$ is the $2 \times 2$ zero matrix, and $I_2$ is the $2 \times 2$ identity matrix. Observe that with this block form, we can think of $A_n^2$ as $A_{n/2}$ views as a matrix in $M_{n/2}(M_2(\mathbb{C}))$.
\item \textbf{When $n$ is odd.} To make the decomposition into $2 \times 2$ blocks possible, we pad the matrix with a row and a column of zeros at the top of the first row and at the start of the first column, respectively. Then we proceed to break the equation as in case (1) 
$${\tiny  \left(\begin{array}{cccc}
   X^t_{11}  & Y^t_{2,1} &\cdots & Y^t_{n-1, 1} \\
    Y^t_{1,2} & X^t_{3,3} & \cdots&X^t_{n-1, 3}\\
    \vdots &\vdots &\vdots &\vdots\\
    Y^t_{1, n-1} & X^t_{3, n-1} &\cdots &X^t_{n-1,n-1}
\end{array}\right)\left(\begin{array}{ccccc}
   O_2  & a &O_2 &\cdots&O_2 \\
    \vdots&\ddots& I_2 &\ddots&\vdots \\
    \vdots&\ddots&\ddots &\ddots&\vdots \\
    O_2&\ddots&\ddots&O_2&I_2\\
    O_2&\cdots &\cdots& O_2&O_2
\end{array}\right)\left(\begin{array}{cccc}
   X_{11}  & Y_{1,2} &\cdots & Y_{1, n-1} \\
    Y_{2,1} & X_{3,3} & \cdots&X_{3, n-1}\\
    \vdots &\vdots &\vdots &\vdots\\
    Y_{n-1, 1} & X_{n-1, 3} &\cdots &X_{n-1,n-1}
\end{array}\right) }$$
 $$ = {\tiny \begin{pmatrix}O_2  & a &O_2 &\cdots&O_2 \\
    \vdots&\ddots&I_2 &\ddots&\vdots \\
    \vdots&\ddots&\ddots &\ddots&\vdots \\
    O_2&\ddots&\ddots&O_2&I_2\\
    O_2&\cdots &\cdots& O_2&O_2
    \end{pmatrix}},$$
where $O_2$, $I_2$, and $X_{i, j}$ as above except $X_{11} = {\footnotesize\left(\begin{array}{cc}
         0& 0 \\
         0&x_{1,1} 
    \end{array}\right)};$ whereas $Y_{1,i} = {\footnotesize\left(\begin{array}{cc}
         0& 0 \\
         x_{1,i} &x_{1,i+1} 
    \end{array}\right)}; \ Y_{i,1} = {\footnotesize\left(\begin{array}{cc}
         0& x_{i, 1} \\
         0 &x_{i+1, 1} 
    \end{array}\right)}$ and $a= {\footnotesize\left(\begin{array}{cc}
         0& 0 \\
         0 & 1 
    \end{array}\right)}$.
Let us call the block matrices with the zeros added as $\tilde{A}_{n+1}$ and for a solution $\tilde{X}$; one can check that the South-East $n \times n$ block of $\tilde{X}$ gives a solution of the original system. Conversely, if we have a solution $X$, of the system $X^tA_{n}^2X=A_{n}^2$, then the padded matrix $\tilde{X}$, satisfies $\tilde{X}^t \tilde{A}_{n+1} \tilde{X} = \tilde{A}_{n+1}$. As a result, it is equivalent to discussing the padded system. Furthermore, as there is no scope for confusion, to reduce the burden of notation, we will continue to denote the $\tilde{A}_{n+1}$ as $A^2_{n}$ and assume there is a $2 \times 2$ decomposition. Note that this amounts to the case that for odd $n$ the size of the matrix $A_{n}^2$ is precisely $\frac{n+1}{2}$ in $M_{\frac{n+1}{2}}( M_2(\mathbb{C}))$. 
\end{enumerate}

\subsection{Lie algebra of $\text{Sol}_{A^2_n}$}\label{lie_algebra}

    Let us denote the Lie algebra of $\text{Sol}_{A^2_n}$ by $sol_{A_n^2}$. Then $sol_{A_n^2}= \{X \in M_{n}(\mathbb{C})\ \mid \ X^tA^2_n + A_n^2X = 0\}.$ We will find the structure of $sol_{A_n^2}$ for four different cases of $n$ modulo 4.
    
   \subsubsection{Case 1: \: $\bf{n \equiv 0 (\text{mod } 4)}; \ n \geq 4$}

\begin{proposition}\label{generation_one} The set of solutions of the linear equation $ X^tA_n^2 + A_n^2X = 0$ is given by the following: For $j=1+4k$, where $0 \leq k \leq \frac{n-4}{4}$, put $X_{1, j} = {\footnotesize\left(\begin{array}{cc}
    x_{1,j} & x_{1,j+1} \\
     x_{2, j}& x_{2, j+1}
\end{array}\right)},$ then the structure of the Lie algebra $sol_{A^2_n}$ is given by
{\tiny $$\left(\begin{array}{ccccccc}
        X_{1,1} &O_2&X_{1,5}&O_2 &\cdots & X_{1,n-3}& O_2 \\
    O_2& -X_{1,1}^t&\ddots& O_2&\ddots&\cdots&O_2\\
    O_2 &\ddots &\ddots&\ddots&\ddots&\ddots&\vdots\\
    \vdots&-X_{1,5}^t&\ddots&\ddots&\ddots&X_{1,5}&\vdots\\
    \cdots &\ddots &\ddots&\ddots&\ddots&\ddots&O_2\\
    O_2&\cdots&\ddots &O_2&\ddots&X_{1,1}&O_2\\
     O_2& -X_{1,n-3}^t &\cdots& O_2&-X_{1,5}^t&O_2&-X_{1,1}^t
    \end{array}\right).$$}

\end{proposition}
\begin{proof}
    It is clear from a direct calculation that these are solutions of the equation $X^tA_n^2+A_n^2X=0$. Conversely, we will show that these matrices generate a general solution. Let $X$ be a general solution matrix. Let us call the block matrices for $A_n^2$ to be $B$ and that of $X$ to be $Y$, so the matrices $B$ and $Y$ are in $M_{n/2}(M_2(\mathbb{C}))$. Let us also use the notation $M_{i,j}$ to denote the $ij$-th entry of a matrix $M$ and $R_i(M)$, $C_i(M)$ to denote the $i$-th row and column of $M$ respectively.  Observe that $BY$ is a matrix with $R_i(BY)= R_{i+1}(Y)$ for $i \leq n-1$, and $R_n(Y)= 0$, a zero row vector. Similarly, $C_i(Y^tB)= C_{i-1}(Y^t)$ for $i > 1$, and $C_1(Y^tB)= 0$, a zero column vector. Putting these two together, we have $(BY)_{i,j}=Y_{i+1,j}$ whenever $i \leq n-1$, and zero otherwise. Similarly $(Y^tB)_{ij}=Y_{j-1,i}^t$ whenever $j >1$, and zero otherwise. Then finally we have, $ (Y^tB+BY)_{i,j}= Y_{i+1,j } + Y_{j-1,i}^t $ whenever $i \leq n-1$ and $j > 1$, $(Y^tB+BY)_{n,j}= Y_{j-1, n}^t$ for $j > 1$, $(Y^tB+BY)_{i,1}= Y_{i+1,1}$ for $i \leq n-1$;, and $(Y^tB+BY)_{n,1}=0$. Setting $Y^tB+BY=0$, we obtain the equations:
   \begin{equation}\label{equation1}
   \begin{split}
    Y_{i+1,j } + Y_{j-1,i}^t=0 & \mbox{ whenever } i \leq n-1 \mbox{ and }j > 1 \\
    Y_{j-1, n}^t=0 &\mbox{ for } j > 1\\
     Y_{i+1,1}=0 & \mbox{ for } i \leq n-1.
\end{split}
    \end{equation}
    Note that the first equation of (\ref{equation1}) can be written as $Y_{i,j}=-Y^t_{j-1,i-1}$ whenever $ 2 \leq i \leq n-1$ and $j >1$.
Let us fix the entries $Y_{1,j}$ for $1 \leq j \leq n$ and show that all other entries can be written in terms of these. 
\begin{itemize}
    \item \textbf{Case I} $i <j$ and $i$ is even. In this case, we have $Y_{i,j}=-Y^t_{j-1,i-1}= Y_{i-2, j-2} = \cdots$, at the end we will stop after an odd number of equalities with $-Y^t_{j-i+1,1}$, and using the third equation from (\ref{equation1}) this is equal to zero. Hence in this case $Y_{i,j}=0$
    \item \textbf{Case II} $i < j$ and $i$ is odd. We will stop after an even number of steps with $Y_{1,j-i+1}$.
    \item \textbf{Case III} $ i > j$. In this case, we write $Y_{i,j}= - Y^t_{j+1, i+1}= Y_{i+2, j+2} = \cdots$, we will end with $Y_{n, j+ n-i}$ if $n-i$ is even; and if $n-i$ is odd, we end with $-Y_{j + n-i,n}=0$ from the second equation of (\ref{equation1}). 
\end{itemize}
    As a result of the above, the only possibly nonzero elements are of the form $Y_{i,j}=Y_{1,j-i+1}$, or $Y_{i,j}=Y_{n,j+n-i}$. So far, we have shown that all the entries are either zero or equal to some entries or its negative transpose on the first or the last row.  Now, we will prove the following claims, clarifying the elements of the first and the last row.
    \begin{enumerate}
        \item $Y_{1,l} =0$, whenever $l$ is even. 
        \item $Y_{n, l}=0 $, whenever $l$ is odd.
        \item $Y_{1,l}+ Y_{n, n+1-l}^t=0$, for $1 \leq l \leq n$.
    \end{enumerate}
    Once we prove the above three claims, the proof of the proposition is complete. Now to prove the above claims. 
    Note that one can rewrite the first statement of (\ref{equation1}) in the following form,
    \[ Y_{i,j}= - Y_{j+1,i+1}^t.\]
    \begin{enumerate}
        \item For $l$ even, we can write $Y_{1,l}= -Y_{l+1,2}^t= Y_{3,l+2} =\cdots$, this stops with $Y_{n-l+1, n}$, but as $Y_{j,n}=0$ for $1 \leq j \leq n-1$, we have the result.
        \item For $l$ odd, one can write $Y_{n,l}= -Y_{l-1,n-1}^t= Y_{n-2, l-2} = \cdots $, this stops with $Y_{n-l+1, 1}=0$. This proves the point.
        \item Let us start with $Y_{1,l}$ and write it as $Y_{1,l}= - Y_{l+1,2}^t = Y_{3,l+2} =\cdots$, this will stop with $Y_{1,l}= -Y_{n, n+1-l}^t $ when $l$ is odd, thus we have $Y_{1,l}+ Y_{n, n+1-l}^t=0$ ; and for $l$ even we get $Y_{1,l}= -Y_{n, n+1-l}$ but using above two points both quantities are zero, thus it proves the statement.
    \end{enumerate}

\end{proof}

We will use below elementary results about the Lie brackets of elementary matrices while proving the following proposition; proofs of these are routine and are omitted.
\begin{lemma}\label{zero} Let $E_{ij}$ denote the elementary matrix with 1 at the $ij$-th place and zero elsewhere. 
\begin{itemize}
    \item  If $S_1$ and $S_2$ are two distinct sets of indices. Then $\left[\sum_{i \in S_1} E_{i,i}, \sum_{i \in S_2} E_{i,i}\right]=0.$
    \item The Lie bracket $[E_{ij}, E_{lk}]= E_{ik} \delta_{jl} - E_{lj} \delta_{ik},$ where $\delta_{ij}$ is the Kronecker delta symbol.
    \item  The Lie bracket $\left[\sum_{i \in S} E_{i,i}, E_{l,k}\right]= \epsilon_{l, S} E_{l,k} - \epsilon_{k,S} E_{l,k} $, where $\epsilon_{x, S}=1$ if $x \in S$ and zero otherwise.
    \item   Let $S$ be a set of indices and $m$ be an index. Then $$\left[ \sum_{i \in S} E_{i,i} , \sum_{i \in S_m}E_{i, i+m}\right]=0,$$ where $S_m$ is a subset of indices such that $i+m \in S$. 
    \item Let $a = \sum_{k,l} E_{k,l}$ and $b= \sum_{s,t} E_{s,t} $. If $k, s$ runs over even numbers (odd numbers) and $l, t$ runs over odd numbers (even numbers), then $[a,b]=0$.
\end{itemize}
   
\end{lemma}

    \begin{proposition}\label{brackets_one}
        The Lie algebra $sol_{A_n^2}$ is generated by the following elements:
        $$h_1 = \sum_{l =0}^{\tfrac{n-4}{4}}\left(E_{1+4l,1+4l} - E_{3+4l,3+4l}\right); \ \ \ \  h_2 =  \sum_{l =0}^{\tfrac{n-4}{4}}\left(E_{2+4l,2+4l} - E_{4+4l,4+4l}\right);$$
        $$e = \sum_{l =0}^{\tfrac{n-4}{4}}\left(E_{1+4l,2+4l} - E_{4+4l,3+4l}\right); \ \ \ f =\sum_{l =0}^{\tfrac{n-4}{4}}\left(E_{2+4l,1+4l} - E_{3+4l,4+4l}\right);$$
        
        For $k =1, 2, \ldots, \frac{(n-4)}{4}$,
         $$a_k = \sum_{l=0}^{\tfrac{n-4-4k}{4}}\left(E_{1+4l,1+4l+4k} - E_{3+4l+4k,3+4l}\right); \ \ \ b_k=  \sum_{l=0}^{\tfrac{n-4-4k}{4}}\left(E_{1+4l,2+4l+4k} - E_{4+4l+4k,3+4l}\right);$$
    $$c_k = \sum_{l=0}^{\tfrac{n-4-4k}{4}}\left(E_{2+4l,1+4l+4k} - E_{3+4l+4k,4+4l}\right); \ \ d_k = \sum_{l=0}^{\tfrac{n-4-4k}{4}}\left(E_{2+4l,2+4l+4k} - E_{4+4l+4k,4+4l}\right),$$
    with the following Lie bracket relations $[X,Y]$ where $X$ is from the vertical margin and $Y$ is from the horizontal margin:
         
        \begin{center}
        {\tiny \begin{tabular}{c|cccccccc}
        & $h_1$ & $h_2$ & $e$ & $f$ & $a_{k'}$ & $b_{k'}$ & $c_{k'}$ & $d_{k'}$ \\
        \hline
        $h_1$ & $0$ & $0$ & $e$ & $-f$ & $0$ & $b_{k'}$ & $-c_{k'}$&$0$\\
        $h_2$ & $0$ & $0$& $-e$& $f$ & $0$ & $-b_{k'}$ & $c_{k'}$ & $0$\\
        $e$ & $-e$ & $e$ & $0$ & $h_1 - h_2$ & $-b_{k'}$ & $0$ & $a_{k'}-d_{k'}$ & $b_{k'}$ \\
        $f$ & $f$ & $-f$ & $-(h_1 - h_2)$ & $0$ & $c_{k'}$ & $-(a_{k'}-d_{k'})$ & $0$&$-c_{k'}$ \\
        $a_k$ & $0$ & $0$ & $b_k$ & $-c_{k}$ & $0$ & $b_{k+k'}$ & $-c_{k+k'}$ & $0$\\
        $b_k$ & $-b_k$ & $b_k$ & $0$ & $a_{k}-d_{k}$ & $0$ & $0$ & $a_{k+k'}-d_{k+k'}$ & $b_{k+k'}$ \\
        $c_k$ & $c_k$ & $-c_{k}$ & $a_k-d_k$ & $0$ & $-c_{k+k'}$ & $-(a_{k+k'}-d_{k+k'})$ & $0$ & $-c_{k+k'}$ \\
        $d_k$ & $0$ & $0$ & $-b_k$ &$-c_k$& $0$ & $-b_{k+k'}$ & $c_{k+k'}$ & $0$ \\
       \end{tabular}}
       \end{center}
    \end{proposition}
\begin{proof}

We proved the generation in the Proposition \ref{generation_one} above. Now, we sketch the proof of the Lie brackets as below.\\    
 The $[h_1,h_2]$ statement follows from the Lemma \ref{zero} as we have distinct sets of indices. For $[h_1,e]$, note that one can write $h_1= E_{1+4l, 1+4l}- E_{3+4l,3+4l}$ ignoring the summation, and $e= E_{1+4l, 2+4l}- E_{4+4l, 3+4l}$, ignoring the summation over $l$ running on suitable range. Now the Lemma \ref{zero} gives the required results. Similarly write $a_k$ as $E_{1+4l, 1+4l+4k}- E_{3+4l+4k, 3+4l}$ ignoring the summation over appropriate range. As a result of this $[h_1,a_k]= a_k$ and $[h_2,a_k]=0$ by applying the Lemma \ref{zero}. 

  $[h_2,e]= [E_{2+4l,2+4l}- E_{4+4l,4+4l}, E_{1+4l',2+4l'}-E_{4+4l',3+4l'}]$ suppressing the summation over the variables $l,l'$ in their appropriate range. Simple calculation of this gives us $[h_2,e]= E_{4+4l',3+4l'}-E_{1+4l',2+4l'}=-e$. A similar calculation gives $[h_2,f]=f$. 
  For the bracket $[e,f]$ observe, $[E_{1+4l,2+4l},E_{2+4l,1+4l} ] = h_1 $ by using the Lemma \ref{zero}. Similarly, one can observe that $ [E_{4+4l,3+4l},E_{3+4l,4+4l}]= -h_2 $ with the rest of the brackets being zero; this completes the assertion. Lemma \ref{zero} allows us to write $[e,a_k]= [E_{1+4l,2+4l} - E_{4+4l,3+4l},E_{1+4l',1+4l'+4k} - E_{3+4l'+4k,3+4l'} ]$ $ = - E_{1+4l', 1+4l'+4k} E_{1+4l, 2+4l}+ E_{4+4l, 3+4l}E_{3+4l'+4k, 3+4l'}$. As all other products are zero. Now allowing $1+4l'+4k=1+4l$ and $3+4l=3+4l'+4k$ respectively in the two products followed by rearranging the sum, we obtain $[e,a_k]= - E_{1+4l, 2+4l+4k}+E_{4+4l+4k,3+4l}=-b_k $. The brackets with $b_k, c_k, d_k$ are verified using a similar technique.
  For the bracket $[f,b_k]$, withholding the summation, we can write it as the Lie bracket, $[E_{2+4l,1+4l} - E_{3+4l,4+4l}$, $E_{1+4l',2+4l'+4k} - E_{4+4l'+4k,3+4l'}]$, which is further equal to $[E_{2+4l,1+4l},E_{1+4l',2+4l'+4k} ]$ $ +[E_{3+4l,4+4l}, E_{4+4l'+4k,3+4l'} ] $ because all other brackets will be equal to zero by the Lemma \ref{zero}. Now $[E_{2+4l,1+4l},E_{1+4l',2+4l'+4k} ]=E_{2+4l,1+4l}E_{1+4l', 2+4l'+4k} -E_{1+4l', 2+4l'+4k}E_{2+4l,1+4l}$ $=E_{2+4l, 2+4l+4k}-E_{1+4l, 1+4l+4k}$ after rearranging the sum and applying Lemma \ref{zero}. For $[E_{3+4l,4+4l}, E_{4+4l'+4k,3+4l'} ]$ a similar calculation and gathering the terms together gives $[f,b_k]= d_k-a_k$.
  The rest of the Lie bracket statements also follow easily with such calculations.
\end{proof}

\begin{theorem}
    With the notations as above and for a fixed $k$, let $\mathfrak{b}_k$ be the span of $a_k,b_k,c_k,d_k$, and put $\mathfrak{a} =\bigoplus_{k} \mathfrak{b}_k$. Then the following are true:
\begin{enumerate}
     \item For any $k$, $\mathfrak{b}_k$ is a Lie ideal in $sol_{A_n^2}$.  
    \item $[\mathfrak{b}_k , \mathfrak{b}_k] \subseteq \mathfrak{b}_{k+1}$.
    \item The ideal $\mathfrak{a}$ is a nilpotent ideal of $sol_{A_n^2}$.
    \item The Lie algebra spanned by $h_1, h_2, e, f$ is isomorphic to $sl_2 \oplus \mathbb{C}.$
\end{enumerate}
\end{theorem}

\begin{proof} These assertions follow from the table of Lie brackets in Proposition \ref{brackets_one},
    \begin{enumerate}
        \item immediately follows from the table.
        \item follows because for all $X,Y \in \mathfrak{b}_k$, clearly $[X,Y]$ is generated by $a_{k+1}, b_{k+1}, c_{k+1}, d_{k+1}$ and hence lands in $\mathfrak{b}_{k+1}.$
        \item This follows from (2).
        \item Observe that the Lie algebra spanned by $h_1-h_2, e, f$ is a copy of $sl_2$, and the supplementary subspace generated by the element $h_1+h_2$, which gives a copy of $\mathbb{C}$ in $sol_{A_n^2}/\mathfrak{a}$. 
        \end{enumerate}
       
\end{proof}
As a consequence of above theorem, we get the following corollaries:
\begin{corollary}\label{radical_0}
    $\mathfrak{a} \oplus \langle h_1+h_2 \rangle$ is the solvable radical of $sol_{A_n^2}$, and we have the exact sequence,
    \[ 0 \to \mathfrak{a}\oplus \mathbb{C}(h_1+h_2) \to sol_{A_n^2} \to sl_2 \to 0\]
\end{corollary}
By exponentiating the above exact sequence, we get the following analogous sequence in the groups:

\begin{corollary}
Let $K$ be the maximal solvable subgroup of $Sol_{A_n^2}$. Then
\begin{itemize}
    \item $ 1 \to K \to (Sol_{A_n^2})^\circ \to SL_2 \to 1$
 \item The dimension of the group $Sol_{A^2_n}$ is $n$.
\end{itemize}
\end{corollary}
 \subsubsection{Case 2: \: $\bf{n \equiv 2 (\text{mod } 4)}; \ n \geq 6$}

 \begin{proposition}\label{generation_two} 
 For $j=3+4k$, where $0 \leq k \leq \frac{n-6}{4}$, put $X_{1, j} = {\footnotesize\left(\begin{array}{cc}
    x_{1,j} & x_{1,j+1} \\
     x_{2, j}& x_{2, j+1}
\end{array}\right)}.$ Then the structure of the solution of the equation $X^tA_{n}^2 + A_{n}^2X = 0$ in this case is given by,\\
     {\tiny  $$\left(\begin{array}{cccccccc}
        X_{1,1} &X_{1,3}&O_2&X_{1,7}&O_2 &\cdots & X_{1,n-3}& O_2 \\
    O_2& -X_{1,1}^t&O_2&\ddots& O_2&\ddots&\cdots&O_2\\
    O_2 &-X_{1,3}^t&\ddots&\ddots&\ddots&\ddots&\ddots&\cdots\\
    O_2& \ddots&\ddots&\ddots&\ddots&\ddots&X_{1,7}&O_2\\
    O_2&-X_{1,7}^t&\ddots&\ddots&\ddots&\ddots&\ddots&O_2\\
    \cdots &\ddots &\ddots&\ddots&\ddots&X_{1,1}&X_{1,3}&O_2\\
    O_2&\cdots&\ddots &O_2&\ddots&O_2&-X_{1,1}^t&O_2\\
     O_2& -X_{1,n-3}^t &\cdots& O_2&-X_{1,7}^t&O_2&-X_{1,3}^t&X_{1,1}
    \end{array}\right),$$}
    where $X_{1,1} = {\footnotesize\left(\begin{array}{cc}
    x_{1,1} & x_{1,2} \\
     x_{2, 1}& x_{2, 2}
\end{array}\right)}$.
  \end{proposition}
 \begin{proof}
     The proof of this case is exactly similar to the even case $n \equiv 0 \mbox{ (mod 4)}$, so it is skipped.
 \end{proof}

    \begin{proposition}\label{brackets_two}
        In this case $sol_{A^2_n}$ is generated by the following elements:\\
        $$h_1 = \sum_{l=0}^{\tfrac{n-2}{4}}\left(E_{1+4l,1+4l} - E_{3+4l,3+4l}\right); \ \ \ h_2 =  \sum_{l=0}^{\tfrac{n-2}{4}}\left(E_{2+4l,2+4l} - E_{4+4l,4+4l}\right);$$
        $$e = \sum_{l=0}^{\tfrac{n-2}{4}}\left(E_{1+4l,2+4l} - E_{4+4l,3+4l}\right); \ \ \ f =\sum_{l=0}^{\tfrac{n-2}{4}}\left(E_{2+4l,1+4l} - E_{3+4l,4+4l}\right);$$
        \bigskip
        
        For $k = 1, 2, \ldots, \frac{(n-2)}{4}$,
        $$a_k = \sum_{l=0}^{\tfrac{n-2-4k}{4}}\left(E_{1+4l,-1+4l+4k} - E_{1+4l+4k,3+4l}\right); \ \ b_k=  \sum_{l=0}^{\tfrac{n-2-4k}{4}}\left(E_{1+4l,4l+4k} - E_{2+4l+4k,3+4l}\right);$$
    $$c_k = \sum_{l=0}^{\tfrac{n-2-4k}{4}}\left(E_{2+4l,-1+4l+4k} - E_{1+4l+4k,4+4l}\right); \ \ d_k = \sum_{l=0}^{\tfrac{n-2-4k}{4}}\left(E_{2+4l,4l+4k} - E_{2+4l+4k,4+4l}\right),$$
      with the following Lie bracket relations:
         
        \begin{center}
        {\tiny \begin{tabular}{c|cccccccc}
        & $h_1$ & $h_2$ & $e$ & $f$ & $a_{k'}$ & $b_{k'}$ & $c_{k'}$ & $d_{k'}$ \\
        \hline
        $h_1$ & $0$ & $0$ & $e$ & $-f$ & $2a_{k'}$ & $b_{k'}$ & $c_{k'}$&$0$\\
        $h_2$ & $0$ & $0$& $-e$& $f$ & $0$ & $b_{k'}$ & $c_{k'}$ & $2d_{k'}$\\
        $e$ & $-e$ & $e$ & $0$ & $h_1 - h_2$ & $0$ & $a_{k'}$ & $a_{k'}$ & $b_{k'} + c_{k'}$ \\
        $f$ & $f$ & $-f$ & $-(h_1 - h_2)$ & $0$ & $b_{k'}+c_{k'}$ & $d_{k'}$ & $d_{k'}$ &$0$ \\
        $a_k$ & $-2a_{k}$ & $0$ & $0$ & $-(b_k+ c_{k})$ & $0$ & $0$ & $0$ & $0$\\
        $b_k$ & $-b_k$ & $-b_k$ & $-d_{k}$ & $-d_{k}$ & $0$ & $0$ & $0$ & $0$ \\
        $c_k$ & $-c_k$ & $-c_{k}$ & $-a_k$ & $-d_k$ & $0$ & $0$ & $0$ & $0$ \\
        $d_k$ & $0$ & $-2d_k$ & $-(b_k+ c_{k})$ &$0$& $0$ & $0$ & $0$ & $0$ \\
       \end{tabular}}
       \end{center}

    \end{proposition}
   \begin{proof}
       This proof is similar to the proof of the Proposition \ref{brackets_one}, so we can skip the proof.
   \end{proof}
   \begin{theorem}\label{ideal-two}
       Let us denote $\mathfrak{b}_l$ to be the span of the elements $a_l,b_l,c_l,d_l$, and $\mathfrak{a}_k = \bigoplus_{l \geq k} \mathfrak{b}_l$, we have,

       \begin{enumerate}
           \item For any $k$, $\mathfrak{b}_k$, $\mathfrak{a}_k$ are Lie ideals.
           \item $[\mathfrak{a}_k , \mathfrak{a}_k] \subseteq \mathfrak{a}_{k+1}$.
           \item $\mathfrak{a}=\mathfrak{a}_1$ is an abelian Lie ideal of $Sol_{A^2_n}$, and thus it is also nilpotent.
           \item The Lie algebra generated by $h_1, h_2, e,f$ is $sl_2 \oplus \mathbb{C}$.
       \end{enumerate}
   \end{theorem}
   \begin{proof}
\begin{enumerate}
    \item Note that for any $X \in Sol_{A_n^2}$ and $Y \in \mathfrak{b}_k$, we have $[X,Y]$ generated by $a_k, b_k, c_k, d_k$. So $[X,Y] \in \mathfrak{b}_k$, implies $\mathfrak{b}_k$ is an ideal. As a consequence, the direct sum $\mathfrak{a}_k = \oplus_{ l \geq k } \mathfrak{b}_l$, is an ideal.
    \item It follows from the observation that $[\mathfrak{b}_k, \mathfrak{b_k}]=0$.
    \item One can observe from the table of Proposition \ref{brackets_two} that for any $X,Y \in \mathfrak{b}_k$, $[X,Y]=0$ as a result, the ideal $\mathfrak{b}_k$ is abelian. Thus the direct sum ideal $\mathfrak{a}$ is abelian. 
    \item The proof is clear from the table in the Proposition \ref{brackets_two}, and observation that $h=h_1-h_2, e, f$ is a copy of $sl_2$.
\end{enumerate}       
   \end{proof}
   As a result, we have the following corollaries:
   \begin{corollary}\label{radical_2}
       We have an exact sequence of Lie algebras,
       \[0 \to \mathfrak{a} \oplus \mathbb{C} \to sol_{A^2_n} \to sl_2 \to 0\]
   \end{corollary}
\begin{proof}
Clearly, as the ideal $\mathfrak{a}\oplus \mathbb{C}(h_1+h_2)$ is the solvable radical of the algebra $\mathfrak{g}$ and the quotient is generated by $h_1-h_2, e, f$.
\end{proof}
   \begin{corollary}
   \begin{enumerate}
       \item At the group level, we have 
       \[ 1 \to K \to (Sol_{A_n^2})^\circ \to SL_2 \to 1\]
 where $K$ is the largest solvable subgroup of the group $Sol_{A^2_n}$, and it is a semidirect product of an abelian subgroup and a one-dimensional torus. 
 \item The dimension of the group $Sol_{A^2_n}$ is $n+2$.
   \end{enumerate}
       \end{corollary}
\subsubsection{Case 3: \:$\bf{n \equiv 1 (\text{mod } 4)}, \ \ n \geq 5$}

\begin{proposition}\label{generation_three}
    The structure of the solution of the equation $X^tA^2_n + A^2_nX = 0$ in this case is given by,\\
      $${\tiny \left(\begin{array}{cccccccc}
        X_{11} &Y_{1,2}&O_2&Y_{1,6}&O_2 &\cdots & Y_{1,n-3}& O_2 \\
    O_2& X_{2,2}&O_2&Y_{2,6}& O_2&\ddots&Y_{2,n-3}&O_2\\
    O_2 &-Y_{1,2}^t&-X_{2,2}^t&\ddots&O_2&\ddots&\ddots&O_2\\
    O_2& O_2&\ddots&\ddots&\ddots&\ddots&Y_{1,6}&O_2\\
    O_2&-Y_{1,6}^t&-Y_{2,6}^t&\ddots&\ddots&\ddots&Y_{2,6}&O_2\\
    \vdots &\vdots&\ddots&\ddots&\ddots&\ddots&Y_{1,2}&O_2\\
    O_2&O_2&\ddots &O_2&O_2&O_2&X_{2,2}&O_2\\
     O_2& -Y_{1,n-3}^t &-Y_{2,n-3}^t& \cdots&-Y_{1,6}^t&-Y_{2,6}^t&-Y_{1,2}^t&-X_{2,2}^t
    \end{array}\right)},$$
    
    where  $X_{11} = {\footnotesize\left(\begin{array}{cc}
         0& 0 \\
         0&x_{1,1} 
    \end{array}\right)}; \ Y_{1,i} = {\footnotesize\left(\begin{array}{cc}
         0& 0 \\
         x_{1,i} &x_{1,i+1} 
    \end{array}\right)}$ for $i = 2+4k, \ 0 \leq k \leq \tfrac{n-5}{4}$; $Y_{2,j} = {\footnotesize\left(\begin{array}{cc}
         x_{2,j} &x_{2,j+1} \\
         0&0
    \end{array}\right)}$ for $j = 6+4k, \ 0 \leq k\leq \tfrac{n-9}{4}$; and \ $X_{2,2} = {\footnotesize\left(\begin{array}{cc}
         x_{2,2}& x_{2,3} \\
         0&-x_{1,1} 
    \end{array}\right)}$. Add the first row and first column with an appropriate number of zeros.  
    
\end{proposition}
\begin{proof}
    One can check these are indeed the solutions of the equation $X^tA_2^n +A_2^nX = 0$. Conversely, we will show that these matrices generate a general solution. Let $X$ be an $n \times n$ general solution matrix, let $Y$ be the block matrix $Y = {\footnotesize\left(\begin{array}{cc}
         0& 0 \\
         0& X
    \end{array}\right)}$, is the $(n+1) \times (n+1) $ matrix with a row of zeros added at the beginning and a column of zeros added at the beginning. Similarly, we obtain $B$ by padding the matrix $A^2_n$. One can observe that $X^t A^2_n+A^2_nX=0$ if and only if $Y^t B + BY=0$. Now convert both matrices $Y$ and $B$ into $2 \times 2$ blocks to get the matrices of size $m :=\tfrac{(n+1)}{2}$ as follows: 
  
    $Y_{1,1} = {\footnotesize\left(\begin{array}{cc}
         0& 0 \\
         0& x_{1,1}
    \end{array}\right)}$, $Y_{i,1} = {\footnotesize\left(\begin{array}{cc}
         0& x_{2i-2, 2j-1} \\
         0& x_{2i-1, 2j-1}
    \end{array}\right)} \text{ for all } i\geq 2$, $Y_{1,j} = {\footnotesize\left(\begin{array}{cc}
         0& 0 \\
         x_{2i-1, 2j-2}& x_{2i-1, 2j-1}
    \end{array}\right)}$ for all $j \geq 2,$ and $Y_{i,j} = {\footnotesize\left(\begin{array}{cc}
         x_{2i-2, 2j-2}& x_{2i-2, 2j-1} \\
         x_{2i-1, 2j-2}& x_{2i-1, 2j-1}
    \end{array}\right)}$ for all $i,j \geq 2,$ where $Y_{i,j}$ denotes the $ij$-entry of $Y$. Whereas, $B_{i,j} = \begin{cases}
        a = {\footnotesize\left(\begin{array}{cc}
         0& 0 \\
         0& 1
    \end{array}\right)}, & i=1, j=2;\\
    I_2, & j = i+1, i\geq 2;\\
    O_2, \text{ otherwise}.
    \end{cases}$\\
   With these conversations, the relation $Y^tB +BY = 0$, implies the following:
   \begin{enumerate}
       \item $Y_{j-1, i}^t + Y_{i+1, j} = 0$ for $ 1 < i, j < m$
       \item $Y_{i, m}^t = 0$ for all $1<i<m$
       \item $Y_{i, 1}=0$ for $3\leq i \leq m$
       \item $Y_{i,1}^t + aY_{2,i+1} = 0$ for $1 \leq i<m$
       \item $Y_{1, j}^ta = Y_{1,j}$ for all $j$.
       \item $Y_{i,j}= Y_{i+1,j+2}$
       
   \end{enumerate}
   As a result of (3) and (4), we have $Y_{2, j} = {\footnotesize\left(\begin{array}{cc}
         0& 0 \\
         \ast& \ast
    \end{array}\right)}$, and  $Y_{2,2} = {\footnotesize\left(\begin{array}{cc}
         \ast& \ast \\
         0& -x_{1,1}
    \end{array}\right)}$. Since $Y_{i,i}=Y_{i+2,i+2}$ and $Y_{i,i}= -Y_{i+1,i+1}^t$, $i >1$, we have the desired structure for the diagonal blocks. To tackle the non-diagonal entries,
    let us build another matrix $W$ of size $m \times m$ as follows:
   $W_{i,i}= 0 $ for all $1 \leq i \leq m$, $W_{i,j}= Y_{i,j}$ if $i < j$ and $W_{ij}= Y_{i+1, j+1}$. Note that the upper triangle of $W$ coincides with the upper triangle of $Y$. We observe that $W$ is skew-symmetric, indeed, if $i <j$, the $ij$th entry of $W+W^t$ is $Y_{i,j}+Y_{j+1,i+1}^t$, which is zero from the above observation. Now if $i > j$, we have the $ij$th entry of $W+W^t$ is $Y_{i+1,j+1}+Y_{j,i}^t$ which is again zero by our observation above. As a result, we can ignore the lower triangle of $W$ and analyze the relations in the upper triangle. So let us fix $i < j$, using $Y_{i,j} = Y_{i+2,j+2}$, we see that whenever $j$ is even $Y_{i,j}=0$. This is because, after a sufficient number of $l$ steps of reapplying this relation, we obtain $Y_{i,j}= Y_{i + 2l, m} =0$. Now, whenever $j$ is odd, we obtain that $Y_{i,j}= Y_{i+2,j+2}$, as a result if we take $Y_{i,1}$ as the free parameters, we establish the structure of the statement of the Proposition. This completes the proof.

\end{proof}
    
\begin{proposition}\label{brackets_three}
        $sol_{A^2_n}$ is generated by the following elements:\\
        $$h_1 = \sum_{l=0}^{\tfrac{n-1}{4}}\left(E_{1+4l,1+4l} - E_{3+4l,3+4l}\right); \ \ \ h_2 =  \sum_{l=0}^{\tfrac{n-1}{4}}\left(E_{2+4l,2+4l} - E_{4+4l,4+4l}\right);$$
        $$e = \sum_{l=0}^{\tfrac{n-1}{4}}\left(E_{1+4l,2+4l} - E_{4+4l,3+4l}\right); \ \ \ f= \sum_{l=0}^{\tfrac{n-1}{4}}\left(E_{1+4l,3+4l} - E_{5+4l,3+4l}\right);$$
        $$ g = \sum_{l=0}^{\tfrac{n-1}{4}}\left(E_{2+4l,3+4l} - E_{5+4l,4+4l}\right);$$
       \bigskip
        
        For $k = 1, 2, \ldots, \frac{(n-5)}{4}$,
        $$a_k = \sum_{l=0}^{\tfrac{n-1-4k}{4}}\left(E_{1+4l,2+4l+4k} - E_{4+4l+4k,3+4l}\right); \ \ \ b_k = \sum_{l=0}^{\tfrac{n-1-4k}{4}}\left(E_{1+4l,3+4l+4k} - E_{5+4l+4k,3+4l}\right);$$
        $$c_k = \sum_{l=0}^{\tfrac{n-1-4k}{4}}\left(E_{2+4l,2+4l+4k} - E_{4+4l+4k,4+4l}\right); \ \ \ d_k = \sum_{l=0}^{\tfrac{n-1-4k}{4}}\left(E_{2+4l,3+4l+4k} - E_{5+4l+4k,4+4l}\right).$$
        with the following Lie brackets relations:

        \begin{center}
        {\tiny \begin{tabular}{c|ccccccccc}
        & $h_1$ & $h_2$ & $e$ & $f$ & $g$ & $a_{k'}$ & $b_{k'}$ & $c_{k'}$ & $d_{k'}$ \\
        \hline
        $h_1$ & $0$ & $0$ & $e$ & $f$ & $g$ & $a_{k'}$& $b_{k'}$ & $0$ & $d_{k'}$\\
        $h_2$ & $0$ & $0$& $-e$& $0$ & $g$ & $-a_{k'}$ & $0$ & $0$ & $d_{k'}$\\
        $e$ & $-e$ & $e$ & $0$ & $0$ & $f$ & $0$ & $0$ & $a_{k'}$ & $b_{k'}$\\
        $f$ & $-f$ & $0$ & $0$ & $0$ & $0$ & $0$ & $0$ & $0$ & $0$\\
        $g$ & $-g$& $-g$& $-f$ & $0$ & $0$ & $-b_{k'}$ & $0$ & $-d_{k'}$ & $0$ \\
        $a_k$ & $-a_{k}$ & $a_k$ & $0$ & $0$ & $b_{k}$ & $0$ & $0$ & $a_{k+k'}$ & $b_{k+k'}$\\
        $b_k$ & $-b_k$ & $0$ & $0$ & $0$ & $0$ & $0$ & $0$ & $0$ & $0$ \\
        $c_k$ & $0$ & $0$ & $-a_{k}$ & $0$ & $d_{k}$ & $-a_{k+k'}$ & $0$ & $0$ & $d_{k+k'}$ \\
        $d_k$ & $-d_k$ & $-d_k$ & $-b_{k}$ & $0$ & $0$ & $-b_{k+k'}$ & $0$ & $-d_{k+k'}$ & $0$ \\
       \end{tabular}}
       \end{center}

        \end{proposition}
        \begin{proof}
Proof of these bracket relations is exactly similar to the proof of the Proposition \ref{brackets_one}, so it is skipped.             
        \end{proof}

        \begin{theorem}\label{radical_1} With the notations as above, and for a fixed $k$, let $\mathfrak{b}_k$ be the span of $a_k,b_k,c_k,d_k$. Put $\mathfrak{a} = \langle e, f, g \rangle \oplus \mathfrak{b}$ where $\mathfrak{b}=\oplus_{k \geq 1} \mathfrak{b}_k$, then the following are true:
\begin{enumerate}
     \item For each $k$, $\mathfrak{b}_k$ is a Lie ideal in $sol_{A^2_n}$,  and  furthermore $\mathfrak{a}$ is also a Lie ideal.
    \item The Lie algebra $sol_{A^2_n}$ is solvable.
    
    \end{enumerate}
    
   \end{theorem}\label{nilpotent_1}
   \begin{proof}
       The first statement is obvious from the table. For the second statement, observe from the table of the Proposition \ref{brackets_three} that $[sol_{A^2_n}, sol_{A^2_n}] = \mathfrak{a}$. Now observe from the table that $ [ \mathfrak{a}, \mathfrak{a}] = \mathfrak{b} \oplus \langle f \rangle$ again from the table one has $[\mathfrak{b} \oplus \langle f \rangle , \mathfrak{b} \oplus \langle f \rangle ] \subseteq \oplus_{k \geq 2} \mathfrak{b}_k$. Continuing this way, we see that $sol_{A^2_n}$ is solvable.
       \end{proof}
       As a consequence, we have the following results:
\begin{corollary}  We have an exact sequence,
    \[ 0 \to  \mathfrak{a} \to sol_{A^2_n} \to \mathbb{C}^2 \to 0\]
    \end{corollary}
At the level of groups we have, 
\begin{corollary} 
Let be a subgroup $K$ corresponding to the Lie ideal $\mathfrak{a}$, then we have
\begin{itemize}
    \item $1 \to K \to (Sol_{A_n^2})^\circ \to \mathbb{C}^{*2} \to 1$
    \item The dimension of the group $Sol_{A^2_n}$ is $n$.
\end{itemize}
\end{corollary}

\subsubsection{Case 4: \: $\bf{n \equiv 3 (\text{mod}\ 4)}, \ \ n \geq 7$}

\begin{proposition}\label{generation_four}
The structure of the solution of the equation $X^tA^2_n + A^2_nX = 0$ in this case is given by,
 {\tiny $$\left(\begin{array}{ccccccccc}
        X_{11} &O_2&Y_{1,4}&O_2 &Y_{1,8}&O_2&\cdots & Y_{1,n-3}& O_2 \\
    X_{21}& X_{2,2}&Z_{2,4}& O_2&Y_{2,8}&O_2&\cdots&Y_{2,n-3}&O_2\\
    O_2 &O_2&-X_{2,2}^t&\ddots &\ddots&\ddots&\ddots&\ddots&O_2\\
    O_2&-Y_{1,4}^t&-Z_{2,4}^t&\ddots&\ddots&\ddots&\ddots&Y_{1,8}&O_2\\
    O_2&O_2&O_2 &\ddots &\ddots&\ddots&\ddots&Y_{2,8}&O_2\\
    \vdots &-Y_{1,6}^t&-Y_{2,8}^t&\ddots&\ddots&\ddots&\ddots&Y_{1,4}&O_2\\
    O_2&\vdots&\vdots &\ddots&\ddots&\ddots&X_{2,2}&Z_{2,4}&O_2\\
    O_2&O_2&O_2&O_2&O_2&O_2&O_2&-X_{2,2}^t&O_2\\
     O_2& -Y_{1,n-3}^t &-Y_{2,n-3}^t&\cdots& -Y_{1,8}^t&-Y_{2,8}^t&-Y_{1,4}^t&-Z_{2,4}^t&X_{2,2}
    \end{array}\right).$$}
Where, $X_{21} = {\footnotesize\left(\begin{array}{cc}
         0& x_{2,1} \\
         0&0 
    \end{array}\right)}; \  Z_{2,4} = {\footnotesize\left(\begin{array}{cc}
         x_{2,4}& x_{2,5} \\
         -x_{2,1}&0 
    \end{array}\right)};$ and $X_{11}, X_{2,2}, Y_{1,i}, Y_{2, i}$ as before for $i = 4+4k, \ k\leq \tfrac{n-7}{4}$ and $j = 8+4k, \ k\leq \tfrac{n-11}{4}$. In this case, add a zero row and zero column at the start.

\end{proposition}
\begin{proof}
    The proof of this case is exactly similar to the odd case $n \equiv 1 \mbox{ (mod 4)}$, so it is skipped.
\end{proof}

    \begin{proposition}\label{brackets_four}
        $sol_{A_n^2}$ is generated by the following elements:\\
        $$h_1 = \sum_{l=0}^{\tfrac{n-3}{4}}\left(E_{1+4l,1+4l} - E_{3+4l,3+4l}\right); \ \ \ h_2 =  \sum_{l=0}^{\tfrac{n-3}{4}}\left(E_{2+4l,2+4l} - E_{4+4l,4+4l}\right);$$
        $$e = \sum_{l=0}^{\tfrac{n-3}{4}}\left(E_{2+4l,1+4l} - E_{3+4l,4+4l}\right); \ \ \ f = \sum_{l=0}^{\tfrac{n-3}{4}}\left(E_{2+4l,3+4l} - E_{5+4l,4+4l}\right);$$
        
        For $k = 1, 2, \ldots, \frac{(n-3)}{4}$,
        $$a_k =\sum_{l=0}^{\tfrac{n-3-4k}{4}}\left(E_{1+4l,4l+4k} - E_{2+4l+4k,3+4l}\right); \ \ \ b_k = \sum_{l=0}^{\tfrac{n-3-4k}{4}}\left(E_{1+4l,1+4l+4k} - E_{3+4l+4k,3+4l}\right);$$
        $$c_k = \sum_{l=0}^{\tfrac{n-3-4k}{4}}\left(E_{2+4l,4l+4k} - E_{2+4l+4k,4+4l}\right); \ \ \ d_k = \sum_{l=0}^{\tfrac{n-3-4k}{4}}\left(E_{2+4l,1+4l+4k} - E_{3+4l+4k,4+4l}\right),$$
         with the following Lie brackets:
         
        {\tiny \begin{center}
        \begin{tabular}{c|cccccccc}
        & $h_1$ & $h_2$ & $e$ & $g$ & $a_{k'}$ & $b_{k'}$ & $c_{k'}$ & $d_{k'}$ \\
        \hline
        $h_1$ & $0$ & $0$ & $-e$ & $g$ & $a_{k'}$ & $0$ & $0$ & $d_{k'}$\\
        $h_2$ & $0$ & $0$& $e$& $g$ & $0$ & $0$ & $c_{k'}$ & $d_{k'}$\\
        $e$ & $e$ & $-e$ & $0$ & $0$ & $c_{k'}$ & $d_{k'}$ & $0$ & $0$ \\
        $g$ & $-g$ & $-g$ & $0$ & $0$ & $0$ & $a_{k'}$ & $0$ & $c_{k'}$ \\
        $a_k$ & $-a_{k}$ & $0$ & $-c_{k}$ & $0$ & $0$ & $-a_{k+k'}$ & $0$ & $-c_{k+k'}$\\
        $b_k$ & $0$ & $0$ & $-d_{k}$ & $-a_{k}$ & $a_{k+k'}$ & $0$ & $0$ & $-d_{k+k'}$ \\
        $c_k$ & $0$ & $-c_{k}$ & $0$ & $0$ & $0$ & $0$ & $0$ & $0$ \\
        $d_k$ & $d_k$ & $-d_k$ & $0$ & $-c_{k}$ & $c_{k+k'}$ & $d_{k+k'}$ & $0$ & $0$ \\
       \end{tabular}
       \end{center}}
    \end{proposition}
\begin{proof}
Proof of these bracket relations is same as the proof of the Proposition \ref{brackets_one}, so it is skipped.             
        \end{proof}
    \begin{theorem}\label{radical_3}
         With the notation as above and for a fixed $k$, let $\mathfrak{b}_k$ be the span of $a_k,b_k,c_k,d_k$, and let $\mathfrak{a}$ be $\langle e, g \rangle \oplus \mathfrak{b}$ where $\mathfrak{b}=\oplus_{k} \mathfrak{b}_k$, then the following are true:
\begin{enumerate}
     \item For each $k$, $\mathfrak{b}_k$ is a Lie ideal in $sol_{A_n^2}$,  and  furthermore $\mathfrak{a}$ is also a Lie ideal.
    \item The Lie algebra $sol_{A_n^2}$ is solvable.
    
    \end{enumerate}
    \end{theorem}
\begin{proof}
 Observe that $[sol_{A_n^2}, sol_{A_n^2}] \subset \mathfrak{a}$ and also $[\mathfrak{a}, \mathfrak{a}] \subset \mathfrak{b}$, and $[\mathfrak{b}, \mathfrak{b}]\subset \mathfrak{b}_1$, and $[\mathfrak{b}_i, \mathfrak{b}_i] \subset \mathfrak{b}_{i+1}$ which implies $sol_{A_n^2}$ is solvable. This completes the proof.
    
\end{proof}
As a result, we have the following corollaries:
\begin{corollary}
  We have an exact sequence,
    \[ 0 \to  \mathfrak{a} \to sol_{A^2_n} \to \mathbb{C}^2 \to 0\]
    \end{corollary}
    \begin{corollary} At the level of groups, we have
    \begin{itemize}
    \item
    $1 \to K \to (Sol_{A_n}^2)^\circ \to \mathbb{C}^{*2} \to 1$
    \item The dimension of the group $Sol_{A^2_n}$ is $n+1$.
\end{itemize}
   
   \end{corollary}
   \begin{remark}\label{Invariant_odd_remark}
        One can observe that for each $m \geq 1$, there exist a Lie algebra morphism $\phi : \text{sol}_{A^2_{5+4m}} \longrightarrow \text{sol}_{A^2_{3+4m}}$ defined as $\phi(h_1) = h_2, \phi(h_2) = h_1, \phi(e) = e, \phi(g) = -f, \phi(f) = 0, \phi(a_k) = d_k, \phi(b_k) = c_k, \phi(c_k) = b_k,$ and $\phi(d_k) = a_k.$ It implies that the group $\text{Sol}_{A^2_{5+4m}}$ is a semidirect product of $\mathbb{C}^{*}$ and $ \text{Sol}_{A^2_{3+4m}}$. 
   \end{remark}

    \subsection{Structure of the group}
    This section thoroughly describes the group structure of $Sol_{A^2_n}$ for all $n\geq 4$. We will consider different cases of $n$ modulo 4 in the following subsections and describe the group in each case. We will show that in all these cases $Sol_{A^2_n}$ splits into a semi-direct product of two subgroups $D$ and $N$, coming naturally from the Lie algebra of $Sol_{A^2_n}$. As a result of the treatment, one can write down the parametric equations of the group $Sol_{A^2_n}$ as an algebraic group. 
\begin{lemma}
Define $$D = \begin{cases} \text{diag}(g, g^{-t}, \cdots, g, g^{-t}) & \text{ when } n = 0 (\text{mod } 4),\\
\text{diag}(g, g^{-t}, \cdots, g, g^{-t}, g) & \text{ when } n = 2(\text{mod } 4),\\
\text{diag}(g_0, g_1, g_1^{-t}, \cdots, g_1, g_1^{-t}) & \text{ when } n = 1(\text{mod } 4),\\
\text{diag}(g_0, g_1, g_1^{-t}, \cdots, g_1, g_1^{-t}, g_1) & \text{ when } n = 3(\text{mod } 4),
\end{cases}$$ where $g \in GL_2(\mathbb{C}), g_0 = {\footnotesize\begin{pmatrix} 0&0\\
0&\beta
\end{pmatrix}},$and $ g_1= {\footnotesize\begin{pmatrix}\alpha &\ast \\
0 & \beta^{-1} \end{pmatrix}}$, for $\alpha, \beta \neq 0$. Then $D$ is a subgroup of $\text{Sol}_{A^2_n}.$
\end{lemma}
\begin{proof} For any matrix $d \in D$, we need to show that $d^tA_n^2d = A_n^2$. After decomposing each matrix into $2 \times 2$ blocks (by amending with zeros if required), we can visualize each one of them as an matrix in $M_{\tfrac{n}{2}}(M_2(\mathbb{C}))$ (or in $M_{\tfrac{n+1}{2}}(M_2(\mathbb{C}))$. It boils down to show $d^tA_{n/2}d = A_{n/2}$ (or $\tilde{d^t}\tilde{A_{n/2}}\tilde{d} = \tilde{A_{n/2}}$). For $n \equiv 0(\text{mod }n),$ we have
\begin{equation*}
    \begin{split}
        d^tA_n^2d &= {\tiny \begin{pmatrix}
        g^t&  & && \\
    &g^{-1}& && \\
    &&\ddots && \\
    &&&g^t&\\
    &&& &g^{-1}
    \end{pmatrix}\begin{pmatrix}
        O_2& g^{-t} & && \\
    &O_2& g&& \\
    &&\ddots &\ddots& \\
    &&&O_2&g^{-t}\\
    &&& &O_2
    \end{pmatrix}} = A_n^2.
    \end{split}
\end{equation*}
The case $n \equiv 2(\text{mod} 4)$ is similarly followed. Now for $n \equiv 1(\text{mod }4)$, we have
\begin{equation*}
    \begin{split}
        \tilde{d^t}\tilde{A_n^2}\tilde{d} &= {\tiny \begin{pmatrix}
        g_0^t &  & & & \\
    &g_1^t& && \\
    &&g_1^{-1}&&\\
    &&&\ddots &\\
    &&&& g_1^{-1}
    \end{pmatrix}\begin{pmatrix}
        O_2& a& && \\
    &O_2& I_2&& \\
    &&\ddots &\ddots& \\
    &&&\ddots&I_2\\
    &&& &O_2
    \end{pmatrix} \begin{pmatrix}
        g_0&  & && \\
    &g_1& &&\\
    &&g_1^{-t}&&\\
    &&&\ddots & \\
    &&& &g_1^{-t}
    \end{pmatrix}}\\
    &= {\tiny \begin{pmatrix}
        g_0^t&  & && \\
    &g_1^t& && \\
    &&g_1^{-1}&&\\
    &&&\ddots &\\
    &&& &g_1^{-1}
    \end{pmatrix}\begin{pmatrix}
        O_2& ag_1 & && \\
    &O_2& g^{-t}&& \\
    &&\ddots &\ddots& \\
    &&&O_2&g^{-t}\\
    &&& &O_2
    \end{pmatrix}} = {\tiny \begin{pmatrix}
        O_2& g_0^tag_1& && \\
    &O_2& I_2&& \\
    &&\ddots &\ddots& \\
    &&&\ddots&I_2\\
    &&& &O_2
    \end{pmatrix} }\\
    &= \tilde{A_n^2} \ \ (\text{because }g_0^tag_1 = a).
    \end{split}
\end{equation*}
Similarly one can verify the case $n \equiv 3(\text{mod} 4)$. 
 \end{proof}

Now we will construct the subgroup $N$ in the following subsections regarding its local matrix coordinates. Note that a non-zero, non-constant coordinate function near identity will contribute to the Lie algebra. Thus, we can expect the equations of the group to involve only the non-zero Lie algebra coordinates. As in the Lie algebra case, we must deal with four cases depending on $n$ modulo four. 

   \subsubsection{The case of $n$ congruent to zero mod 4}\
   
   \textbf{Construction of $N$.} For $1 \leq i \leq \frac{n}{2}$, $C_i$ denotes the column vector with entries from $M_2(\mathbb{C})$, in which 1 is at the $i^\text{th}$-place and is defined as follows
\begin{equation}
    \text{For $i$ even}, \ \ C_i = {\footnotesize \begin{bmatrix}
     0  \\
     \vdots\\
     0\\
     1\\
     0\\
     -x_1^t + x_1^\ast\\
     \vdots\\
     0\\
     -x_{\tfrac{n-i}{2}}^t + x_{\tfrac{n-i}{2}}^\ast
\end{bmatrix}}; \ \ \text{For $i$ odd}, \ \ C_i = {\footnotesize\begin{bmatrix}
     x_{\tfrac{i-1}{2}}  \\
     0\\
     \vdots\\
     x_1\\
     0\\
     1\\
     0\\
     \vdots\\
     0
\end{bmatrix}},
\end{equation}
where $x_l^\ast = \sum_{r+s = l} x^t_s (x_r^t - x_r^\ast).$
Thus the transpose of columns will look like
$$\text{For $i$ even}, \ \ C_i^t = {\footnotesize\begin{bmatrix}
     0  & \cdots& 0&1& 0& -x_1 + x_1^{\ast t}& \cdots&0&-x_{\tfrac{n-i}{2}} + x_{\tfrac{n-i}{2}}^{\ast t}
\end{bmatrix}},$$
and $$\text{for $i$ odd}, \ \ C_i^t = {\footnotesize\begin{bmatrix}
     x_{\tfrac{i-1}{2}}^t & 0& \cdots& x_1^t& 0& 1& 0& \cdots& 0
\end{bmatrix}}.$$
Call $A_n^2C_i = \Tilde{C_i},$ then $\Tilde{C_i}$ will be the shifted column of the matrix above by one index and append two by two zero matrix in the last entry, such that $1$ appears at $(i-1)^\text{th}$-place, and it looks as follows:
$$\text{For $i$ even}, \ \ \widetilde{C_i} = {\tiny \begin{bmatrix}
     0  \\
     \vdots\\
     0\\
     1\\
     0\\
     -x_1^t + x_1^\ast\\
     \vdots\\
     0\\
     -x_{\tfrac{n-i}{2}}^t + x_{\tfrac{n-i}{2}}^\ast\\
     0
\end{bmatrix}}; \ \ \text{For $i$ odd}, \ \ \widetilde{C_i} = {\tiny\begin{bmatrix}
     0\\
     x_{\tfrac{i-3}{2}}  \\
     0\\
     \vdots\\
     x_1\\
     0\\
     1\\
     0\\
     \vdots\\
     0
\end{bmatrix}}.$$
Define $N = \left\{\begin{bmatrix} C_1& C_2& \ldots& C_{\tfrac{n}{2}}\end{bmatrix}\right\}$. It is a collection of $\tfrac{n}{2} \times \tfrac{n}{2}$-matrices whose columns are $C_i$ and $x_i$s belong to $M_2(\mathbb{C})$.
\begin{lemma}\label{divisible_by_four} $N$ defined as above is a subset of the solution set $\{X \ | \ X^t A^2_n X = A_n^2\}$ such that $Lie(N)\simeq \mathfrak{a}$.
\end{lemma}
\begin{proof}
Clearly, the $ij$-th entry of the matrix $\begin{bmatrix} C_1& C_2& \ldots& C_{\tfrac{n}{2}}\end{bmatrix}^t A_n^2 \begin{bmatrix} C_1& C_2& \ldots& C_{\tfrac{n}{2}}\end{bmatrix}$ is $C_i^t \widetilde{C}_j.$
\begin{itemize}
    \item[] \textbf{Claim 1.} $C_i^t \widetilde{C}_{i+1} = 1$ for all $i$.\\   
    When $i$ is even, $i+1$ is odd. Then $\widetilde{C}_{i+1}$ is the column in which 1 appears at $i$-th place (even place).
Hence \[ C_i^t \widetilde{C}_{i+1} = {\tiny \begin{matrix} &&i^{\text{th}}& \text{place}&&&&&\\
&&&\big\downarrow &&&&&\\
 \Big[ 0 &\cdots &0 & 1 & 0 &-x_1 + x_1^{\ast t} & \cdots &0&-x_{\tfrac{n-i}{2}} + x_{\tfrac{n-i}{2}}^{\ast t} \Big] \\
 &&&&&&&&\end{matrix}
\begin{bmatrix}
     0\\
     x_{\tfrac{i-2}{2}}  \\
     0\\
     \vdots\\
     x_1\\
     0\\
     1\\
     0\\
     \vdots\\
     0
\end{bmatrix} \begin{matrix}
\\ \\ \\ \\ \\ \\ \longleftarrow \\ \\ \\ \\
\end{matrix}\begin{matrix}
\\ \\ \\ \\ \\ i^{\text{th}} \text{ place} \\ \\ \\ \\ 
\end{matrix}
= 1} \]

When $i$ is odd, $i+1$ is even. Then
\[C_i^t \widetilde{C}_{i+1} = {\tiny \begin{matrix} &&&&i^{\text{th}}& \text{place}&&&\\
&&&&&\big\downarrow &&&\\
     \Big[x_{\tfrac{i-1}{2}}^t  &  0& \cdots& x_1^t&0& 1 & 0 & \cdots& 0\Big]\\
     &&&&&&&&
\end{matrix}\begin{bmatrix} 0 \\ \vdots \\ 0 \\ 1 \\ 0 \\-x_1^t + x_1^{\ast} \\ \vdots \\0\\ -x_{\tfrac{n-i-1}{2}}^t + x_{\tfrac{n-i-1}{2}}^{\ast} \\0 \end{bmatrix}\begin{matrix}
\\ \\  \\ \\ \longleftarrow \\ \\ \\ \\ \\ \\ \\ \\
\end{matrix}\begin{matrix}
\\ \\ i^{\text{th}} \text{ place} \\  \\ \\ \\ \\ \\  \\
\end{matrix}  = 1.}\]

\item[] \textbf{Claim 2.} $C_i^t \widetilde{C}_j = 0$ whenever $j-i \neq 1.$\\

\textbf{Subcase a.} $i \equiv j \, (\text{ mod } 2)$. If both $i$ and $j$ are even, then 
\[C_i^i \widetilde{C}_j = {\tiny \begin{matrix} &&i^{\text{th}}& \text{place}&&&&&\\
&&&\big\downarrow &&&&&\\
     \Big[0 &\cdots &0 & 1 & 0 &-x_1 + x_1^{\ast t} & \cdots &0&-x_{\tfrac{n-i}{2}} + x_{\tfrac{n-i}{2}}^{\ast t} \Big]\\
     &&&&&&&&
\end{matrix} \begin{bmatrix} 0 \\ \vdots \\ 0 \\ 1 \\ 0 \\-x_1^t + x_1^{\ast} \\ \vdots \\0\\ -x_{\tfrac{n-j-1}{2}}^t + x_{\tfrac{n-j-1}{2}}^{\ast}\\0 \end{bmatrix}\begin{matrix}
\\ \\  \\  \\ \longleftarrow \\ \\ \\ \\ \\ \\ \\ \\
\end{matrix}\begin{matrix}
\\ \\ (j-1)^{\text{th}} \text{ place} \\  \\ \\ \\ \\ \\  \\
\end{matrix}} \]
which is clearly zero because the last entry of $C_i^t$ is non-zero and in $\widetilde{C}_j$ it is zero. Also, the zero and non-zero positions are alternate, thus they will never match as $i$ is even and $j-1$ is odd.\\

If $i$ and $j$ both are odd, then  \[C_i^i \widetilde{C}_j = {\tiny \begin{matrix}&&&&i^{\text{th}}& \text{place}&&&\\
&&&&&\big\downarrow &&&\\
     \Big[x_{\tfrac{i-1}{2}}^t  &  0& \cdots& x_1^t&0& 1&0& \cdots& 0\Big]\\
      &&&&&&&&
\end{matrix}\begin{bmatrix}
     0\\
     x_{\tfrac{j-2}{2}}  \\
      0\\
     \vdots\\
     x_1\\
     0\\
     1\\
     0\\
     \vdots\\
     0
\end{bmatrix} \begin{matrix}
\\ \\  \\ \\ \\ \\ \longleftarrow \\ \\ \\ \\
\end{matrix}\begin{matrix}
\\ \\ \\  \\ \\ \\ (j-1)^{\text{th}} \text{ place}  \\  \\ \\
\end{matrix}}\]
which is also clearly zero as first entry of $C_i^t$ is non-zero and it is zero in $\widetilde{C}_j$. Also, the zero and non-zero entries are alternate.\\

\textbf{Subcase b.} $i \not\equiv j (\text{ mod } 2)$ such that $j -i \neq 2$. Without loss of generality, one can assume, $i$ is odd and $j$ is even, and $i>j$. Then $C_i^t \widetilde{C}_j = $ 
{\tiny \[ \begin{matrix}
\\ \\  (j-1)^{\text{th}} \text{ place} \\ \\ \\ \\ \\ i^{\text{th}} \text{ place}  \\ \\  \\ 
\end{matrix}\begin{matrix}
\\ \\ \\ \longrightarrow  \\ \\ \\ \\ \\ \longrightarrow \\ \\ \\ 
\end{matrix}\begin{bmatrix}
     x_{\tfrac{i-1}{2}}  \\ 0\\ \vdots\\ x_{\tfrac{i-j-1}{2}} \\ 0\\\vdots \\ x_1 \\ 0\\ 1 \\  0\\ \vdots\\ 0
\end{bmatrix}^t\begin{bmatrix} 0 \\ \vdots \\ 0 \\ 1 \\ 0 \\-x_1^t + x_1^{\ast} \\ \vdots \\ 0\\ -x_{\tfrac{i-j-1}{2}}^t + x_{\tfrac{i-j-1}{2}}^\ast \\ \vdots \\-x_{\tfrac{n-j-1}{2}}^t + x_{\tfrac{n-j-1}{2}}^{\ast} \\0 \end{bmatrix} \begin{matrix}
\\ \\ \longleftarrow \\ \\ \\  \\ \\ \\ \longleftarrow \\ \\ \\ 
\end{matrix}\begin{matrix}
\\ \\ (j-1)^{\text{th}} \text{ place}  \\ \\ \\ \\ \\  i^{\text{th}} \text{ place}  \\  \\ \\
\end{matrix}\]}
After taking the dot product of vectors, we get
$$C_i^t \widetilde{C}_j = x_{\tfrac{j-i-1}{2}}^t + \sum_{r+s = \tfrac{i-j-1}{2}} x^t_s (-x_r^t + x_r^\ast) + -x_{\tfrac{j-i-1}{2}}^t + x_{\tfrac{j-i-1}{2}}^\ast$$
$$=0 \ \ (\text{using the definition of } x_l^\ast)$$
\end{itemize}

This completes the proof that $N \subseteq Sol_{A^2_n}$.
\end{proof}

\begin{proposition}
    $N$ is a normal subgroup of $Sol_{A^2_n}.$
\end{proposition}
\begin{proof} 
This is done in the following steps:
\begin{itemize}
    \item[]\textbf{Step 1.} For $X, Y \in N$, we need to show $XY \in N.$ For that define the rows $R_i$ (with 1 at $i^{\text{th}}$ place) of any matrix $X$ of $N$ as follows:
$$ \text{For } i \text{ odd }, \ \ R_i = {\footnotesize\begin{bmatrix} 0 &\cdots & 0 & 1 & 0 & x_1 & 0 & \cdots & x_{\tfrac{n-i-1}{2}} & 0 \end{bmatrix}}$$
and $$\text{For } i \text{ even }, \ \ R_i = {\footnotesize\left[\begin{array}{cccccccccccc} 0 & -x_{\tfrac{i-2}{2}}^t + x_{\tfrac{i-2}{2}}^\ast & \cdots & 0 & -x_2^t + x_2^\ast & 0 & -x_1^t  & 0 & 1 & 0 & \cdots&0
\end{array}\right]}.$$
The columns are the same as given in the equation (3). Then $(ij)$th entry of matrix $XY$ is $R_i C_j.$ Clearly for all $i$, \ $R_i C_i = 1.$ Also for $i$ even (or odd) and $j$ odd (or even), $R_iC_j = 0.$\\
If $i, j$ both are odd and $i < j$, then $$R_iC_j =  x_{\tfrac{j-i}{2}}+y_{\tfrac{j-i}{2}} + \sum_{l+k = \tfrac{j-i}{2}} x_l y_k$$

If $i, j$ both are even and $i>j$, then $$R_iC_j = - \left(x_{\tfrac{i-j}{2}}^t +  y_{\tfrac{i-j}{2}}^t\right) +  \left(x_{\tfrac{i-j}{2}}^\ast + y_{\tfrac{i-j}{2}}^\ast\right)+ \sum_{l+k = \tfrac{i-j}{2}} (-x_l^t+ x_l^\ast)(-y_k^t+ y_k^\ast).$$
From the structure of any matrix in $N$, clearly, the matrix $XY$ is determined by the first row, which is $R_1C_j$. Put $z_{\tfrac{j-1}{2}} = R_1C_j$ for $j = 3, 5, \ldots, n-1$.\\
\textbf{Claim.} The last row $R_nC_j$ of $XY$ should be related to first row as follows: $R_nC_{n-j+1} = -z_{\tfrac{j-1}{2}}^t + z_{\tfrac{j-1}{2}}^\ast$ for $j = 3, 5, \ldots, n-1$. Here $z_{-}^\ast$ is given by equation \ref{def_of_star}. That is,

\begin{equation*}
\begin{split}
   - \left(x_{\tfrac{j-1}{2}}^t +  y_{\tfrac{j-1}{2}}^t\right) +  \left(x_{\tfrac{j-1}{2}}^\ast + y_{\tfrac{j-1}{2}}^\ast\right)+ \sum_{l+k = \tfrac{j-1}{2}} (-x_l^t+ x_l^\ast)(-y_k^t+ y_k^\ast) &= -z_{\tfrac{j-1}{2}}^t + z_{\tfrac{j-1}{2}}^\ast\\
 = - x_{\tfrac{j-1}{2}}^t -  y_{\tfrac{j-1}{2}}^t - \sum_{l+k = \tfrac{j-1}{2}} y_k^t x_l^t + z_{\tfrac{j-1}{2}}^\ast.
\end{split}
\end{equation*}
That is,
\begin{equation}
    \label{Prod_multiplicative}
z_n^\ast = x_n^\ast + y_n^\ast + \sum_{l+k = n} y_k^t x_l^t + \sum_{l+k = n} (-x_l^t+ x_l^\ast)(-y_k^t+ y_k^\ast).
\end{equation}
\textbf{Proof of claim.} We prove it by induction on $n$. Clearly, $z_1^\ast = 0$ as $x_1^\ast = y_1^\ast=0$, and $$z_2^\ast = z_1^{2t} = (x_1^t + y_1^t)^2 = x_2^\ast + y_2^\ast + y_1^1x_1^t + x_1^ty_1^t.$$ 
For show general $n$ case as follows: For a matrix variable $x$, define $L(x)= x^t -x^{*}$, where $*$ is defined as above. Note that we would like to prove  \[ z_n^{*} = x_n^* + y_n^* + \sum_{a+b=n} y_a^t x_b ^t + \sum_{a+b=n} (-x_a^t +x_a^*)(-y_b^t + y_b^*).\] Equivalently we would like to prove,
\[ L(z_n)=L(x_n)+L(y_n) - \sum_{a+b=n} L(x_a)L(y_b).\] Note that using the definition of $z_n^{*}= \displaystyle\sum_{a+b=n} z_a^t L(z_b)$, one can rewrite the previous equation, as 
\[\sum_{a+b=n} z_a^t L(z_b) = x_n^{*} + y_n^{*} + \sum_{a+b=n} y_a^t x_b^t + \sum_{a+b=n} L(x_a)L(x_b),\]
which follows from the Theorem \ref{star_one}. \\

\item[] \textbf{Step 2.} For any $X \in N$, we show the existence of its inverse $X^{-1}$ such that $X X^{-1} = I.$ Define $X^{-1}$ whose columns are as follows
$$\text{For $i$ even}, \ \ \widetilde{C_i} = {\tiny \begin{bmatrix}
     0  \\
     \vdots\\
     0\\
     1\\
     0\\
     x_1^t\\
     \vdots\\
     0\\
     x_{\tfrac{n-i}{2}}^t\\
     0
\end{bmatrix}}; \ \ \text{For $i$ odd}, \ \ \widetilde{C_i} = {\tiny\begin{bmatrix}
     0\\
     -x_{\tfrac{i-3}{2}} + x_{\tfrac{i-3}{2}}^{\ast t}\\
     0\\
     \vdots\\
     -x_1+x_1^{*t}\\
     0\\
     1\\
     0\\
     \vdots\\
     0
\end{bmatrix}}.$$
One can easily check for $i \neq j$, $(XX^{-1})_{ij} = 0$\\

\item[] \textbf{Step 3.} Once we ignore the non-linear terms in the coordinates of $N$, we get that $Lie(N)= \mathfrak{a}$. the normality of the subgroup $N$ follows from this fact, and $\mathfrak{a}$ is a Lie ideal in $sol_{A_n^2}$.
\end{itemize}

\end{proof}

We have a closed-form formula for the non-linear part of the coordinate functions of the group $N$ given by the following remark:
   \begin{remark}
        $$ x_l^{\ast t} = (-1)^l \text{det}{\tiny\left(\begin{array}{ccccc}
        x_1 & x_2 & \cdots &x_{l-1} &0\\
        1 & x_1 &x_2 &\cdots & x_{l-1}\\
        0 & \ddots &\ddots &\ddots & \vdots \\
        \vdots &\ddots&\ddots&\ddots&x_2\\
        0&\cdots&0&1&x_1
    \end{array}\right)}_{l \times l}.$$
   \end{remark}
   
Clearly, one observes that the determinant of any element in $X \in Sol_{A^2_n}$ is one. Thus we have
\begin{corollary}
    The group $Sol_{A^2_n}$ is a subgroup of $SL_n$.
\end{corollary}

\subsubsection{The case of $n$ congruent to 2 modulo four}\

\textbf{Construction of $N$.} Put $m = \frac{n-2}{4}$. Define $N = \{I+X \ | \ X \in \mathfrak{a}\}$ where $\mathfrak{a}$ is defined in Theorem \ref{ideal-two}.
\begin{lemma}
    For all $X \in \mathfrak{a},$ we have $X^t A X = 0$.
\end{lemma}
\begin{proof}
$X$ is a matrix in $M_{n/2}(M_2(\mathbb{C}))$ whose columns are given as follows:
$$\text{For $i$ even}, \ \ C_i = {\footnotesize \begin{bmatrix}
     x_{\tfrac{i}{2}}  \\
     0\\
     \vdots\\
     x_1\\
     0\\
-x_1^t\\
\vdots\\
     0\\
     -x^t_{m+1-\tfrac{i}{2}}
\end{bmatrix}}; \ \ \text{For $i$ odd}, \ \ C_i = {\footnotesize\begin{bmatrix}
     0\\
     \vdots\\
     0\\
     0\\
     \vdots\\
     0
\end{bmatrix}}.$$
Thus the truncated column $\widetilde{C}_i$ would be
$$\text{For $i$ even}, \ \ \widetilde{C}_i ={\footnotesize \begin{bmatrix}
     0\\
     x_{\tfrac{i-2}{2}}  \\
     \vdots\\
     x_1\\
     0\\
-x_1^t\\
0\\
\vdots\\
     -x^t_{m+1-\tfrac{i}{2}}\\
     0
\end{bmatrix}}; \ \ \text{For $i$ odd}, \ \ C_i = {\footnotesize\begin{bmatrix}
     0\\
     \vdots\\
     0\\
     0\\
     \vdots\\
     0
\end{bmatrix}}.$$
From the shape of the columns, it can be obtained that $C_i^t \widetilde{C}_j = 0$ for all $i,j$.
\end{proof}

\begin{proposition}
    The set $N= I + \mathfrak{a}$ is a normal subgroup of $Sol_{A^2_n}$ with $Lie(N) \simeq \mathfrak{a}$, and where $I$ is the diagonal matrix with 1 and -1 alternatively on the diagonal.
\end{proposition}
\begin{proof}
    Let $I+X$ and $I+Y$ be elements of $N$. Then the product $(I+X)(I+Y)= I + X+Y + XY$. Observe that, through a direct computation, $(XY)^tA^2_n(XY)= Y^tX^tA_n^2XY=0$, as a result $((1+X)(1+Y))^tA_n^2(1+X)(1+Y)=A_n^2$. This shows that $N$ is closed under multiplication. Clearly the elements $(I+X)$ are invertible with inverse $I-X$ in $N$. Since the Lie algebra of $N$ coincides with $\mathfrak{a}$, which is a Lie ideal of the Lie algebra of $Sol_{A^2_n}$, which implies $N$ is normal.
\end{proof}

\subsubsection{The case of $n$ congruent to 1 modulo 4}\

\textbf{Construction of $N$.} Define $a = {\footnotesize\begin{pmatrix}
0 &0\\
0&1
\end{pmatrix}},$ $x_i =  {\footnotesize\begin{pmatrix}
0 &0\\
\ast&\ast
\end{pmatrix}};  \text{ if $i$ is odd}$, and
$ x_i = {\footnotesize\begin{pmatrix}
\ast &\ast\\
0&0
\end{pmatrix}};  \text{ if $i$ is even}$. Then for $2 \leq i \leq \frac{n+1}{2}$, $C_i$ denotes the $i^\text{th}$ column vector, in which 1 is at the $i^\text{th}$-place and is defined as follows
\begin{equation} \label{columns_onemodfour} C_1 = {\footnotesize\begin{bmatrix}
     a\\
     0\\
     \vdots\\
     0
\end{bmatrix}}; \ \text{For $i$ even}, \ \ C_i = {\footnotesize \begin{bmatrix}
     x_{i-1}  \\
     \vdots\\
     x_1\\
     1\\
     -x_1^t + x_1^\ast\\
     0\\
     -x_3^t + x_3^\ast\\
     \vdots\\
     0\\
     -x_{\tfrac{n+1}{2} -i}^t + x_{\tfrac{n+1}{2}-i}^\ast
\end{bmatrix}}; \ \ \text{For $i$ odd, $i \geq 3$}, \ \ C_i = {\footnotesize\begin{bmatrix}
     0\\
     \vdots\\
     0\\
     1\\
     0\\
     -x_2^t + x_2^\ast\\
     \vdots\\
     0\\
     -x_{\tfrac{n+1}{2} -i}^t + x_{\tfrac{n+1}{2}-i}^\ast
\end{bmatrix}}.\end{equation}
Define $$ x_{2n}^\ast = -\sum_{\substack{r+s = 2n \\
r,s \text{ even}}} x_r^t x_s^\ast + \sum_{\substack{r+s = 2n \\
r,s \text{ even}}}x_r^t x_s^t,$$
where $x^\ast_2 = 0;$ and
$$ x_{2n+1}^\ast = -\sum_{l-k = 2n+1} x_l^t x_k + \sum_{\substack{l+k = 2n+1 \\
k \text{ odd}}}x_l^t(x_k^t - x_k^\ast).$$
Observe the following lemma:
\begin{lemma}\label{star_new}
   For $a,b$, not congruent to each other modulo 2, then $x_a^t x_b$, $x_a^ty_b$, $x_a^*y_b$, and $x_ay_b^* $ are zero. 
\end{lemma}
\begin{proof}
    We can prove this using an easy induction argument on the indices $a,b$, and observe the fact that $x_a^ty_b$ for such pairs $a,b$ is zero.
\end{proof}
Now by the Lemma \ref{star_new}, we reduce the expression further to, \[ x_{2n+1}^\ast =\sum_{\substack{l+k = 2n+1 \\
k \text{ odd}}}x_l^t(x_k^t - x_k^\ast).\]
As a result, both the odd and even cases get unified into a simple, $$ x_{n}^\ast = -\sum_{\substack{r+s = n 
}} x_r^t x_s^\ast + \sum_{\substack{r+s = n }}x_r^t x_s^t,$$

which gives $x_1^* = 0$. Let $N$ be the subset of matrices given by the columns as described above, namely the matrices $\big[ C_1, C_{2}, C_{3}, \ldots C_{(n+1)/2}\big]$, where $x_i \in M_2(\mathbb{C})$.

\begin{lemma} $N$ is a subset of $\text{Sol}_{A_n^2}.$ \end{lemma}
\begin{proof}
Recall from the previous subsection that $\tilde{C_i}$ stands for the vector $A_n^2C_i$.
\begin{itemize}
\item[Case 1.] $C_i^t \widetilde{C}_{i+1} = 1,$ for all $i$. For $i$ even, we have $C_i^i \widetilde{C}_{i+1} = $
{\tiny \[\begin{matrix} &&i^{\text{th}}& \text{place}&&&&&\\
&&&\big\downarrow &&&&&\\
     \Big[x_{i-1}^t &\cdots &x^t_1 & 1 &-x_1 + x_1^{\ast t} &0& \cdots &0&-x_{\tfrac{n+1}{2}-i} + x_{\tfrac{n+1}{2}-i}^{\ast t}\Big]\\
     &&&&&&&&
\end{matrix}\begin{bmatrix} 0 \\ \vdots \\ 0 \\ 1 \\ 0 \\-x_2^t + x_2^{\ast} \\ \vdots \\0\\ -x_{\tfrac{n+1}{2}-j}^t + x_{\tfrac{n+1}{2}-j}^{\ast}\\0 \end{bmatrix}\begin{matrix}
\\ \\  \\  \\ \longleftarrow \\ \\ \\ \\ \\ \\ \\ \\
\end{matrix}\begin{matrix}
\\ \\ (j-1)^{\text{th}} \text{ place} \\  \\ \\ \\ \\ \\  \\
\end{matrix}\]}

and for $i$ odd, $C_i^t \widetilde{C}_{i+1} = $
{\tiny \[\begin{matrix}
&&i^{\text{th}}& \text{place}&&&&&\\
&&&\big\downarrow &&&&&\\
     \Big[0 &\cdots &0 &1 &0&-x_2 + x_2^{\ast t} & \cdots &0&-x_{\tfrac{n+1}{2}-i} + x_{\tfrac{n+1}{2}-i}^{\ast t}\Big] \\
     &&&&&&&&
\end{matrix}\begin{bmatrix} x_{i-2} \\ \vdots  \\ x_1\\ 1  \\-x_1^t + x_1^{\ast} \\ \vdots \\0\\ -x_{\tfrac{n+1}{2}-j}^t + x_{\tfrac{n+1}{2}-j}^{\ast}\\0 \end{bmatrix}
\begin{matrix}
\\ \\  \\  \\ \\ \longleftarrow \\ \\ \\ \\ \\ \\ \\ \\
\end{matrix}\begin{matrix}
\\ \\ \\ (j-1)^{\text{th}} \text{ place} \\ \\  \\ \\ \\  \\
\end{matrix}\]}

In both cases, the value is 1.\\

    \item[Case 2.] Consider $i$ even and $j$ odd such that $i > j$. Then $C_i^t \widetilde{C}_j =$
   {\tiny  \[  \begin{matrix}
(j-1)^{\text{th}} \text{ place} \\ \\ \\ \\ i^{\text{th}} \text{ place}  \\ \\  \\ \\ \\ \\ \\ 
\end{matrix}\begin{matrix}
 \longrightarrow \\ \\ \\ \longrightarrow \\ \\ \\ \\ \\ \\ \\
\end{matrix}\begin{bmatrix}
     x_{i-1}  \\ \vdots\\ x_{i-j+1} \\ \vdots\\x_1 \\ 1 \\ -x_1^t +x_1^\ast \\ 0  \\  -x_3^t +x_3^\ast \\0\\ \vdots\\ 0\\-x_{\tfrac{n+1}{2} -i}^t + x_{\tfrac{n+1}{2}-i}^\ast
\end{bmatrix}^t\begin{bmatrix} 0 \\ \vdots \\ 0 \\ 1 \\ 0 \\-x_2^t + x_2^{\ast} \\ \vdots \\ 0\\ -x_{i-j+1}^t + x_{i-j+1}^{\ast} \\ \vdots\\ 0\\-x_{\tfrac{n+1}{2} -j}^t + x_{\tfrac{n+1}{2}-j}^\ast \\0 \end{bmatrix}\begin{matrix}
\\ \\ \longleftarrow \\ \\ \\  \\ \\ \longleftarrow \\ \\ \\ \\
\end{matrix}\begin{matrix}
\\ (j-1)^{\text{th}} \text{ place}  \\ \\ \\ \\ \\  i^{\text{th}} \text{ place}  \\  \\ \\ \\
\end{matrix}\]}

$$= x_{i-j+1}^t + \sum_{l=1}^{\tfrac{i-j-1}{2}} x^t_{i-j-2l+1} (-x^t_{2l} + x^\ast_{2l})+  (-x^t_{i-j+1} + x^\ast_{i-j+1}) $$
$$= 0, \text{ using the definition of } x^\ast_{2n}.$$

\item[Case 3.] Consider $i$ even and $j$ even such that $i \geq j$. Then $C_i^t \widetilde{C}_j =$
    {\tiny \[ \begin{matrix}
\\ \\ \\ (j-1)^{\text{th}} \text{ place} \\ \\ \\  i^{\text{th}} \text{ place}  \\ \\  \\ \\ \\ \\ \\ \\ \\ \\
\end{matrix}\begin{matrix}
 \\ \\ \\  \longrightarrow \\ \\  \longrightarrow \\ \\ \\ \\ \\ \\ \\ \\ \\
\end{matrix}\begin{bmatrix}
     x_{i-1}  \\ \vdots\\ x_{i-j+2} \\ x_{i-j+1} \\ \vdots \\ 1 \\ -x_1^t +x_1^\ast \\ 0  \\  -x_3^t +x_3^\ast \\0\\ \vdots\\ 0\\-x_{\tfrac{n+1}{2} -i}^t + x_{\tfrac{n+1}{2}-i}^\ast
\end{bmatrix}^t\begin{bmatrix} x_{j-2} \\ \vdots \\ x_1 \\ 1 \\ 0 \\-x_1^t + x_1^{\ast} \\ \vdots\\ -x_{i-j+1}^t + x_{i-j+1}^{\ast} \\ 0\\\vdots\\ 0\\-x_{\tfrac{n+1}{2} -j}^t + x_{\tfrac{n+1}{2}-j}^\ast \\0 \end{bmatrix}\begin{matrix}
\\  \\  \longleftarrow \\ \\ \\  \\ \\ \longleftarrow \\ \\ \\ \\ \\
\end{matrix}\begin{matrix}
\\ \\ (j-1)^{\text{th}} \text{ place}  \\ \\ \\ \\  i^{\text{th}} \text{ place}  \\  \\ \\ \\
\end{matrix}\]}
$$= \sum_{l-k = i-j+1} x_l^t x_k + x^t_{i-j+1} + \sum_{\substack{l+k = i-j+1, \\ k \text{ odd}}} x^t_l (-x_k^t+x_k^\ast)+ (-x^t_{i-j+1}+ x^\ast_{i-j+1}) $$
$$= 0, \text{ using the definition of } x^\ast_{2n+1}.$$
\end{itemize}
\end{proof}
\begin{proposition}
    $N$ is a normal subgroup of $Sol_{A^2_n}.$
\end{proposition}
\begin{proof}
    First we prove that it is closed under multiplication. That is, for $X, Y \in N$ we need to show $XY \in N.$\\
Clearly, the rows of any matrix of $N$ are as follows:
$$R_1 = \big[a, x_1, 0, x_3, 0, \cdots, x_{\tfrac{n-3}{2}} \big],$$
$$R_i = \begin{cases}
    \big[0, -x_{i-2}^t + x_{i-2}^\ast, \cdots, -x_2^t + x_2^\ast, -x_1^t + x_1^\ast, 1, x_1, 0, \cdots, x_{\tfrac{n-1}{2} - i}, 0\big], & \text{ when } i>1 \text{ is odd}\\
    \big[0, \cdots, 0, 1, 0, x_2, 0, \cdots, x_{\tfrac{n-1}{2} - i}, 0\big] & \text{ when } i \text{ is even}
\end{cases}.$$ The columns will be referred from the equation \ref{columns_onemodfour}. Then one can easily verify $R_1C_1 = a; \ R_iC_i = 1; \ R_iC_j = 0$ for $i<j$ whenever $i,j$-odd or $i$-even, $j$-odd. Furthermore, whenever $i, j$-even,  then $$R_iC_j = x_{j-i} +y_{j-i} + \sum_{\substack{a+b=j-i \\ a \text{ even}}} x_ay_b$$
and when $i$-odd, $j$-even, then
\begin{equation*}
\begin{split}R_iC_j &= x_{j-i} +y_{j-i} + \sum_{\substack{a+b=j-i \\ a \text{ even}}} x_ay_b + \sum_{b-a=j-i} { (-x_a^t +x_a^*)y_b + x_a(-y_b^t +y_b^*)}\\
&= x_{j-i} +y_{j-i} + \sum_{\substack{a+b=j-i \\ a \text{ even}}} x_ay_b \\ 
& (\text{because } x_a^ty_b = 0, \text{ and } x_a^ty_b^\ast =0 \text{ when } a \not\equiv b (\text{mod } 2) \text{ from the Lemma } \ref{star_new})\end{split}\end{equation*}
Similarly, for $i>j$  the only non-zero case will be $i$-odd, $j$-even or $i,j$-odd and in that case
$$R_iC_j = -x_{i-j}^t+x_{i-j}^* -y_{i-j}^t+ y_{j-i}^* + \sum_{\substack{a+b=i-j \\ a \text{ even}}} (- x_a^t + x_a^*)(-y_b^t + y_b^*).$$
To show $XY \in N$, we need to show $-(R_1C_j)^t + (R_1C_j)^\ast = R_nC_{n-j+1}.$ After substituting the values it will boil down to show the following:
$$z_k^* = x_k^* + y_k^* + \sum_{\substack{a+b=k \\ a \text{ even}}} y_a^tx_b^t + \sum_{\substack{a+b=k \\ a \text{ even}}} (- x_a^t + x_a^*)(-y_b^t + y_b^*),$$
which follows from the Theorem \ref{star_one}. Now for any $X \in N,$ define the inverse $X^{-1}$ with the following columns
$$ C_1 = {\tiny\begin{bmatrix}
     a\\
     0\\
     \vdots\\
     0
\end{bmatrix}}; \ \text{For $i$ even}, \ \ C_i = {\tiny \begin{bmatrix}
     -x_{i-1}+x_{i-1}^{*t}  \\
     \vdots\\
     -x_1+x_1^{*t}\\
     1\\
     x_1^t\\
     0\\
     x_3^t\\
     \vdots\\
     0\\
     x_{\tfrac{n+1}{2} -i}^t
\end{bmatrix}}; \ \ \text{For $i$ odd, $i \geq 3$}, \ \ C_i = {\tiny\begin{bmatrix}
     0\\
     \vdots\\
     0\\
     1\\
     0\\
     x_2^t\\
     \vdots\\
     0\\
     x_{\tfrac{n+1}{2} -i}^t
\end{bmatrix}}.$$

One can check $XX^{-1} = I$. Furthermore, $N$ is normal because $\text{Lie}(N)$ is an ideal.
\end{proof}

\subsubsection{The case when $n$ is congruent to 3 mod 4}\

\textbf{Construction of $N$.} Define $a = {\footnotesize\begin{pmatrix}
0 &0\\
0&1
\end{pmatrix}}, \ b = {\footnotesize\begin{pmatrix}
0 &\lambda\\
0&0
\end{pmatrix}}$, $x_i =  {\footnotesize\begin{pmatrix}
0 &0\\
\ast&\ast
\end{pmatrix}};  \text{ if $i$ is odd}$, $x_2 = {\footnotesize\begin{pmatrix}
\ast &\ast\\
-\lambda&0
\end{pmatrix}}$, and
$ x_i = {\footnotesize\begin{pmatrix}
\ast &\ast\\
0&0
\end{pmatrix}};  \text{ if $i>2$ is even}$. Then for $2 \leq i \leq \frac{n+1}{2}$, $C_i$ denotes the $i^\text{th}$ column of $X$, in which 1 (the $2 \times 2$ identity matrix) is at the $i^\text{th}$-place and is defined as follows
\begin{equation} \label{3 mod4 columns}
C_1 = {\footnotesize\begin{bmatrix}
     a\\
     b\\
     0\\
     \vdots\\
     0
\end{bmatrix}}; \ \text{For $i$ odd, $i \geq 3$}, \ \ C_i = {\footnotesize\begin{bmatrix}
     x_{i-1}  \\
     \vdots\\
     x_1\\
     1\\
     -x_1^t + x_1^\ast\\
     0\\
     -x_3^t + x_3^\ast\\
     \vdots\\
     0\\
     -x_{\tfrac{n+1}{2} -i}^t + x_{\tfrac{n+1}{2}-i}^\ast
\end{bmatrix}}; \ \ \text{For $i$ even,}, \ \ C_i = {\footnotesize\begin{bmatrix}
     0\\
     \vdots\\
     0\\
     1\\
     0\\
     -x_2^t + x_2^\ast\\
     \vdots\\
     0\\
     -x_{\tfrac{n+1}{2} -i}^t + x_{\tfrac{n+1}{2}-i}^\ast
\end{bmatrix}}.\end{equation} Now define $$ x_{2n}^\ast = -\sum_{\substack{r+s = 2n \\
s \text{ even}}} x_r^t (-x_s^t+x_s^\ast),$$
where $x^\ast_2 = 0;$ and
$$ x_{2n+1}^\ast = b \ x_{2n+2}-\sum_{\substack{l+k = 2n+1 \\ k \text{ odd}}} x_l^t(-x_k^t + x_k^\ast).$$
Before constructing the set $N$, we observe the following: Given any variables $x_i \in M_2(\mathbb{C})$
    \begin{itemize}
        \item For $r, s \neq 1$ and $r \not\equiv s (\text{mod} \ 2)$, $x_r x_s=0$.
        \item $b x_{2k+1} = b^tx_{2k} = 0$.
        \item $-x_{2k}^tx_1 = bx_{2k} = x_{2k}^tb^t$.
        \item $x^t_{2k+1}a = 0$.
        \item $b^tb = \lambda^2 a$.
    \end{itemize}
    Using these we prove the following lemma:
    \begin{lemma}\label{relations}
        For $l, k > 1$, $x^*_{2k+1}x_{2l} = x_{2k}^*x_{2l+1} = x_{2l}x^*_{2k+1} =  x_{2k}^*x_1 = 0.$
    \end{lemma}
    \begin{proof}
        We prove these identities by induction. Observe that, $$ x_{2n+1}^\ast = b \ x_{2n+2}-\sum_{\substack{l+k = 2n+1 \\ k \text{ odd}}} x_l^t(-x_k^t + x_k^\ast).$$
        Now, assume by induction, we have $x_{2r+1}^* x_{2s} =0$ for all $2r=1 < 2n+1$. So $x_{2n+1}^\ast x_{2s} =  b \ x_{2n+2} x_{2s}-\sum_{\substack{l+k = 2n+1 \\ k \text{ odd}}} x_l^t(-x_k^t + x_k^\ast) x_{2s}$. This is equal to $b \ x_{2n+2} x_{2s}-\sum_{\substack{l+k = 2n+1 \\ k \text{ odd}}} x_l^t(-x_k^t + x_k^\ast) x_{2s}$. Now, using the induction hypothesis and the fact that $x_sx_t$ where $t \neq t$ modulo 2 is zero, we obtain, $x_{2n+1}^* x_{2s}=b \ x_{2n+2} x_{2s}$. Now $b x_{2n+2}= x_{2n+2}^tb^t$, but $b^t x_{2s}=0$ so we obtain the result. $x_{2k}^*x_{2l+1}=0$ is proved by a similar induction argument and observation that $x_{2k}^*$ only involves variables with even subscripts. The next two observations also follow similarly. 
    \end{proof}
Now let us define $N := \{ [C_1, C_2, \ldots, C_{\tfrac{n+1}{2}}] \ | \ x_i \in M_2(\mathbb{C})\}$.\\
\begin{lemma}
    $N$ is a subset of $Sol_{A^2_n}$. 
    \end{lemma}
    \begin{proof}
    For any $g \in N$, we want $g^tA^2_ng = A^2_n,$ which is same as checking $C_i^t\widetilde{C}_j = 1$ for $j = i+1$, and otherwise zero.
$$C_1^t\widetilde{C}_2 = {\tiny \begin{bmatrix}
    a^t& b^t &0 &\cdots &0
\end{bmatrix}\begin{bmatrix}
    1 \\ 0\\ -x_2+x_2^* \\ \vdots \\ -x_{\tfrac{n-1}{2}}^t + x_{\tfrac{n-1}{2}}^*\\0
\end{bmatrix}} = a^t = a.$$
For $i>3$ and $i$ is odd,
$$C_i^t\widetilde{C}_{i+1} = {\tiny \begin{matrix} &&i^{\text{th}}& \text{place}&&&&&\\
&&&\big\downarrow &&&&&\\
     \Big[x_{i-1}^t &\cdots &x^t_1 & 1 &-x_1 + x_1^{\ast t} &0& \cdots &0&-x_{\tfrac{n+1}{2}-i} + x_{\tfrac{n+1}{2}-i}^{\ast t}\Big]\\
     &&&&&&&&
\end{matrix}\begin{bmatrix} 0 \\ \vdots \\ 0 \\ 1 \\ 0 \\-x_2^t + x_2^{\ast} \\ \vdots \\0\\ -x_{\tfrac{n+1}{2}-i}^t + x_{\tfrac{n+1}{2}-i}^{\ast}\\0 \end{bmatrix}\begin{matrix}
\\ \\  \\  \\ \longleftarrow \\ \\ \\ \\ \\ \\ \\ \\
\end{matrix}\begin{matrix}
\\ \\ i^{\text{th}} \text{ place} \\  \\ \\ \\ \\ \\  \\
\end{matrix}}$$
$$= 1,$$

and for $i$ even, $$C_i^t \widetilde{C}_{i+1} = {\tiny\begin{matrix}
&&i^{\text{th}}& \text{place}&&&&&\\
&&&\big\downarrow &&&&&\\
     \Big[0 &\cdots &0 &1 &0&-x_2 + x_2^{\ast t} & \cdots &0&-x_{\tfrac{n+1}{2}-i} + x_{\tfrac{n+1}{2}-i}^{\ast t}\Big] \\
     &&&&&&&&
\end{matrix}\begin{bmatrix} x_{i-2} \\ \vdots \\ x_1\\ 1  \\-x_1^t + x_1^{\ast} \\ \vdots \\0\\ -x_{\tfrac{n+1}{2}-i}^t + x_{\tfrac{n+1}{2}-i}^{\ast}\\0 \end{bmatrix}
\begin{matrix}
\\ \\  \\  \\   \longleftarrow  \\ \\ \\ \\ \\ \\ \\
\end{matrix}\begin{matrix}
\\ \\ \\ i^{\text{th}} \text{ place} \\ \\  \\ \\ \\  \\
\end{matrix}}$$
$$ = 1.$$
Next we need to check $C_i^t \widetilde{C}_j$ is zero whenever $j \neq i+1.$ \\
\textbf{Case 1.} $i$ is odd, $j$ is even such that $i<j$, then it is clearly zero. For $i>j$, we have
$$C_i^t\widetilde{C}_j = {\tiny \begin{matrix}
(j-1)^{\text{th}} \text{ place} \\ \\ \\ \\ i^{\text{th}} \text{ place}  \\ \\  \\ \\ \\ \\ \\ 
\end{matrix}\begin{matrix}
 \longrightarrow \\ \\ \\ \longrightarrow \\ \\ \\ \\ \\ \\ \\
\end{matrix}\begin{bmatrix}
     x_{i-1}  \\ \vdots\\ x_{i-j+1} \\ \vdots\\x_1 \\ 1 \\ -x_1^t +x_1^\ast \\ 0  \\  -x_3^t +x_3^\ast \\0\\ \vdots\\ 0\\-x_{\tfrac{n+1}{2} -i}^t + x_{\tfrac{n+1}{2}-i}^\ast
\end{bmatrix}^t\begin{bmatrix} 0 \\ \vdots \\ 0 \\ 1 \\ 0 \\-x_2^t + x_2^{\ast} \\ \vdots \\ 0\\ -x_{i-j+1}^t + x_{i-j+1}^{\ast} \\ \vdots\\ 0\\-x_{\tfrac{n+1}{2} -j}^t + x_{\tfrac{n+1}{2}-j}^\ast \\0 \end{bmatrix}\begin{matrix}
\\ \\ \longleftarrow \\ \\ \\  \\ \\ \longleftarrow \\ \\ \\ \\
\end{matrix}\begin{matrix}
\\ (j-1)^{\text{th}} \text{ place}  \\ \\ \\ \\ \\  i^{\text{th}} \text{ place}  \\  \\ \\ \\
\end{matrix}}$$
$$ = x^t_{i-j+1} +\sum_{\substack{l+k = i-j+1\\ l \text{ even}}} x_k^t(-x_l^t + x_l^*) - x^t_{i-j+1} + x^*_{i-j+1}$$
$$ = 0 \ (\text{ by definition of } x_{2N}^*)$$
\textbf{Case 2.} $i$ is even, $j$ is odd such that $i>j$ it is clearly zero. But for $i<j$, we have
$$C_i^t\widetilde{C}_j = {\tiny \begin{matrix}
\\ \\ \\ \\ i^{\text{th}} \text{ place} \\ \\ \\ \\ \\ (j-1)^{\text{th}} \text{ place}  \\ \\  \\ \\ \\ \\ 
\end{matrix}\begin{matrix}
 \\ \\ \\ \longrightarrow \\ \\ \\ \\ \\  \longrightarrow \\ \\ \\ \\ \\ 
\end{matrix}\begin{bmatrix}
     0  \\ \vdots\\0 \\ 1 \\ 0\\-x_2^t +x_2^\ast \\ 0  \\  \vdots\\-x_{j-i+1}^t + x_{j-i+1}^*\\ \vdots\\ 0\\-x_{\tfrac{n+1}{2} -i}^t + x_{\tfrac{n+1}{2}-i}^\ast
\end{bmatrix}^t\begin{bmatrix} x_{j-2} \\ \vdots \\x_{j-i+1}\\\vdots\\ x_1 \\ 1 \\ 0 \\-x_1^t + x_1^{\ast}  \\ 0\\ \vdots\\ 0\\-x_{\tfrac{n+1}{2} -j}^t + x_{\tfrac{n+1}{2}-j}^\ast \\0 \end{bmatrix}\begin{matrix}
 \\  \longleftarrow \\ \\ \\   \longleftarrow \\ \\ \\ \\ \\ \\ \\
\end{matrix}\begin{matrix}
\\  i^{\text{th}} \text{ place}  \\ \\ \\   (j-1)^{\text{th}} \text{ place}  \\  \\ \\ \\ \\ \\
\end{matrix}}$$
$$ = x_{j-i+1} +\sum_{\substack{l+k = j-i+1\\ l \text{ even}}} (-x_l + x_l^{*t})x_k - x_{j-i+1} + x^{*t}_{j-i+1}$$
$$ = 0 \ (\text{ by definition of } x_{2n}^{*t})$$
\textbf{Case 3.} When both $i$ and $j$ are even, then it is easy to see $C^t_i\widetilde{C}_j = 0.$\\
\textbf{Case 4.} When both $i$ and $j$ are odd. If $i>j$, then we have
$$C_i^t\widetilde{C}_j =  {\tiny \begin{matrix}
(j-1)^{\text{th}} \text{ place} \\ \\ \\ \\ i^{\text{th}} \text{ place}  \\ \\  \\ \\ \\ \\ \\ 
\end{matrix}\begin{matrix}
 \longrightarrow \\ \\ \\ \longrightarrow \\ \\ \\ \\ \\ \\ \\
\end{matrix}\begin{bmatrix}
     x_{i-1}  \\ \vdots\\ x_{i-j+1} \\ \vdots\\x_1 \\ 1 \\ -x_1^t +x_1^\ast \\ 0  \\  -x_3^t +x_3^\ast \\0\\ \vdots\\ 0\\-x_{\tfrac{n+1}{2} -i}^t + x_{\tfrac{n+1}{2}-i}^\ast
\end{bmatrix}^t\begin{bmatrix} x_{j-2} \\ \vdots \\ x_1 \\ 1 \\-x_1^t + x_1^{\ast} \\0\\ \vdots \\-x_{i-j+1}^t + x_{i-j+1}^{\ast} \\ \vdots\\ 0\\-x_{\tfrac{n+1}{2} -j}^t + x_{\tfrac{n+1}{2}-j}^\ast \\0 \end{bmatrix}\begin{matrix}
\\ \\ \\ \longleftarrow \\ \\ \\  \\  \\ \longleftarrow \\ \\ \\ \\ \\
\end{matrix}\begin{matrix}
\\ (j-1)^{\text{th}} \text{ place}  \\ \\ \\ \\ \\  i^{\text{th}} \text{ place}  \\  \\ \\ \\
\end{matrix}}$$
$$= \sum_{l-k=i-j+1}x_l^tx_k + x_{i-j+1}^t +\sum_{\substack{l+k = i-j+1 \\ k \text{ odd}}}x_l^t(-x_k^t + x_k^*) - x_{i-j+1}^t + x_{i-j+1}^*$$
$$ = x_{i-j+2}^tx_1 + \sum_{\substack{l+k = i-j+1 \\ k \text{ odd}}}x_l^t(-x_k^t + x_k^*) + x_{i-j+1}^* $$
$$= 0 \ \ (\text{using definition of } x^*_{2k+1}).$$
\textbf{Case 5.} When both $i$ and $j$ are odd. If $i<j$, then we have
$$C_i^t\widetilde{C}_j =  {\tiny\begin{matrix}
\\ \\ \\  (j-1)^{\text{th}} \text{ place}  \\ \\ \\ \\ i^{\text{th}} \text{ place}  \\ \\  \\ \\ \\
\end{matrix}\begin{matrix}
\\ \\ \\ \\  \longrightarrow  \\ \\ \\ \\ \\  \longrightarrow \\ \\ \\ \\ \\
\end{matrix}\begin{bmatrix}
     x_{i-1} \\ \vdots\\x_1 \\ 1 \\ -x_1^t +x_1^\ast \\ 0  \\\vdots\\  -x_{j-i+1}^t +x_{j-i+1}^\ast \\\vdots\\0\\-x_{\tfrac{n+1}{2} -i}^t + x_{\tfrac{n+1}{2}-i}^\ast
\end{bmatrix}^t\begin{bmatrix} x_{j-2} \\ \vdots \\x_{j-i+1}\\\vdots\\ x_1 \\ 1 \\ 0 \\-x_1^t + x_1^{\ast}  \\ 0\\ \vdots\\ 0\\-x_{\tfrac{n+1}{2} -j}^t + x_{\tfrac{n+1}{2}-j}^\ast \\0 \end{bmatrix}\begin{matrix}
 \\  \longleftarrow \\ \\ \\   \longleftarrow \\ \\ \\ \\ \\ \\ \\
\end{matrix}\begin{matrix}
\\  i^{\text{th}} \text{ place}  \\ \\ \\   (j-1)^{\text{th}} \text{ place}  \\  \\ \\ \\ \\ \\
\end{matrix}}$$
$$= \sum_{k-l=j-i+1}x_l^tx_k + x_{j-i+1} +\sum_{\substack{l+k = j-i+1 \\ k \text{ odd}}}(-x_k + x_k^{*t}) x_l- x_{j-i+1} + x_{j-i+1}^{*t}$$
$$ = x_1^tx_{j-i+2} + \sum_{\substack{l+k = j-i+1 \\ k \text{ odd}}}(-x_k + x_k^{*t}) x_l + x_{j-i+1}^{*t} $$
$$= 0 \ \ (\text{using the definition of } x^{*t}_{2k+1}).$$
\end{proof}
\begin{proposition}
    $N$ is a normal subgroup of $Sol_{A^2_n}.$
\end{proposition}
\begin{proof}
    For this, first, we prove for any $X$ and $Y$ in $N$, we have $XY \in N$. The columns of the $Y$ matrix is given by equation \ref{3 mod4 columns}, and the rows of the matrix $X$ are as follows:
$$R_1 = \big[a, 0, x_2, 0, x_4, 0, \cdots, x_{\tfrac{n-3}{2}} \big],$$
$$R_2 = \big[b, 1, x_1, 0, x_3, 0, \cdots, x_{\tfrac{n-5}{2}} \big],$$
$$R_i = \begin{cases}
    \big[0, -x_{i-2}^t + x_{i-2}^\ast, \cdots, -x_2^t + x_2^\ast, -x_1^t + x_1^\ast, 1, x_1, 0, \cdots, x_{\tfrac{n-1}{2} - i}, 0\big], & \text{ when } i>2 \text{ is even}\\
    \big[0, \cdots, 0, 1, 0, x_2, 0, \cdots, x_{\tfrac{n-1}{2} - i}, 0\big], & \text{ when } i>1 \text{ is odd}
\end{cases}.$$
Then
$$R_1C_j = \begin{cases}
    a^2= a & \text{ if } j=1,\\
    0 & \text{ if } j \text{ is even,}\\
    x_{j-1} + y_{j-1}+ \sum_{\substack{k+l = j-1 \\ k-\text{even}}} x_ky_l :=  z_{j-1} & \text{ if } j \geq 3 \text{ is odd}.
\end{cases}$$
$$R_2C_j = \begin{cases}
    ab+a=b & \text{ if } j=1,\\
    1 & \text{ if } j=2,\\
    0 & \text{ if } j \geq 4 \text{ is even,}\\
    by_{j-1} +x_{j-2} + y_{j-2}+ \sum_{\substack{k+l = j-2 \\ k-\text{odd}}} x_ky_l :=  z_{j-2} & \text{ if } j \geq 3 \text{ is odd}.
\end{cases}.$$
To get $XY \in N$, the last row of the product matrix $XY$ should behave in the following way:
$$R_{\tfrac{n+1}{2}}C_{2i} = -(R_1C_{\tfrac{n+3}{2}-2i})^t + (R_1C_{\tfrac{n+3}{2}-2i})^*$$
That is, for even $m$
$$\sum_{\substack{l+k=m \\ l,k-\text{even}}}(-x_l^t + x_l^*)(-y_k^t + y_k^*) - x_m^t + x_m^* -y_m^t + y_m^* = -z_m^t + z_m^*$$
or
\begin{equation}
  z_m^* =    x_m^* + y_m^* +\sum_{\substack{k+l = m \\ k-\text{even}}} y_l^tx_k^t + \sum_{\substack{l+k=m \\ l,k-\text{even}}}(-x_l^t + x_l^*)(-y_k^t + y_k^*),
\end{equation}
which holds by Theorem \ref{star_one}. Whereas, for the odd entries, we get
$$R_{\tfrac{n+1}{2}}C_{2i+1} = -\left(R_2C_{\tfrac{n+3}{2}-2i}\right)^t + \left(R_2C_{\tfrac{n+3}{2}-2i}\right)^*.$$
That is, for odd $m$
$$\sum_{\substack{l+k=m \\ k-\text{odd}}}(-x_l^t + x_l^*)(-y_k^t + y_k^*) - x_m^t + x_m^* -y_m^t + y_m^* +\sum_{\substack{k-l=m \\ l<\tfrac{n-1}{2}-m}}(-x_k^t + x_k^*)y_l= -z_m^t + z_m^*$$
or
\begin{equation}
  z_m^* =    y_{m+1}^tb^t + x_m^* + y_m^*  +  \sum_{\substack{k+l = m \\ l-\text{odd}}} y_k^tx_l^t+ \sum_{\substack{l+k=m \\ k-\text{odd}}}(-x_l^t + x_l^*)(-y_k^t + y_k^*) + \sum_{\substack{k-l=m \\ l<\tfrac{n-1}{2}-m}}(-x_k^t + x_k^*)y_l.
\end{equation}

    Using the above lemma, $z_m^*$ for $m$ odd becomes,
    $$z_m^* =   b x_{m+1} + by_{m+1} + x_m^* + y_m^*  +  \sum_{\substack{k+l = m \\ l-\text{odd}}} y_k^tx_l^t+ \sum_{\substack{l+k=m \\ k-\text{odd}}}(-x_l^t + x_l^*)(-y_k^t + y_k^*) + x^*_{n+1}y_1$$
    $$= bz_{m+1} + x_m^* + y_m^*  +  \sum_{\substack{k+l = m \\ l-\text{odd}}} y_k^tx_l^t+ \sum_{\substack{l+k=m \\ k-\text{odd}}}(-x_l^t + x_l^*)(-y_k^t + y_k^*)$$
    $$= bz_{m+1} -\sum_{\substack{l+k = m \\ k \text{ odd}}} z_l^t(-z_k^t + z_k^\ast),$$
    the last equality follows from Theorem \ref{star_one}. Furthermore, for any $X \in N$, define the columns of the inverse matrix $X^{-1}$ as follows:
    $$
C_1 = {\tiny\begin{bmatrix}
     a\\
     b\\
     0\\
     \vdots\\
     0
\end{bmatrix}}; \ \text{For $i$ odd, $i \geq 3$}, \ \ C_i = {\tiny\begin{bmatrix}
    -x_{i-1}+x_{i-1}^{*t}  \\
     \vdots\\
     -x_1+x_1^{*t}\\
     1\\
     x_1^t\\
     0\\
     x_3^t\\
     \vdots\\
     0\\
     x_{\tfrac{n+1}{2} -i}^t
\end{bmatrix}}; \ \ \text{For $i$ even,}, \ \ C_i = {\tiny\begin{bmatrix}
     0\\
     \vdots\\
     0\\
     1\\
     0\\
     x_2^t\\
     \vdots\\
     0\\
     x_{\tfrac{n+1}{2} -i}^t
\end{bmatrix}}.$$
Clearly, $XX^{-1} = I.$ Also $N$ is normal as $\text{Lie}(N)$ is an ideal.
\end{proof}
    We now conclude the section with the main theorem, i.e., the complete structure of the group described as a semi-direct product of the diagonal part and the part coming from the subgroups $N$.
\begin{theorem}
    The group $Sol_{A^2_n}$ is a semi-direct product of $D$ and $N$.
\end{theorem}
\begin{proof}
Observe that since $N$ is a normal subgroup, we look at the map $\phi: D \to Aut(N)$, and it gives a semi-direct product structure of the group due to the short exact sequence.
\[ 1 \to N \to Sol_{A^2} \to D \to 1\]
\end{proof}

After proving $G$ as a semi-direct product of $D$ and $N$, we can clearly write the shape of matrix $g$ by multiplying $d \in D$ and $n \in N$. In $1 (\text{mod } 4),$ and $3 (\text{mod } 4)$ case we will ignore the first row and first column of $g$. Now if we assume that the field $k$ is $\mathbb{R}$ or $\mathbb{C}$, we can observe the following result.
\begin{lemma}
 The groups, $N,D,$ and $Sol_{A^2_n}$ are path-connected. The fundamental group of $Sol_{A^2_n}$ is isomorphic to that of $N$.
\end{lemma}

\begin{proof}
    The path connectedness of $D$ is an automatic consequence of the fact that it is isomorphic to $GL_2(k)$) (or a Borel subgroup of $GL_2(k)$) and hence contractible. For $N$, one observes that given a parameter $t \in [0,1]$ and an element $X \in N$, one defines the family of elements $X_t$, by changing the coordinates from $x_1,x_2, \ldots, x_n$ to $tx_1,tx_2, \ldots, tx_n$. Now observe that $X_t$ is an element in $N$ for all $t$ and for all $n$. This continuous map gives a path from the arbitrary element $X$ to the identity element. Once we have the path-connectedness for $D$ and $N$, the result for $Sol_{A^2_n}$ follows from the fact that it is a semidirect product of $D$ and $N$. Let us denote the group $Sol_{A^2_n}$ as $G$ and write down the long exact sequence of the fundamental groups we obtain,
    \[ \cdots\to \pi_2(D) \to \pi_1( N) \to \pi_1(G) \to \pi_1(D) \to 1 \]
    As $\pi_1(D),\pi_2(D)$ are trivial, we obtain the following result:
\end{proof}
\begin{theorem}\label{fundamental}
 The group $Sol_{A^2_n}$ is simply-connected.
 
\end{theorem}

\begin{proof}
    We will prove that the subgroup $N$ is simply connected, and thus, by the lemma above, we have the result. Let us construct a fibration from $N$ as follows. Let $N_1$ be the subgroup of $N$ consisting of the elements $X$, such that the ``local coordinate matrix" $x_1$ is zero. It is easily observed that this is a subgroup. We have the following fibration.
    \[ M_2(k) \to N \to N_1\]
 Observe that $M_2(k)$ is simply connected. Thus, from the long exact sequence, one obtains $\pi_1(N) \simeq \pi_1(N_1)$. Now from $N_1$, we construct another fibration with the ``matrix coordinate" $x_2$ equated to the zero matrix. Thus continuing we get the result.
\end{proof}

\section{Invariants of $Sol_{A}$}\label{Sec6}
In this section, we discuss the invariants of $Sol_{A}$, for a coefficient matrix $A$. The main result of this chapter is to show that the invariants of $Sol_{A_n}$ and $Sol_{A^2_n}$ are finitely generated. We will also comment on the solutions of the equation $X^tAX=B$. Note that $Sol_{A}$ acts on the vector space $M_n(k)$ by left multiplication; we will call this action and its restriction to $GL_n$ to be the natural action(s).  Let us say the orbits of the left action of $Sol_{A}$ on the vector space of matrices $M_n(k)$ are $\mathfrak{O}_{\alpha}$, where $\alpha$ belongs to an indexing set $\Lambda$. Let us pick a set of orbit representatives $X_{\alpha}$. Now if for some $\alpha$, $X_{\alpha}$ does not satisfy the equation $X^t AX=B$, the whole orbit $\mathfrak{O}_{\alpha}$ does not belong to the solutions $Sol(A,B)$. Let us say $\Gamma$ is the subset of the indices $\Lambda$ containing all $\alpha$ such that  $X_{\alpha}$ is a solution of the general congruence equation $X^tAX=B$. As a result, we can write down the solutions of the general congruence equation $X^tAX=B$ in the following manner,
\[Sol(A,B)=\bigcup\limits_{\alpha \in \Gamma}^{} \mathfrak{O}_{\alpha}.\]
The orbit space of actions of a Lie group is a historical field of study \cite{Procesi_history}. The case for a linearly reductive group, i.e., a group $G$, which has no non-trivial solvable radical and has the property that any linear representation is completely reducible, was solved by the efforts of several mathematicians starting with D. Hilbert in his famous paper that paved the way for commutative algebra \cite{Hilbert}. It is showed \cite{reductive} that for a linearly reductive group $G$ acting on the polynomials $k[\underline{x}]=k[x_1,x_2,x_3, \ldots , x_n]$, the algebra generated by the invariants $k[\underline{x}]^G$, is a finitely generated $k$ algebra. (Also see \cite{Procesi_matrix}). This gives a way to construct the space of orbits of the action of the group $G$ in geometric terms, which can be interpreted as the $G$ quotient of the affine variety $\mathbb{A}^n$ exists in the category of affine varieties (see \cite{GIT}). This case is aligned with our situation as we observed that we have to find the space of $G=Sol_A\times Sol_B$ orbits of $M_n(k)$, in other words, to find the quotient $M_n(k)/G$. When $G$ is not reductive, it poses some difficulty as the sub-algebra of invariants may not be finitely generated. One of the main approaches in this line is that of the compactification method, which is described for particular types of compactifications in \cite{Kirwan}, where the group $G$ is embedded in $\mathbb{C}^* \times G$. More generally, there are various other methods to deal with the non-reductive case as documented in the book  \cite{Kirwan_Book}, and in the survey article \cite{Kirwan_overview}. This choice of technique will be apparent when we have a structural understanding of the group or its Lie algebra. \\

Moreover, if we know the structure of $Sol_B$ as well, we can make use of the right action of $Sol_B$ on the matrices to write down the above decomposition in terms of the orbits $\mathfrak{O}_{\alpha}$ of the action of $Sol_A \times Sol_B$ as
\[ Sol(A,B)= \bigcup \limits_{\alpha \in \Gamma'}^{} \mathfrak{O}_{\alpha}.\]
So, it is essential to understand the individual orbits of both left and right actions of $Sol_{A}$ on the space of matrices, for any $A$.\\

\subsection{Finite generation of $Sol_{A_n}$ and $Sol_{A_n^2}$.} Let us recall a few terminologies and results from \cite{observable}. For a closed subgroup $H$ of an affine algebraic group $G$, consider the natural action of $H$ on the ring $k[G]$ of regular functions of $G$. Let us call $H'$ to be the set of invariants of this action, and further consider $H''$ to be the set of elements $g \in G$ such that $g$ leaves all the invariants in $H'$ invariant. Naturally, $H \subseteq H''$. If these two coincide, we call the subgroup $H$ to be an observable group. 

\begin{theorem}[\cite{observable}, Theorem 2, \textsection 3 ] Let $G$ be a connected afffine algebraic group. For $H$ an algebraic subgroup (closed) of $G$, the following conditions are equivalent:
\begin{itemize}
    \item $H=H''$
    \item there is a finite-dimensional representation $\rho: G \to GL(V)$, and a vector $v \in V$ such that,
    \[ H= S_G(\rho, v)= \{ g \in G \, \mid \, \rho(g)v=v\}.\]
\end{itemize}
    Furthermore, if $H$ satisfies either of the above conditions, then there is a finite-dimensional representation $\rho: G \to GL(V)$ and a vector $v \in V$ such that $H= S_G(\rho,v)$ and $G/H \simeq G\cdot v.$
\end{theorem}

\begin{theorem}[\cite{observable}, Theorem 1, \textsection 4]\label{non_reductive} Let $G$ be a connected affine algebraic group and $H$ be a closed connected subgroup of $G$. If the radical of $H$ is nilpotent, then $H$ is observable in $G$.
    
\end{theorem}

Consequently, finite generation of $k[G]^H$ as a $k$ algebra is implied by the observability of $H$ as a subgroup of $G$. Thus we have a following theorem:

\begin{theorem}\label{invariant_main}
    The algebra of invariants $k[GL_n]^{G}$, where $G$ action on $GL_n$ is the natural action by the inclusion $G \subseteq GL_n$, is finitely generated when $G=Sol_{A_n}$ for all $n$; and $Sol_{A^2_n}$ for $n$ even. 
\end{theorem}
\begin{proof}
   The main tool for us to prove this is the Theorem \ref{non_reductive}; we show that $G$ is an observable subgroup of $GL_n$ as follows. For $G=Sol_{A_n}$ from Proposition \ref{radical_A}, we know that the solvable radical of the group is nilpotent, and hence $G$ is observable. For $G=Sol_{A_n^2}$, in the case $n$ is divisible by $4$, the solvable radical of the Lie algebra is nilpotent by Corollary \ref{radical_0}. For $n$ is congruent to 2 modulo 4, the radical is nilpotent by Corollary \ref{radical_2}.
\end{proof}
The radicals of the groups $Sol_{A_n}^2$, when $n$ congruent to 1 or 3 modulo 4, are not nilpotent. As a result, one can not directly use the results of \cite{observable} to derive the finite generation of the $Sol_{A_n}^2$ action on $\text{GL}_n$. In fact, the groups in these cases are solvable. Invariant theory of solvable groups is a difficult and largely unknown area. But we can say that from the Remark \ref{Invariant_odd_remark}, the problem of finding the invariants for the odd cases reduces to finding that of the case when $n$ is congruent to 3 modulo 4. \\

Recall from \cite{observable}, if $H$ is an observable subgroup of $G$, then are not only the invariants of $H$ action on $G$ finitely generated but also the invariants of $H$ action on any affine variety finitely generated. As a result, the algebra of $Sol_{A_n}$, for any $n$, and  $Sol_{A_n^2}$, for $n$ even, invariants for the action of $G$ on $M_n(K)$ is also finitely generated.  \\

Even though the space of invariants is finitely generated, finding them can still be challenging (there are a few algorithms. See, for example, \cite{algo_non_reductive}). Thus, in the examples below, we are trying to find the orbits and realize in one example $Sol_{A_6}$ we can find all the invariants, whereas, in the example of $Sol_{A_8^2}, $ it is cumbersome. 
\begin{enumerate}
    \item Recall $Sol_{A_6}$ consists of matrices $X$ where  
$$ X = {\tiny\begin{pmatrix}
x_0&0&x_0x_1&0 & x_0x_2&0\\
0&x_0^{-1}&0&0&0&0\\
0& 0&x_0&0& x_0x_1 & 0\\
0 & -x_0^{-1}x_1 & 0 & x_0^{-1} & 0& 0\\
0&0 & 0& 0 & x_0 & 0 \\
0 & -x_0^{-1}(x_2 - x_1^{2}) & 0& -x_0^{-1}x_1 & 0& x_0^{-1}
\end{pmatrix}}.$$
Then the action of $Sol_{A_6}$ on $M_6(\mathbf{C})$ by left multiplication as $g \cdot Y = gY$. Under this action $M_6(\mathbf{C}) = \bigcup_{Y \in M_6(\mathbf{C})} Sol_{A_6}[Y]$, where $Sol_{A_6}[Y]=$ $Sol_{A_6}$-orbit of $Y = Sol_{A_6}/ Stab[Y]$. In this example, we obtain the stabilizers as follows:
\begin{itemize}
    \item if $R_1(Y)$ or $R_6(Y)$ is not a zero vector, then $Stab[Y] = N$.
    \item if $R_2(Y)$ or $R_5(Y)$ is not a zero vector, then $Stab[Y] = I_6$.
    \item if $R_3(Y)$ or $R_4(Y)$ is not a zero vector, then $Stab[Y]$ consists of the matrices $${\tiny \begin{pmatrix}
1&0&0&0 &x_2&0\\
0&1&0&0&0&0\\
0& 0&1&0& 0 & 0\\
0 & 0 & 0 & 1 & 0& 0\\
0&0 & 0& 0 & 1 & 0 \\
0 & -x_2 & 0& 0 & 0& 1
\end{pmatrix}}$$
\end{itemize}
\item For $n = 8$, let us consider the matrix $X = {\footnotesize\begin{pmatrix}
a&b&c&d\\
e&f&g&h\\
i&j&k&l\\
m&n&o&p
\end{pmatrix}} \in M_{8}(\mathbb{C}),$
where $a, b, \cdots, o, p \in M_2(\mathbb{C}).$ Now using the structure of the group $Sol_{A^2_8}$, any element of this group looks like $Y = {\footnotesize\begin{pmatrix}
s & 0& x&0\\
0&s^{-t} &0&0\\
0&0&s&0\\
0&-s^{-t}xs^{-t}&0&s^{-t}
\end{pmatrix}}$ for some $x \in M_2(\mathbb{C})$. From the discussion at the beginning of this section, we would like to understand the orbits of the group action. Now as a consequence of the orbit stabilizer Lemma, it is equivalent to understanding the stabilizers of individual points. The stabilizer $G_X = \{Y \ | \ YX = X\},$ that is, 
$${\footnotesize\begin{pmatrix}
s & 0& x&0\\
0&s^{-t} &0&0\\
0&0&s&0\\
0&-s^{-t}xs^{-t}&0&s^{-t}
\end{pmatrix} \begin{pmatrix}
a&b&c&d\\
e&f&g&h\\
i&j&k&l\\
m&n&o&p
\end{pmatrix}}  = {\footnotesize\begin{pmatrix}
a&b&c&d\\
e&f&g&h\\
i&j&k&l\\
m&n&o&p
\end{pmatrix}}.$$
This leads to the following matrix equations:
\begin{equation*}
   (1-s)a = xi; \ \ \ \ \ (1-s)b = xj; \ \ \ \ \ (1-s)c = xk; \ \ \ \ \ (1-s)d = xl
\end{equation*}
\begin{equation*}
   (s^{-t} -1)e; \ \ \ \ \ (s^{-t} -1)f =0; \ \ \ \ \ (s^{-t} -1)g =0; \ \ \ \ \ (s^{-t} -1)h =0
\end{equation*}
\begin{equation*}
   (s-1)i = 0; \ \ \ \ \ (s-1)j = 0; \ \ \ \ \ (s-1)k = 0; \ \ \ \ \ (s-1)l = 0
\end{equation*}
\begin{equation*}
   (1-s^{-t})m= ye; \ \ \ \ \ (1-s^{-t})n= yf; \ \ \ \ \ (1-s^{-t})o= yg; \ \ \ \ \ (1-s^{-t})p= yh,
\end{equation*}
where $y = -s^{-t}xs^{-t}.$

\begin{lemma}
    If any linear combination of $e,f,g,h$ or $i,j,k,l$ is invertible, then $G_X = (1).$
\end{lemma}
\begin{proof}
Let $\alpha e + \beta f + \gamma g + \delta h$ be invertible. From the equations, we get
$$(s^{-t} -1)(\alpha e + \beta f + \gamma g + \delta h) = 0,$$ which implies $s^t = 1.$ Then $y (\alpha e + \beta f + \gamma g + \delta h) = 0,$ gives $y=0$ but $y = -s^{-t}xs^{-t}.$ Therefore, $G_X = (1).$ 
\end{proof}
Once we have a trivial stabilizer, we can write down the orbit in this case. Let us discuss an example following this case, since we want some linear combination of $e,f,g,h$ or $i,j,k,l$, invertible, for a simple example, we have taken the case $e=I_2$, the identity matrix.
Let $g = {\footnotesize\begin{pmatrix}
0&0&0&0\\
I_2 & f&g&h\\
i&j&k&l\\
0&0&0&0
\end{pmatrix}}$, and let us take the case when $B = g^tA^2_8g$. Then $\ \alpha g \in \text{Sol}(A_8^2, B)$ for all $\alpha \in \text{Sol}_{A_8^2}$ or the following types of matrices are solutions of the general congruence equation.
$$\alpha g = {\footnotesize\begin{pmatrix}
s &0&x&0\\
0&s^{-t}&0&0\\
0&0&s&0\\
0& - s^{-t}xs^{-t}&0&s^{-t}
\end{pmatrix}\begin{pmatrix}
0&0&0&0\\
I_2 & f&g&h\\
i&j&k&l\\
0&0&0&0
\end{pmatrix}} $$
$$={\footnotesize\begin{pmatrix}
xi &xj&xk&xl\\
s^{-t}&s^{-t}f&s^{-t}g&s^{-t}h\\
si&sj&sk&sl\\
 - s^{-t}xs^{-t}&- s^{-t}xs^{-t}f&- s^{-t}xs^{-t}g&- s^{-t}xs^{-t}h.
\end{pmatrix}}.$$
The remaining case is when every linear combination of $e,f,g,h$ or $i,j,k,l$, is singular. Such a space is called a non-singular subspace of matrices. Let us make one small observation for singular subspaces in $M_2(\mathbb{C})$.
 \begin{lemma}
    For $x_1, x_2, \ldots, x_n \in M_2(k).$ If the $\text{span}\{x_1, x_2, \ldots, x_n\}$ is a singular subspace, then the span must contain vectors of the form of any one of the following four types:
    $${\footnotesize\begin{pmatrix}
    \ast & 0\\
    \ast & 0
    \end{pmatrix}}, \ \ \ \ {\footnotesize\begin{pmatrix}
    \ast & \ast\\
    0 & 0
    \end{pmatrix}}, \ \ \ \ {\footnotesize\begin{pmatrix}
    0 & \ast\\
    0 & \ast
    \end{pmatrix}}, \ \ \ \ {\footnotesize\begin{pmatrix}
    0 & 0\\
    \ast & \ast
    \end{pmatrix}}$$
    \end{lemma}
    If we assume the span of $\{e, f, g, h\}, \{i, j, k, l\}$ are singular subspaces. If any of these $\{e, f, g, h\}, \{i, j, k, l\}$ is not zero then $ 1 \in \sigma(S)$. If $e, f, g, h, i, j, k, l = 0.$ Then     
    $$(1-s)a = 0; \ \ \ \ \ (1-s)b = 0; \ \ \ \ \ (1-s)c = 0; \ \ \ \ \ (1-s)d = 0$$
    $$(1-s^{-t})m= 0; \ \ \ \ \ (1-s^{-t})n= 0; \ \ \ \ \ (1-s^{-t})o= 0; \ \ \ \ \ (1-s^{-t})p= 0.$$
    Any of $\text{span}\{a, b, c, d\}$ or $\text{span}\{m, n, o, p\}$ must not contain an invertible element else $G_X = (1)$. On the other hand, if $\text{span}\{a, b, c, d\}$ and $\text{span}\{m, n, o, p\}$ are singular subspaces, then $1 \in \sigma(S)$. Now, one deals with the remaining cases through the Jordan canonical form of $s$.
    
    Let the Jordon form of $s$ be
    $${\footnotesize\begin{pmatrix}
    1 & 0\\
    0 & \mu
    \end{pmatrix}} \ \ \mu \neq 0; \ \ {\footnotesize\begin{pmatrix}
    1 & 1\\
    0 & 1
    \end{pmatrix}}.$$

Without loss of generality $\text{span}\{a, b, c, d\}$ be ${\footnotesize\begin{pmatrix}
    \ast & 0\\
    \ast & 0
    \end{pmatrix}}$. Then 
    $${\footnotesize\left(I_2 - \begin{pmatrix}
    1 & 0\\
    0 & \mu
    \end{pmatrix}\right) \begin{pmatrix}
    \ast & 0\\
    \ast & 0
    \end{pmatrix}} = {\footnotesize\begin{pmatrix}
    0 & 0\\
    0 & 1-\mu
    \end{pmatrix} \begin{pmatrix}
    \ast & 0\\
    \ast & 0
    \end{pmatrix}}.$$
    This implies, if the $(2,1)$ entry is non-zero then $\mu =1$. It gives $s =1$, which further gives $x =0$, which again leads to a trivial stabilizer. Therefore we can assume the  $\text{span}\{a, b, c, d\} = {\footnotesize\begin{pmatrix}
    \ast & 0\\
    0 & 0
    \end{pmatrix}}$. Similarly, the other cases can be reduced to the following shapes:
    $${\footnotesize\begin{pmatrix}
    0 & \ast\\
    0 & \ast
    \end{pmatrix}} \rightarrow{\footnotesize \begin{pmatrix}
    0 & \ast\\
    0 & 0
    \end{pmatrix}}; \ \ \ \ {\footnotesize\begin{pmatrix}
    \ast & \ast\\
    0 & 0
    \end{pmatrix}} \to {\footnotesize\begin{pmatrix}
    \ast & \ast\\
    0 & 0
    \end{pmatrix}}; \ \ \ {\footnotesize\begin{pmatrix}
    0 & 0\\
    \ast & \ast
    \end{pmatrix}} \rightarrow {\footnotesize\begin{pmatrix}
    0 & 0\\
    0 & 0
    \end{pmatrix}}.$$
Similarly, for $m, n, o, p$, there will be multiple subcases. It is evident that this procedure will be quite extensive, even for relatively small values of $n$. For larger values of $n$ and other parities of $n (\text{mod } 4)$, the approach is far from straightforward. 
\end{enumerate}

\section{Appendix}\label{Sec7}
 \subsection{An algebraic identity}
 Let us consider the noncommutative $k$ algebra generated by $x_1,x_2,x_3 \ldots x_n, y_1,y_2,y_3, \ldots y_n, x_1^t, x_2^t, x_3^t \ldots x_n^t, y_1^t, y_2^t, y_3^t, \ldots y_n^t$  \[k \langle \underline{x}, \underline{y}\rangle=k\langle x_1,x_2,x_3 \ldots x_n, y_1,y_2,y_3, \ldots y_n, x_1^t, x_2^t, x_3^t \ldots x_n^t, y_1^t, y_2^t, y_3^t, \ldots y_n^t\rangle.\] Let $t$ be an involution operator that has the following properties:
 \begin{itemize}
     \item For any $f,g \in k\langle\underline{x}, \underline{y}\rangle$ we have $ t(fg)=t(g)t(f)$,
     \item $t(x_i)=x_i^t$ and $t(y_i)=y_i^t$ for all $i$.
 \end{itemize}
Let us denote an algebra of the above kind as an algebra with a transpose. In a non-commutative algebra with transpose as in the above, let us define a function $*: k\langle\underline{x}, \underline{y}\rangle \to k\langle\underline{x}, \underline{y}\rangle $ inductively as, 
\begin{equation}\label{def_of_star}
    x_l^\ast = \sum_{r+s = l} x^t_s x_r^t -  \sum_{\substack{r+s = l }}x^t_s x_r^\ast.
\end{equation}  Let us also define $L(x)=x^*-x^t$ for any element $x$ in $k\langle\underline{x}, \underline{y}\rangle$.

\begin{lemma}\label{rearrangement}
    Let $x_i,y_i,z_i$, $1 \leq i \leq n$ be any elements in a non-commutative algebra we have, 
    \[ \sum_{a+b=n} x_a \sum_{s+t=b} y_sz_t = \sum_{a+b=n} (\sum_{s+t=a} x_sy_t) z_b.\]
\end{lemma}
\begin{proof}
    This follows from opening the sum and re-arranging.
\end{proof}
 \begin{theorem}\label{star_one}
   Let $x_i, y_i$ belong to non-commutative algebra with transpose $k\langle\underline{x}, \underline{y}\rangle$, let us define $z_m$ inductively as $z_m= x_m+y_m+ \displaystyle\sum_{a+b=m}x_ay_b $, then we have the identity,
    
   \[\sum_{a+b=n} z_a^t L(z_b) = x_n^{*} + y_n^{*} + \sum_{a+b=n} y_a^t x_b^t + \sum_{a+b=n} L(x_a)L(x_b).\]
\end{theorem}
\begin{proof}
    Let us substitute the form of $Z_a^t$, and the form of $L(Z_b)$ that we know by induction. The left-hand side becomes,
    \[(x_a^t + y_a^t + \sum_{s+t=a} y_s^t x_t^t)(L(x_b)+L(y_b) - \sum_{u+v=b}L(x_u)L(y_v)).\]
    After opening the product, we obtain the following summed over $a+b=n$ we obtain the left-hand side as sum of the following four expressions summed over the same indices.
        \begin{equation*}
    \begin{split}
    I&= x_a^tL(x_b)+ y_a^tL(y_b)\\
    II &=x_a^tL(y_b) - x_a ^t\sum_{u+v=b} L(x_u)L(y_v)\\
    III &=- y_a^t\sum_{u+v=b} L(x_u)L(y_v) + \sum_{s+t=a} y_s^t x_t^t L(y_b) - \sum_{s+t=a} y_s^tx_t^t \sum_{u+v=b} L(x_u)L(x_v)\\
    IV &= y_a^tL(x_b) +\sum_{s+t=a} y_s^t x_t^t L(x_b)
    \end{split}
    \end{equation*}
    Now observe that the sum over the indices $a,b$ following $a+b=n$ of the expression $I$, yields, $x_n^*+y_n^*$. The sum of $II$ is by Lemma \ref{star_two} is equal to the sum $\sum_{a+b=n} L(x_a)L(y_b)$. The sum of $III$ is zero by the Lemma \ref{star_three}. The fourth term, $IV$ is equal to $\sum_{a+b=n} y_a^t x_b^t$ from the Lemma \ref{star_four}. This completes the proof of the claim.
\end{proof}
 \begin{lemma}\label{star_two}
    With the notations above, we have,
    \[ \sum_{a+b=n} L(x_a)L(y_b)= \sum_{a+b=n} x_a^t L(y_b) - \sum_{a+b=n} x_a^t \sum_{s+t=b} L(x_s)L(y_t).\]
\end{lemma}
\begin{proof}
    Let us look at the right hand side $\displaystyle \sum_{a+b=n} x_a^t \sum_{s+t=b} L(x_s)L(y_t)$, once we apply the re-arrangement Lemma \ref{rearrangement}, we obtain,
    $\displaystyle \sum_{a+b=n} \sum_{s+t=a} x_s^t L(x_t) L(Y_b)$. Now the sum $\sum_{s+t=a} x_s^t L(x_t)$, is $x_a^{*}$ by definition, now after we observe that $L(x_a)= x_a^t- x_a^*$, we obtain the result.
\end{proof}
\begin{lemma}\label{star_three}
    With the notations as in the above, when we sum the following over $a+b=n$, we obtain the following:
    \[  y_a^t \sum_{s+t=b} L(x_s)L(y_t) - \sum_{s+t=a} y_s^t x_t^t L(y_b) +  \sum_{s+t=a} y_s^t x_t ^t \sum_{u+v=b} L(x_u)L(y_v)=0,\]
    and
    \[ y_a^t x_b^* - \sum_{s+t=a} Y_s^tx_t^t L(x_b)=0. \]
\end{lemma}
\begin{proof}
    First, observe that, $ y_a^t x_b^*  - \sum_{s+t=a} Y_s^tx_t^t L(x_b)=0$, this follows from the re-arrangement Lemma \ref{rearrangement}. For the next,
    \[ y_a^t \sum_{s+t=b} L(x_s)L(y_t) - \sum_{s+t=a} y_s^t x_t^t L(y_b) +  \sum_{s+t=a} y_s^t x_t ^t \sum_{u+v=b} L(x_u)L(y_v)\]
    which, when summed over $a+b=n$ we must prove to be zero. Now let us start with,
    \[\sum_{a+b=n}\sum_{s+t=a} y_s^t x_t ^t \sum_{u+v=b} L(x_u)L(y_v),\] using the re-arrangement Lemma \ref{rearrangement}, one can rewrite this as,
    \[\sum_{a+b=n}\sum_{s+t=a} y_s^t x_t ^t \sum_{u+v=b} L(x_u)L(y_v)=\sum_{a+b=n} y_a^t \sum_{s+t=b} x_s^* L(y_t).\]
    We get the result when we substitute and observe that $L(x_s)= x_s^t -x_s^*$.
\end{proof}
\begin{lemma}\label{star_four}
    With the notations as in the above, we have the following equality:
    \[ \sum_{a+b=n} y_a^tL(x_b) +\sum_{s+t=a} y_s^t x_t^t L(x_b) = \sum_{a+b=n} y_a^tx_a^t.\]
\end{lemma}
\begin{proof}
    If we start with the second term, we obtain the result and use the re-arrangement Lemma \ref{rearrangement}.
\end{proof}

\textbf{Acknowledgment:}
 HM is partially supported by CDRF, BITS Pilani, India, grant no. CDRF, C1/23/185. GS is partially supported by DST-SERB POWER Grant No. SPG/2022/001738.

\bibliographystyle{plain}
\bibliography{linear}
\end{document}